\documentclass[11pt,a4paper]{article}
\usepackage[T1]{fontenc}
\usepackage{lmodern}
\usepackage[margin=27mm]{geometry}
\usepackage{amsmath,amssymb,amsthm,mathrsfs,mathtools,booktabs,xcolor,hyperref}
\usepackage{graphicx,float,bm}
\usepackage{tikz}
\usetikzlibrary{arrows.meta}

\hypersetup{hidelinks,pdftitle={Randomly shifted lattice rules for SDE simulation with error bounds uniform in the time discretization and dimension}} 

\newtheorem{theorem}{Theorem}[section]
\newtheorem{lemma}[theorem]{Lemma}
\newtheorem{proposition}[theorem]{Proposition}
\newtheorem{corollary}[theorem]{Corollary}
\theoremstyle{definition}
\newtheorem{assumption}[theorem]{Assumption}
\newtheorem{definition}[theorem]{Definition}
\theoremstyle{remark}
\newtheorem{remark}[theorem]{Remark}
\newcommand{\E}{\mathbb E}
\newcommand{\Prob}{\mathbb P}
\newcommand{\R}{\mathbb R}
\newcommand{\N}{\mathbb N}
\newcommand{\Z}{\mathbb Z}
\newcommand{\PP}{\mathcal P}
\newcommand{\norm}[1]{\left\|#1\right\|}
\newcommand{\dd}{\,\mathrm d}
\numberwithin{equation}{section}
\title{Randomly shifted lattice rules for SDEs { with constant diffusion}: error bounds uniform in time discretization and dimension}
\author{Du Ouyang\thanks{Centre de Math\'{e}matiques Appliqu\'{e}es, \'{E}cole Polytechnique, Palaiseau 91128, France (\texttt{du.ouyang@polytechnique.edu}).}
     \and Josef Dick\thanks{School of Mathematics and Statistics, The University of New South Wales, Sydney, NSW 2052,
Australia (\texttt{josef.dick@unsw.edu.au}).}
\and Ziyang Ye\thanks{Corresponding author. Department of Mathematical Sciences, Tsinghua University, Beijing 100084, People’s Republic
of China (\texttt{yzy23@mails.tsinghua.edu.cn}).}
}
\date{ }

\begin{document}
\maketitle
\begin{abstract}
We study randomly shifted rank-1 lattice rules for computing
$\E[g(\bm{X}_T)]$, where $\bm{X}_T$ is the terminal state of a
stochastic differential equation (SDE) with constant diffusion matrix.
Discretizing the SDE by the Euler--Maruyama scheme with $m$
uniform time steps turns the expectation into an integral whose
dimension grows with $m$ and whose integrand has a nested
composition structure; both complicate the lattice rule error analysis.
We generate the Brownian path by the Brownian bridge construction
and handle the nested structure using the multivariate
Fa\`a di Bruno formula and a rooted-tree expansion.
Under suitable regularity and growth conditions, we prove that an
$N$-point randomly shifted lattice rule with a uniformly random
generating vector achieves a mean absolute sampling error of
$O(N^{-1+\varepsilon})$ for every $\varepsilon\in(0,1/2)$ and every
prime $N$, with the implied constant independent of $m$.
To improve stability and control higher error moments, the
median lattice estimator is used.
The estimator takes the median of $R=2k-1$ independent lattice
estimates and achieves an $L^q$ error of $O(N^{-1+\varepsilon})$
for every $1\le q<k/(1-\varepsilon)$.
We further establish error bounds independent of both the
state dimension and the {Brownian path dimension} under an additional summability
condition on the columns of the diffusion matrix.
Numerical experiments show that median lattice estimators
converge faster than Monte Carlo estimators and remain stable
as the number of time steps $m$ grows.
For models satisfying the additional summability condition,
stability is also observed as the state and Brownian path
dimensions increase.
\end{abstract}



\section{Introduction}\label{sec:introduction}
Stochastic differential equations (SDEs) model systems subject to random
perturbations. They are widely used in biology, chemistry, mechanics,
and finance \cite{higham2001algorithmic}. Many applications require
computing expectations of functions of SDE solutions. For example,
such expectations are used to price financial derivatives
\cite{glasserman2004} and to solve parabolic partial differential
equations \cite{oksendal2003stochastic}.

In this paper, we study the numerical approximation of such
expectations, focusing on SDEs with constant diffusion.
Specifically, let
$\mathcal I=\E[g(\bm{X}_T)]$, where $\bm{X}$ solves
\begin{equation}\label{eq:sde}
 \dd \bm{X}_t=b(t,\bm{X}_t)\dd t+\Sigma\dd \bm{W}_t,\qquad \bm{X}_0=\bm{x}_0,
 \qquad t\in[0,T].
\end{equation}
Here $T>0$ is fixed, $\bm{X}_t\in\R^{d_x}$, and $\bm{W}$ is a standard
$d_w$-dimensional Brownian motion. The drift is
$b:[0,T]\times\R^{d_x}\to\R^{d_x}$, the diffusion matrix
$\Sigma\in\R^{d_x\times d_w}$ is constant, and
$g:\R^{d_x}\to\R$ is the payoff function. The initial state
$\bm{x}_0\in\R^{d_x}$ is deterministic.
SDEs of this form arise in various applications, including
overdamped Langevin dynamics \cite{pavliotis2014stochastic}.

Explicit solutions of SDEs are generally unavailable, so numerical
schemes are used to approximate their sample paths.
We use the Euler--Maruyama (EM) scheme \cite{kloeden1992numerical}
on the uniform grid $\{t_i=ih: i=0,\ldots,m\}$, with $h=T/m$,
and set $s=d_wm$.
For $\bm{z}=(\bm{z}_1,\ldots,\bm{z}_m)\in(\R^{d_w})^m\cong\R^s$, set
$\widetilde{\bm{X}}_{t_0}(\bm{z})=\bm{x}_0$ and define
\begin{equation}\label{eq:em}
	\widetilde{\bm{X}}_{t_{i+1}}(\bm{z})
	=\widetilde{\bm{X}}_{t_i}(\bm{z})
	+h b(t_i,\widetilde{\bm{X}}_{t_i}(\bm{z}))
	+\sqrt h\,\Sigma\bm{z}_{i+1},
	\qquad i=0,\ldots,m-1.
\end{equation}
Let $\bm{Z}\sim N(0,{\operatorname{Id}_s})$ be an $s$-dimensional standard Gaussian
vector. Then $\widetilde{\bm{X}}_{t_m}(\bm{Z})$ is the EM approximation
at time $T$. The corresponding approximation to $\mathcal I$
can be written as
\begin{equation}\label{eq:intro-integral}
	\mathcal I_m
	:=\E[g(\widetilde{\bm{X}}_{t_m}(\bm{Z}))]
	 =\int_{(0,1)^s}g\bigl(\widetilde{\bm{X}}_{t_m}(\Phi^{-1}(\bm{u}))\bigr)\dd \bm{u},
\end{equation}
where $\Phi$ is the standard Gaussian distribution function
and $\Phi^{-1}$ acts componentwise. Under suitable regularity and growth conditions, $\mathcal I_m$
approximates $\mathcal I$ with the classical weak error bound
$|\mathcal I_m-\mathcal I|\le C m^{-1}$, where $C$ is independent
of $m$; see \cite[Theorem~14.5.1]{kloeden1992numerical}.
Thus, increasing the number of time steps $m$ reduces the
discretization error.

To compute $\mathcal I_m$, however, we typically approximate the
integral in \eqref{eq:intro-integral} numerically, since it
generally has no closed-form expression.
In this paper, we focus on the error introduced by numerical
integration, referred to as the \textbf{sampling error},
as distinct from the time discretization error
$|\mathcal I_m-\mathcal I|$.
The integration dimension $s=d_wm$ grows as the time grid is
refined and becomes arbitrarily large as $m\to\infty$.
This growth poses a central challenge in obtaining sampling
error bounds that remain uniform in $m$.
Monte Carlo (MC) provides a straightforward approach.
For independent uniform points $\bm{u}_0,\ldots,\bm{u}_{N-1}\in(0,1)^s$,
the MC estimator is
\begin{equation}\label{eq:intro-quadrature}
	\frac1N\sum_{n=0}^{N-1}
	g\bigl(\widetilde{\bm{X}}_{t_m}(\Phi^{-1}(\bm{u}_n))\bigr).
\end{equation}
The estimator achieves a root mean square error of $O(N^{-1/2})$,
with an implied constant independent of $m$ provided that
the variance of the integrand is uniformly bounded.
This convergence rate is nevertheless slow.
Quasi-Monte Carlo (QMC) replaces independent random points with
low-discrepancy point sets and achieves nearly first-order
convergence under bounded-variation assumptions
\cite{caflisch1998,Niederreiter1992,dick2010}.
However, $\Phi^{-1}$ in \eqref{eq:intro-integral} can make the integrand unbounded near the
cube boundary, placing it outside the scope of the classical theory.
To address this limitation, approaches based on boundary growth
conditions \cite{owen2006a}, projection \cite{ouyang2024achieving},
and weighted function spaces \cite{kuo2014} have been developed. Nevertheless, the boundary-growth and projection-based approaches
do not directly resolve the dimension dependence.

The weighted-space framework of Nichols and Kuo \cite{kuo2014}
allows dimension-independent bounds, but applying it here
requires mixed derivative estimates and summability conditions
on the function-space weights.
The nested composition structure of the EM integrand makes
such estimates difficult to obtain.
It remains unclear whether {the dependence of this integrand} on the Gaussian
inputs permits weights satisfying these conditions uniformly
in $m$.
Recent studies offer insights into possible low-dimensional
structure that may help address this question.
For EM approximations, Ouyang, Pan, and He
\cite{ouyang2026sdeqmc} derive randomized QMC (RQMC) error bounds
of the form $C_{\varepsilon,s}N^{-1+\varepsilon}$,
but the constants grow with $s$. 
Ouyang and Wang \cite{ouyang2026understanding} later establish effective-dimension bounds uniform in $m$.
Numerical experiments in both studies show advantages of
RQMC over MC.
These findings suggest that EM integrands may possess
low-dimensional structure favorable to QMC, raising the
possibility that suitable weights satisfy the required
summability conditions.

It is also worth noting that QMC performance depends on
how Brownian paths are constructed from Gaussian inputs.
Different constructions can be represented by replacing
$\Phi^{-1}({\bm u})$ in \eqref{eq:intro-integral} with
${U}\Phi^{-1}({\bm u})$, where ${U}\in\R^{s\times s}$ is orthogonal.
Such transformations leave the standard MC error unchanged
but alter the structure of the integrand, and hence can
affect the QMC error
\cite{leobacher2012fast,wang2012how}.
In particular, we focus on the Brownian bridge (BB)
construction, which first generates $\bm{W}_T$ and then samples
conditional midpoint values on successively finer intervals
\cite[Section~1.1.2]{morters2010brownian}.
Such a construction assigns the largest-scale variation of
the Brownian path to the Gaussian inputs generating $\bm{W}_T$,
while each subsequent refinement introduces Gaussian inputs
that affect smaller-scale fluctuations with decreasing amplitude.
This hierarchical structure may lead to low-dimensional
structure in the integrand that is favorable to QMC.

In this paper, we use randomly shifted rank-1 lattice rules,
an RQMC method whose points are generated by an integer
generating vector and a uniform random shift.
Combined with the Brownian bridge construction, these rules
are used to approximate $\mathcal I_m$ for SDEs with constant
diffusion.
Within the weighted function-space framework of Nichols
and Kuo \cite{kuo2014}, we establish mixed derivative estimates
and verify the required summability conditions on the weights.
To handle the nested composition structure of the EM integrand,
we use the multivariate Fa\`a di Bruno formula
\cite[Theorem~1]{levy2006why}{ (see also \cite{cons:1996,hardy2006combinatorics})} to expand its mixed derivatives.
We organize the resulting terms and their propagation through the EM
recursion using a rooted-tree representation, following the standard
terminology in \cite[Section~4.2]{levin2021discrete}.
The decay and localized supports of the BB basis functions
then allow us to obtain summability estimates uniform in $m$.

Since the resulting weights are not available in an explicitly
computable form, a component-by-component (CBC) construction
\cite{nuyens2006fast,kuo2014} 
is impractical for these weights.
To avoid this difficulty, we choose generating vectors
at random and use the median of independent lattice estimates \cite{goda2022construction,ye2026median}.
Specifically, we independently generate $R=2k-1$ lattice
estimates, each using a uniformly random generating vector
and an independent uniform shift, and take their median.
This approach requires no explicit knowledge of the weights.
Taking the median also improves stability by reducing
the influence of unusually large errors and increasing
the probability of obtaining an accurate estimate.

Our analysis of the randomly shifted lattice rule and its
median estimator yields the following main results on
sampling error.
\begin{itemize}
	\item
	For an $N$-point randomly shifted lattice rule with a
	uniformly random generating vector, we prove a mean
	absolute sampling error bound of $O(N^{-1+\varepsilon})$
	for every $\varepsilon\in(0,1/2)$ and every prime $N$.
	The implied constant is independent of the number of
	time steps $m$. See Theorem~\ref{thm:main}.
	
	\item
	The median of $R=2k-1$ independent lattice estimates
	achieves an $L^q$ sampling error of $O(N^{-1+\varepsilon})$
	for every $1\le q<k/(1-\varepsilon)$, with an implied
	constant independent of $m$.
	In particular, the median of three estimates achieves
	the same rate in root mean square error. See Theorem~\ref{thm:main}.
	
	\item
	Under an additional summability condition on the columns
	of the diffusion matrix and suitable uniform bounds on
	the model data, the implied constants in the error bounds
	for both estimators are further independent of the state
	and {Brownian path dimensions}. See Corollary~\ref{cor:dimensions}.
\end{itemize}

The remaining time discretization error is covered by the classical
strong and weak convergence theory for Euler--Maruyama approximations;
see, for example, \cite[Chapters~10 and~14]{kloeden1992numerical}
and \cite[Section~5]{higham2001algorithmic}.
We retain a general assumption on the Euler--Maruyama time discretization error and verify it under
two sets of time regularity conditions. Theorems~\ref{thm:euler-strong-half}
and \ref{thm:euler-bias} establish bias rates of order $m^{-1/2}$
and $m^{-1}$, respectively.
The sampling error estimates and the Euler--Maruyama time discretization error bound together give
the total error bounds for both estimators in Corollary~\ref{cor:total}.

The remainder of the paper is organized as follows. Section~\ref{sec:setting} states
the assumptions on the SDE \eqref{eq:sde} and payoff $g$ and describes
the Brownian bridge construction.
Section~\ref{sec:main} states the sampling and total error results.
Section~\ref{sec:regularity} establishes the mixed derivative summability
result, and Section~\ref{sec:lattice} derives the lattice sampling error bounds.
Section~\ref{sec:numerics} presents the numerical experiments.
Section~\ref{sec:discussion} concludes the paper.
Technical proofs are deferred to the appendices.

\section{Notation and preliminaries}\label{sec:setting}
This section states the assumptions used in our analysis
and reviews the Brownian bridge construction, the
Euler--Maruyama scheme, the weighted function spaces of
Nichols and Kuo \cite{kuo2014}, and randomly shifted
lattice rules.
We also introduce the median lattice estimator following
\cite{goda2022construction,ye2026median}.
Before proceeding, we fix some notation.

We write $\Z$ for the integers, $\N=\{1,2,\ldots\}$ for the
positive integers, and $[s]=\{1,\ldots,s\}$ for $s\in\N$.
For a finite set $u$, we denote its cardinality by $|u|$
and the set of its partitions into nonempty blocks by $\PP(u)$,
with $\PP(\varnothing)=\{\varnothing\}$. For $\pi\in\PP(u)$, define
the set of labels occurring as singleton blocks of $\pi$ by
\begin{equation}\label{eq:singleton-blocks}
 \operatorname{sing}(\pi):=\{j\in u:\{j\}\in\pi\}.
\end{equation}
For an interval $I\subseteq\R$, $|I|$ denotes its length.
Throughout the paper, $d_x$ and $d_w$ denote the state and
{Brownian path dimensions}, respectively. An unnamed $C$ denotes a positive
constant that may change from line to line. Named constants retain
the definitions and parameter dependence specified when they are introduced.
We write $\phi(z)=(2\pi)^{-1/2}e^{-z^2/2}$ and $\Phi$ for the
standard Gaussian density and distribution function, respectively,
and $\phi_s(\bm{z})=\prod_{j=1}^s\phi(z_j)$ for the standard Gaussian
density on $\R^s$. The inverse $\Phi^{-1}$ is applied componentwise
to vector arguments.

Vectors are denoted by bold symbols and their scalar components
by unbolded symbols. The $r$th standard basis vector is denoted
by $\bm{e}_r$, with its dimension determined by the context.
For $d\in\N$, $\operatorname{Id}_d$ denotes the $d\times d$
identity matrix.
For $\bm{z}=(z_1,\ldots,z_s)$ and $u\subseteq[s]$, we write
$\bm{z}_u=(z_j)_{j\in u}$ and $\bm{z}_{-u}=(z_j)_{j\notin u}$.
For a sufficiently smooth scalar-valued function $f:\R^s\to\R$,
we define
\[
\partial_u f(\bm{z})
:=\left(\prod_{j\in u}\frac{\partial}{\partial z_j}\right)f(\bm{z}),
\qquad \partial_\varnothing f=f.
\]
Thus $\partial_u$ differentiates once with respect to each
coordinate indexed by $u$.

For $\bm{z}\in\R^s$ and a matrix
$A=(A_{ij})\in\R^{r\times w}$, we use the
Euclidean vector norm, the induced matrix $2$-norm, and the
Hilbert--Schmidt matrix norm:
\[
 |\bm{z}|_2=\left(\sum_{j=1}^s|z_j|^2\right)^{1/2},\qquad
 \|A\|_2=\sup_{\bm{v}\in\R^{w},\,|\bm{v}|_2\le1}|A\bm{v}|_2,\qquad
 \|A\|_{\mathrm{HS}}=\left(\sum_{i=1}^{r}\sum_{j=1}^{w}|A_{ij}|^2\right)^{1/2}.
\]
For $1\le q<\infty$, we define the $L^q$ norms of a scalar-valued
random variable $V$ and a vector-valued random variable $\bm{V}$ by
\[
 \|V\|_{L^q}=\bigl(\E[|V|^q]\bigr)^{1/q},
 \qquad
 \|\bm{V}\|_{L^q}=\bigl(\E[|\bm{V}|_2^q]\bigr)^{1/q},
\]
respectively.

We also fix notation for Fr\'echet derivatives and their norms.
For a smooth function $f:\R^{d_x}\to\R^{d_f}$,
the $n$th Fr\'echet derivative at $\bm{x}\in\R^{d_x}$,
\[
D_x^n f(\bm{x}):(\R^{d_x})^n\to\R^{d_f},
\qquad n\ge1,
\]
is a symmetric $n$-linear map.
For spatial direction vectors
$\bm{v}_1,\ldots,\bm{v}_n\in\R^{d_x}$, it is given by
\begin{equation}\label{eq:derivative-notation}
 D_x^nf(\bm{x})[\bm{v}_1,\ldots,\bm{v}_n]
 =\sum_{i_1,\ldots,i_n=1}^{d_x}
 \partial_{x_{i_1}}\cdots\partial_{x_{i_n}}f(\bm{x})
 (\bm{v}_1)_{i_1}\cdots(\bm{v}_n)_{i_n}.
\end{equation}
We define the operator norm and its supremum over $\bm{x}$ by
\[
 \|D_x^nf(\bm{x})\|_{\mathrm{op}}=\sup_{|\bm{v}_1|_2,\ldots,|\bm{v}_n|_2\le1}|D_x^nf(\bm{x})[\bm{v}_1,\ldots,\bm{v}_n]|_2,
 \qquad \|D_x^nf\|_\infty=\sup_{\bm{x}\in\R^{d_x}}\|D_x^nf(\bm{x})\|_{\mathrm{op}}.
\]

\subsection{Assumptions}
In this subsection, we state the assumptions used in our analysis.
Smoothness and derivative growth conditions on the drift and
payoff allow us to estimate mixed derivatives of the integrand
defining $\mathcal I_m$ \eqref{eq:intro-integral}.
Additional uniform bounds on the model data and a summability
condition on the diffusion columns yield error bounds independent
of the state and {Brownian path dimensions}.

\begin{assumption}\label{ass:space}
The drift $b:[0,T]\times\R^{d_x}\to\R^{d_x}$ is Borel measurable,
$b(t,\cdot)\in C^\infty$ for every $t$, and all its spatial derivatives are
jointly Borel measurable. There are constants $K_0,B_0,B_1>0$ and a
fixed $\tau\in[1,3/2)$ such that
\begin{align}
 \sup_{t\in[0,T]}|b(t,\bm{0})|_2&\le K_0,\label{eq:origin}\\
 \sup_{t\in[0,T],\,\bm{x}\in\R^{d_x}}\norm{D_x^{n}b(t,\bm{x})}_{\mathrm{op}}
 &\le \mathsf B_n:=B_0B_1^n(n!)^\tau,\qquad n\ge1.\label{eq:drift-ass}
\end{align}
The payoff $g:\R^{d_x}\to\R$ is $C^\infty$ and satisfies,
for some constants $G_0,G_1,\beta_g>0$,
\begin{equation}\label{eq:payoff-ass}
 \norm{{D_x^{n}}g(\bm{x})}_{\mathrm{op}}\le \mathsf G_n e^{\beta_g|\bm{x}|_2},\qquad \mathsf G_n:=G_0G_1^n(n!)^\tau,\quad n\ge1.
\end{equation}
\end{assumption}
In particular, the drift $b$ is globally Lipschitz in the spatial
variable, uniformly in time.
The payoff $g$ has at most exponential growth, as follows by
integrating its first derivative along the segment from
$\bm{0}$ to $\bm{x}$. Specifically,
\begin{equation}\label{eq:growth}
\begin{aligned}
 |b(t,\bm{x})-b(t,\bm{y})|_2
 &\le \mathsf B_1|\bm{x}-\bm{y}|_2,\\
 |b(t,\bm{x})|_2
 &\le K_0+\mathsf B_1|\bm{x}|_2,\\
 |g(\bm{x})|
 &\le |g(\bm{0})|+\frac{G_0G_1}{\beta_g}
      \bigl(e^{\beta_g|\bm{x}|_2}-1\bigr)
 \le C_g e^{\beta_g|\bm{x}|_2},
\end{aligned}
\end{equation}
where $C_g:=|g(\bm{0})|+G_0G_1/\beta_g$.

Assumption~\ref{ass:space} will be used to establish RQMC
sampling error bounds uniform in the number of time steps, with constants that
may depend on the {Brownian path dimension} $d_w$.
For bounds independent of $d_w$, we require an additional
summability condition on the column norms of the diffusion matrix.
To state this condition, define
\begin{equation}\label{eq:column-sums}
	\sigma_r:=|\Sigma\bm{e}_r|_2,\qquad
	\Theta_v:=\sum_{r=1}^{d_w}\sigma_r^v,\qquad v>0.
\end{equation}

\begin{assumption}\label{ass:uniform-data}
	Assumption~\ref{ass:space} holds with constants
	$K_0,B_0,B_1,G_0,G_1,\beta_g$ and $\tau$ independent of $d_x,d_w$.
	Moreover, $|\bm{x}_0|_2\le R_0$ for a constant $R_0<\infty$
	independent of $d_x,d_w$.
\end{assumption}

\begin{assumption}\label{ass:columns}
For some constant $M_{2/3}<\infty$ independent of $d_x,d_w$,
\begin{equation}\label{eq:column-uniform}
 \Theta_{2/3}=\sum_{r=1}^{d_w}|\Sigma \bm{e}_r|_2^{2/3}\le M_{2/3}.
\end{equation}
\end{assumption}
{For every $v\ge2/3$, the map $v\mapsto\Theta_v^{1/v}$ is nonincreasing, which gives $\Theta_v\le\Theta_{2/3}^{3v/2}$ without any further assumption. Under Assumption~\ref{ass:columns}, we therefore obtain}
\begin{equation}\label{eq:column-monotonicity}
	\Theta_v\le \Theta_{2/3}^{3v/2}\le M_{2/3}^{3v/2},
	\qquad
	\|\Sigma\|_{\mathrm{HS}}^2=\Theta_2\le M_{2/3}^3.
\end{equation}

\subsection{Brownian bridge construction}\label{sec:bridge}
In this subsection, we review the Brownian bridge construction
and the properties used in our analysis.
Throughout the remainder of the paper, $m$ denotes the number of
time steps and belongs to $\mathcal M:=\{2^J:J=0,1,\ldots\}$.
For each $m\in\mathcal M$, we use the uniform time grid
$\mathscr G:=\{t_i=ih:i=0,\ldots,m\}$ with step size $h=T/m$.
All bounds and suprema over $m$ refer to this family of grids.

\paragraph{Construction.} We use the dyadic Brownian bridge construction, corresponding to
L\'evy's construction of Brownian motion
\cite[Section~1.1.2]{morters2010brownian}.
Let $ \boldsymbol W = (W^{(1)},\dots,W^{(d_w)})^{\top}$. For each Brownian component $r=1,\ldots,d_w$, first set
\[
 W_0^{(r)}=0,\qquad W_T^{(r)}=\sqrt T\,z_r,
\]
where $z_1,\ldots,z_{d_w}$ are independent standard Gaussian inputs.
At each refinement level $\ell=0,\ldots,J-1$, define
\[
 u_{\ell,{{\mu}}}={{\mu}}T2^{-\ell},\qquad
 v_{\ell,{{\mu}}}=({{\mu}}+1)T2^{-\ell},\qquad
 c_{\ell,{{\mu}}}=\frac{u_{\ell,{{\mu}}}+v_{\ell,{{\mu}}}}2,
 \qquad {{\mu}}=0,\ldots,2^\ell-1.
\]
Generate the midpoint values recursively by
\begin{equation}\label{eq:bridge}
 W_{c_{\ell,{{\mu}}}}^{(r)}
 =\frac{W_{u_{\ell,{{\mu}}}}^{(r)}+W_{v_{\ell,{{\mu}}}}^{(r)}}{2}
  +\frac{\sqrt{v_{\ell,{{\mu}}}-u_{\ell,{{\mu}}}}}{2}\,z_{\ell,{{\mu}},r},
\end{equation}
proceeding from level $0$ to level $J-1$ for every component.
The input $z_{\ell,{{\mu}},r}$ corresponds to interval ${{\mu}}$ at level $\ell$
in Brownian component $r$. These detail inputs are independent standard
Gaussian inputs, also independent of the endpoint inputs.

We order the inputs with the $d_w$ endpoint inputs first, followed by
the detail inputs in lexicographic order of $(\ell,{{\mu}},r)$.
Let $j$ denote the position of an input in this ordering. For detail inputs,
\begin{equation}\label{eq:detail-index}
	z_{\ell,{{\mu}},r}=z_j,\qquad j=d_w(2^\ell+{{\mu}})+r,
	\qquad 0\le\ell<J,\quad 0\le {{\mu}}<2^\ell,\quad 1\le r\le d_w.
\end{equation}
Collecting the ordered inputs gives
$\bm{z}=(z_1,\ldots,z_s)\in\R^s$, where $s=d_wm$.
The endpoint and detail index sets are
\begin{equation}\label{eq:bb-index-sets}
	\mathcal E_w=\{1,\ldots,d_w\},\qquad
	\mathcal J_m=\{d_w+1,\ldots,s\}.
\end{equation}
For $j\in\mathcal J_m$, let $(\ell(j),{{\mu}}(j),r(j))$ denote the unique
triple determined by \eqref{eq:detail-index}, so $r(j)$ identifies the
Brownian component constructed using input $z_j$.
For $j\in \mathcal E_w$, set $r(j)=j$. When $J=0$, there are no detail inputs.

Writing the resulting grid values as
$\bm{W}_{t_i}(\bm{z})=(W_{t_i}^{(r)}(\bm{z}))_{r=1}^{d_w}$,
define $\mathcal B_m\in\R^{s\times s}$ by the standardized increments
in time order:
\begin{equation}\label{eq:bb-map}
 \mathcal B_m\bm{z}:=
 \left(\frac{\bm{W}_{t_i}(\bm{z})-\bm{W}_{t_{i-1}}(\bm{z})}{\sqrt h}\right)_{i=1}^{m}.
\end{equation}
Since the standardized increments are independent standard
Gaussian variables under standard Gaussian input,
$\mathcal B_m$ is orthogonal.

\paragraph{Coefficient functions and their properties.}
The BB grid values admit a finite L\'evy--Ciesielski representation
in Faber--Schauder hat functions
\cite[Proposition~1.44 and Remark~1.45]{morters2010brownian}.
To describe its coefficients, define the unit-height hat function
\[
 \Lambda(x)=
 \begin{cases}
 2x,&0\le x\le\tfrac12,\\
 2(1-x),&\tfrac12<x\le1,\\
 0,&\text{otherwise},
 \end{cases}
\]
and the scalar coefficient functions, for $t\in[0,T]$, by
\begin{equation}\label{eq:hat-explicit}
 \vartheta_j(t)=
 \begin{cases}
 \displaystyle\frac{t}{\sqrt T},&j\in \mathcal E_w,\\[6pt]
 \displaystyle\frac{\sqrt T}{2}\,2^{-\ell(j)/2}
 \Lambda\!\left(\frac{2^{\ell(j)}t}{T}-{{\mu}}(j)\right),&j\in\mathcal J_m.
 \end{cases}
\end{equation}
Here $\vartheta_j$ is the coefficient of input $z_j$ in Brownian
component $r(j)$. For $j\in\mathcal J_m$, its supporting interval is
\begin{equation}\label{eq:detail-support}
 I_j=[{{\mu}}(j)T2^{-\ell(j)},({{\mu}}(j)+1)T2^{-\ell(j)}].
\end{equation}
Define the corresponding coefficients in the diffusion term by
\begin{equation}\label{eq:diffusion-coefficients}
 \bm{\xi}_j(t):=\Sigma\bm{e}_{r(j)}\vartheta_j(t),\qquad
 \|\bm{\xi}_j\|_\infty:=\sup_{0\le t\le T}|\bm{\xi}_j(t)|_2.
\end{equation}

For $\bm{z}\in\R^s$ and $i=0,\ldots,m$, the resulting linear
representations are
\cite[Proposition~1.44 and Remark~1.45]{morters2010brownian}
\begin{equation}\label{eq:linearbridge}
 \bm{W}_{t_i}(\bm{z})=\sum_{j=1}^{s}\bm{e}_{r(j)}\vartheta_j(t_i)z_j,
 \qquad
 \Sigma\bm{W}_{t_i}(\bm{z})=\sum_{j=1}^{s}\bm{\xi}_j(t_i)z_j.
\end{equation}

\begin{lemma}\label{lem:bridge}
For $j\in\mathcal J_m$ with $\ell(j)=\ell$, we have
\begin{equation}\label{eq:bridgebounds}
 \begin{aligned}
 |I_j|&=T2^{-\ell},&
 \|\vartheta_j\|_\infty&=\frac{\sqrt T}{2}\,2^{-\ell/2},\\
 \|\bm{\xi}_j\|_\infty
 &=\tfrac12\sigma_{r(j)}\sqrt{|I_j|},&
 \bm{\xi}_j(T)&=0.
 \end{aligned}
\end{equation}
For $j\in \mathcal E_w$,
\[
 \vartheta_j(T)=\|\vartheta_j\|_\infty=\sqrt T,
 \qquad
 {|\bm{\xi}_j(T)|_2}=\|\bm{\xi}_j\|_\infty=\sqrt T\,\sigma_j.
\]
\end{lemma}
\begin{proof}
	The assertions follow directly from
	\eqref{eq:hat-explicit}, \eqref{eq:detail-support},
	and the definition of $\bm{\xi}_j$.
\end{proof}

\subsection{Euler--Maruyama integrand}\label{sec:euler-rules}

Recall the EM recursion \eqref{eq:em}. Under the BB construction,
its standardized increments are given by $\mathcal B_m\bm{z}$.
Define
\begin{equation}\label{eq:bb-euler-state}
 \bm{X}^h_{t_i}(\bm{z})
 :=\widetilde{\bm{X}}_{t_i}(\mathcal B_m\bm{z}),
 \qquad i=0,\ldots,m.
\end{equation}
All input derivatives below are taken with respect to the BB coordinates
$\bm{z}$. By \eqref{eq:em} and \eqref{eq:bb-map}, we have
$\bm{X}^h_{t_0}(\bm{z})=\bm{x}_0$ and
\begin{equation}\label{eq:bb-em}
		\bm{X}^h_{t_{i+1}}(\bm{z})
		=\bm{X}^h_{t_i}(\bm{z})
		+h\,b\bigl(t_i,\bm{X}^h_{t_i}(\bm{z})\bigr)
		+\Sigma\bigl(\bm{W}_{t_{i+1}}(\bm{z})-\bm{W}_{t_i}(\bm{z})\bigr),
		\qquad 0\le i<m.
\end{equation}
We define the corresponding integrand by
\begin{equation}\label{eq:integrand}
	F_m(\bm{z})
	:=g\bigl(\bm{X}^h_{t_m}(\bm{z})\bigr),
	\qquad \bm{z}\in\R^s.
\end{equation}
Since $\mathcal B_m$ is orthogonal, the BB transformation
leaves the Gaussian integral unchanged. Thus,
\begin{equation}\label{eq:bb-integral}
	\mathcal I_m=\E[F_m(\bm{Z})]
	=\int_{\R^s}F_m(\bm{z})\phi_s(\bm{z})\dd\bm{z},
	\qquad \bm{Z}\sim N(0,{\operatorname{Id}_s}).
\end{equation}
We will analyze $F_m$ in the weighted function space
introduced in Section~\ref{sec:spaces-rules}.

\subsection{Randomly shifted lattice rules and median estimators}\label{sec:lattice-estimators}
In this subsection, we review randomly shifted rank-1 lattice rules and
median estimators following \cite{goda2022construction,ye2026median}.

For a prime $N\ge2$, let $\bm{a}\in\{1,\ldots,N-1\}^s$ be a generating
vector and let $\bm{\Delta}$ be uniformly distributed on $[0,1]^s$.
The lattice points are
\begin{equation}\label{eq:intro-lattice-points}
 \bm{u}_n=\{n\bm{a}/N+\bm{\Delta}\},\qquad n=0,\ldots,N-1,
\end{equation}
where braces denote componentwise fractional parts.
For the Gaussian integral
$\mathcal I(F):=\int_{\R^s}F(\bm{z})\phi_s(\bm{z})\dd\bm{z}$,
the randomly shifted lattice estimator is
\[
Q_{\bm{a},\bm{\Delta}}(F)
:=\frac1N\sum_{n=0}^{N-1}
F\bigl(\Phi^{-1}(\bm{u}_n)\bigr).
\]
The uniform shift makes each lattice point uniformly distributed,
so $\E_{\bm{\Delta}}[Q_{\bm{a},\bm{\Delta}}(F)]=\mathcal I(F)$ for integrable $F$. In particular, our estimate of $\mathcal I_m$ is
\begin{equation}\label{eq:Q}
 Q_{\bm{a},\bm{\Delta}}^{N,m}:=Q_{\bm{a},\bm{\Delta}}(F_m)
 =\frac1N\sum_{n=0}^{N-1}
 F_m\bigl(\Phi^{-1}(\{n\bm{a}/N+\bm{\Delta}\})\bigr).
\end{equation}

We now follow the median approach in
\cite{goda2022construction,ye2026median} to combine independent
randomly shifted lattice estimates.
For an odd integer $R=2k-1$, $k\ge1$, draw $R$ independent pairs
$(\bm{A}_1,\bm{\Delta}_1),\ldots,(\bm{A}_R,\bm{\Delta}_R)$,
where each $\bm{A}_{{\nu}}$ is uniform on $\{1,\ldots,N-1\}^s$
and each $\bm{\Delta}_{{\nu}}$ is uniform on $[0,1]^s$,
independently of $\bm{A}_{{\nu}}$.
The median estimator is
\begin{equation}\label{eq:median-def}
 M_{R,N,m}:=\operatorname{median}
 \bigl(Q_{\bm{A}_1,\bm{\Delta}_1}^{N,m},\ldots,Q_{\bm{A}_R,\bm{\Delta}_R}^{N,m}\bigr).
\end{equation}
This estimator requires no dedicated construction of generating
vectors. Taking the median can reduce the probability of a large
integration error; see, for example,
\cite{pan2024super,goda2022construction,ye2026median}.

\subsection{Weighted Gaussian function space}\label{sec:spaces-rules}
In this subsection, we consider a Gaussian specialization of
the weighted function space framework of Nichols and
Kuo~\cite{kuo2014}. 

Fix $\alpha>0$ and define the auxiliary weight function
$\psi(t):=e^{-\alpha|t|}$, $t\in\R$.
Write
$\phi_{-u}(\bm{z}_{-u}):=\prod_{j\notin u}\phi(z_j)$.
For nonnegative subset weights $\bm{\gamma}:=(\gamma_u)_{u\subseteq[s]}$ with
$\gamma_\varnothing=1$, the weighted unanchored Sobolev space
$\mathcal H_{s,\bm{\gamma}}$ consists of locally integrable functions
$F:\R^s\to\R$ with the required weak mixed derivatives for which
the following norm is finite.
\begin{equation}\label{eq:NK-norm}
 \|F\|_{\mathcal H_{s,\bm{\gamma}}}^2
 :=\sum_{u\subseteq[s]}\frac{1}{\gamma_u}
 \int_{\R^{|u|}}
 \left|\int_{\R^{s-|u|}}\partial_uF(\bm{z})\phi_{-u}(\bm{z}_{-u})\dd \bm{z}_{-u}\right|^2
 \prod_{j\in u}\psi(z_j)^2\dd \bm{z}_u.
\end{equation}
{If $\gamma_u=0$, we require the corresponding numerator in \eqref{eq:NK-norm} to vanish for every function in $\mathcal H_{s,\bm\gamma}$,}
and we use the convention $0/0:=0$.
The subset weights $\gamma_u$ describe the relative importance of
coordinate subsets $u$ by scaling their contributions to the norm.

Following Nichols and Kuo~\cite{kuo2014}, we define the
shift-averaged worst-case error of the rules in
Section~\ref{sec:lattice-estimators} by
\begin{equation}\label{eq:shift-averaged-wce}
 [e^{\rm sh}_{s,N}(\bm{a};\bm{\gamma})]^2
 :=\E_{\bm{\Delta}}\sup_{\|G\|_{\mathcal H_{s,\bm{\gamma}}}\le1}
       |Q_{\bm{a},\bm{\Delta}}(G)-\mathcal I(G)|^2.
\end{equation}
The following lemma is a Gaussian specialization of
\cite[Theorem~1, Eq.~(13)]{ye2026median}.

\begin{lemma}
\label{lem:random-wce}
Fix $\lambda\in(1/2,1)$ and choose ${\delta_*}\in(1/(2\lambda),1)$.
Set
\begin{equation}\label{eq:rho}
 C_{\alpha,{\delta_*}}:=\frac{\sqrt{2\pi}\,e^{\alpha^2/(1-{\delta_*})}}
              {\pi^{2{\delta_*}}{\delta_*}(1-{\delta_*})},
 \qquad
 \rho_\lambda:=2C_{\alpha,{\delta_*}}^\lambda\zeta(2{\delta_*}\lambda).
\end{equation}
Here $\zeta$ is the Riemann zeta function.
For prime $N$ and nonnegative weights $\gamma_u$,
let $\bm{A}$ be uniform on $\{1,\ldots,N-1\}^s$. Then
\begin{equation}\label{eq:random-wce-moment}
 \E_{\bm{A}}\bigl[e^{\rm sh}_{s,N}(\bm{A};\bm{\gamma})^{2\lambda}\bigr]
 \le\frac{1}{N-1}
 \sum_{\varnothing\ne u\subseteq[s]}
       \gamma_u^\lambda\rho_\lambda^{|u|}.
\end{equation}
\end{lemma}

\begin{proof}
The constant $C_{\alpha,{\delta_*}}$ in \eqref{eq:rho} is given by
\cite[Example~5]{kuo2010} for the Gaussian density and auxiliary
weight used here. Since $N$ is prime, {exactly $N-1$ integers in $\{1,\ldots,N-1\}$ are coprime to $N$};
\cite[Theorem~1, Eq.~(13)]{ye2026median}, with
$C_j=C_{\alpha,{\delta_*}}$ and $r_j={\delta_*}$ for every coordinate {(in the notation of \cite{ye2026median})}, then gives
\eqref{eq:random-wce-moment}.
\end{proof}

\section{Main results}\label{sec:main}
This section presents the sampling error bounds for the random
lattice \eqref{eq:Q} and median estimators \eqref{eq:median-def}, with constants independent of the
number of time steps $m$. We then state an assumption on the
Euler--Maruyama time discretization error and combine the two error bounds to obtain estimates
for the total error.

\subsection{Sampling error}\label{sec:sampling-results}

\begin{theorem}\label{thm:main}
Suppose Assumption~\ref{ass:space} holds. For every
$\varepsilon\in(0,1/2)$, there exists a constant $C_\varepsilon>0$,
depending only on
$\varepsilon,\tau,T,K_0,B_0,B_1,G_0,G_1,\beta_g,|\bm x_0|_2,
\Theta_1,\Theta_{2/3}$, such that for every $m\in\mathcal M$
and every prime $N\ge2$,
the single random lattice rule satisfies
\begin{equation}\label{eq:single-l1-main}
 \norm{Q_{\bm{A},\bm{\Delta}}^{N,m}-\mathcal I_m}_{L^1}
 \le C_\varepsilon N^{-1+\varepsilon}.
\end{equation}
Let $R=2k-1$ for a fixed integer $k\ge1$. For every
$q$ satisfying $1\le q<k/(1-\varepsilon)$, there exists a constant
$C_{\varepsilon,q,R}>0$, depending only on the same parameters as
$C_\varepsilon$ and on $q,R$, such that the median
estimator in \eqref{eq:median-def} satisfies
\begin{equation}\label{eq:median-sampling-main}
 \norm{M_{R,N,m}-\mathcal I_m}_{L^q}
 \le C_{\varepsilon,q,R}N^{-1+\varepsilon}.
\end{equation}
The expectations are over the generating vectors and random shifts.
\end{theorem}

\begin{remark}
Taking the median of three independent estimates ($k=2$, $R=3$)
gives the root mean square error bound
\[
 \norm{M_{3,N,m}-\mathcal I_m}_{L^2}
 \le C_{\varepsilon,2,3}N^{-1+\varepsilon}.
\]
Increasing $k$ can improve the stability of the median estimator.
\end{remark}

With the other parameters fixed, the constants in Theorem~\ref{thm:main}
may still depend on $d_x,d_w$ through $\Theta_1$ and $\Theta_{2/3}$. The following corollary gives the same convergence rates with
constants independent of $d_x$ and $d_w$ under uniform bounds
on the model parameters and a summability condition on the
columns of $\Sigma$.

\begin{corollary}
\label{cor:dimensions}
Suppose the assumptions of Theorem~\ref{thm:main} hold, together with
Assumptions~\ref{ass:uniform-data} and \ref{ass:columns}.
Then the bounds \eqref{eq:single-l1-main} and \eqref{eq:median-sampling-main}
hold with constants $C_\varepsilon$ and $C_{\varepsilon,q,R}$ independent
of $m,N,d_x,d_w$. These constants may depend on $T$, the uniform model
bounds, and the parameters indicated by their subscripts.
\end{corollary}

Theorem~\ref{thm:main} and Corollary~\ref{cor:dimensions} are proved in
Section~\ref{sec:sampling-proof}.

\begin{remark}\label{rem:rotation}
	For any orthogonal matrix ${U}$, replacing $\Sigma$ by $\Sigma{U}$
	preserves the law of the SDE and its Euler--Maruyama approximation.
	Let $\varsigma_j$ denote the nonzero singular values of $\Sigma$.
	Then
	\[
	\min_{{U}^\top{U}={\operatorname{Id}_{d_w}}}
	\sum_{r=1}^{d_w}|\Sigma{U}\bm e_r|_2^{2/3}
	=
	\sum_{j=1}^{\operatorname{rank}(\Sigma)}\varsigma_j^{2/3}.
	\]
	Indeed, pad the singular values with zeros to obtain $d_w$ values and write
	\[
	\Sigma^\top\Sigma
	=V\operatorname{diag}(\varsigma_1^2,\ldots,\varsigma_{d_w}^2)V^\top.
	\]
	Here $V\in\mathbb R^{d_w\times d_w}$ is orthogonal.
	Set ${\widetilde U}=V^\top{U}=({\tilde u_{jr}})_{1\le j,r\le d_w}$.
	By concavity of $x^{1/3}$ and orthogonality of ${\widetilde U}$,
	\[
	\sum_{r=1}^{d_w}|\Sigma{U}\bm e_r|_2^{2/3}
	=\sum_{r=1}^{d_w}\Bigl(\sum_{j=1}^{d_w}{\tilde u_{jr}}^2\varsigma_j^2\Bigr)^{1/3}
	\ge\sum_{r,j=1}^{d_w}{\tilde u_{jr}}^2\varsigma_j^{2/3}
	=\sum_{j=1}^{d_w}\varsigma_j^{2/3}.
	\]
	Equality holds for ${U}=V$, the matrix of right singular vectors.
	Thus there exists a choice of orthogonal matrices across dimensions
	for which the transformed diffusion matrices satisfy
	Assumption~\ref{ass:columns} if and only if the singular-value sums
	are uniformly bounded.
\end{remark}

\subsection{Euler--Maruyama time discretization error}\label{sec:euler-bias}
Our main focus is on the sampling error relative to $\mathcal I_m$.
To obtain total error bounds relative to $\mathcal I$, we additionally
assume the following rate for the Euler--Maruyama time discretization error.

\begin{assumption}\label{ass:bias}
For some $\theta>0$ and $C_{\rm bias}<\infty$, independent of
$m\in\mathcal M$,
\begin{equation}\label{eq:bias-main}
 |\mathcal I_m-\mathcal I|\le C_{\rm bias}m^{-\theta}.
\end{equation}
\end{assumption}

The following two theorems verify Assumption~\ref{ass:bias} with
$\theta=1/2$ and $\theta=1$, respectively, under different time
regularity conditions on $b$.

\begin{theorem}\label{thm:euler-strong-half}
Suppose Assumption~\ref{ass:space} holds and, for some ${\widetilde L_t}<\infty$,
\begin{equation}\label{eq:strong-time-regularity}
 |b(t,\bm{x})-b(s,\bm{x})|_2
 \le {\widetilde L_t}(1+|\bm{x}|_2)|t-s|^{1/2},
 \qquad s,t\in[0,T],\quad \bm{x}\in\R^{d_x}.
\end{equation}
Let $\bm{X}$ solve the SDE~\eqref{eq:sde}, and let $\bm{X}^h$ be
the Euler--Maruyama approximation defined by \eqref{eq:bb-em},
driven by the same Brownian motion $\bm{W}$. Then the following conclusions hold.
\begin{enumerate}
\item[(i)] For fixed $d_x,d_w$, there are constants
$C_{{\rm str},1/2}$ and $C_{{\rm bias},1/2}$, independent of
$m\in\mathcal M$, such that
\begin{align}
 \left\|\max_{0\le i\le m}
       |\bm{X}^h_{t_i}-\bm{X}_{t_i}|_2\right\|_{L^2}
 &\le C_{{\rm str},1/2}m^{-1/2},\label{eq:strong-half-main}\\
 |\mathcal I_m-\mathcal I|
 &\le C_{{\rm bias},1/2}m^{-1/2}.\label{eq:bias-half-main}
\end{align}
The constant $C_{{\rm str},1/2}$ depends only on
$T,K_0,\mathsf B_1,{\widetilde L_t},|\bm{x}_0|_2$ and
$\|\Sigma\|_{\mathrm{HS}}$; $C_{{\rm bias},1/2}$ may also depend
on $G_0,G_1$ and $\beta_g$.
\item[(ii)] If Assumptions~\ref{ass:uniform-data} and
\ref{ass:columns} also hold, and ${\widetilde L_t}$ in
\eqref{eq:strong-time-regularity} is independent of $d_x,d_w$, then
both constants can be chosen independent of $m,d_x,d_w$.
They depend only on
$T,{\widetilde L_t},K_0,B_0,B_1,G_0,G_1,\beta_g,\tau,R_0$ and $M_{2/3}$.
\end{enumerate}
\end{theorem}

\begin{theorem}\label{thm:euler-bias}
Suppose Assumption~\ref{ass:space} holds and, for some $L_t<\infty$,
\begin{equation}\label{eq:bias-time-regularity}
 |b(t,\bm{x})-b(s,\bm{x})|_2
 \le L_t(1+|\bm{x}|_2)|t-s|,
 \qquad s,t\in[0,T],\quad \bm{x}\in\R^{d_x}.
\end{equation}
Then the following conclusions hold.
\begin{enumerate}
\item[(i)] For fixed $d_x,d_w$, there is a constant $C_{\rm bias}$,
independent of $m\in\mathcal M$, such that
\begin{equation}\label{eq:bias-order-one}
 |\mathcal I_m-\mathcal I|\le C_{\rm bias}m^{-1}.
\end{equation}
The constant depends only on
$T,K_0,\mathsf B_1,\mathsf B_2,L_t,G_0,G_1,\beta_g,
|\bm{x}_0|_2$ and $\|\Sigma\|_{\mathrm{HS}}$.
\item[(ii)] If Assumptions~\ref{ass:uniform-data} and
\ref{ass:columns} also hold, and $L_t$ in
\eqref{eq:bias-time-regularity} is independent of $d_x,d_w$, then
$C_{\rm bias}$ in \eqref{eq:bias-order-one} can be chosen independent
of $m,d_x,d_w$. It depends only on
$T,L_t,K_0,B_0,B_1,G_0,G_1,\beta_g,\tau,R_0$ and $M_{2/3}$.
\end{enumerate}
\end{theorem}

The proofs are given in Appendix~\ref{sec:moments}.

\subsection{Total error}\label{sec:total-error}
Combining the sampling estimates with Assumption~\ref{ass:bias}
gives the following total error bounds.

\begin{corollary}\label{cor:total}
Under Assumptions~\ref{ass:space} and \ref{ass:bias}, the rules in
Theorem~\ref{thm:main} satisfy
\begin{align}
 \norm{Q_{\bm{A},\bm{\Delta}}^{N,m}-\mathcal I}_{L^1}
 &\le C_\varepsilon N^{-1+\varepsilon}
      +C_{\rm bias}m^{-\theta},\label{eq:single-l1-total}\\
 \norm{M_{R,N,m}-\mathcal I}_{L^q}
 &\le C_{\varepsilon,q,R}N^{-1+\varepsilon}
      +C_{\rm bias}m^{-\theta},\label{eq:total-main}
\end{align}
with the same ranges of $m,N,\varepsilon,q,R$ as in the sampling results.
Under the assumptions of Corollary~\ref{cor:dimensions}, the sampling
constants $C_\varepsilon$ and $C_{\varepsilon,q,R}$ {can be chosen independent of $d_x,d_w$}.
All constants are then independent of $d_x,d_w$ if $\theta$ and $C_{\rm bias}$ in
Assumption~\ref{ass:bias} are common to all dimensions.
\end{corollary}
\begin{proof}
Apply the triangle inequality to each sampling estimate and
\eqref{eq:bias-main}.
\end{proof}

\section{Regularity of the Euler--Maruyama integrand}\label{sec:regularity}
In this section, we establish mixed derivative bounds for the Euler--Maruyama
integrand $F_m$ and prove their summability uniformly in $m$.
These estimates will be used in Section~\ref{sec:sampling-proof} to
derive the sampling error bounds.

Unless otherwise stated, constants in the estimates below are independent
of $m$. Their dependence on other parameters
is specified when relevant.

For $\varnothing\ne u\subseteq[s]$, define the normalized derivative
envelope by
\begin{equation}\label{eq:normalized-derivative-envelope}
 \mathcal D_{m,u}:=\sup_{\bm{z}\in\R^s}
 \frac{|\partial_uF_m(\bm{z})|}
 {e^{\beta_g|\bm{X}^h_{t_m}(\bm{z})|_2}}.
\end{equation}

\begin{theorem}
	\label{thm:derivatives}
	Suppose Assumption~\ref{ass:space} holds.
	For every $\kappa>0$ and $p>2/3$, there exists a constant
	$C_{\mathcal D}(\kappa,p;\Theta_{2/3})>0$, depending only on
	$\kappa,p,\tau,T,B_0,B_1,G_0,G_1,\Theta_{2/3}$ and nondecreasing in
	$\kappa$ and $\Theta_{2/3}$, such that
	\begin{equation}\label{eq:derivative-sum}
		\sup_{m\in\mathcal M}
		\sum_{\varnothing\ne u\subseteq[s]}
		\kappa^{|u|}\mathcal D_{m,u}^p
		\le C_{\mathcal D}(\kappa,p;\Theta_{2/3}).
	\end{equation}
\end{theorem}

The proof is completed in Section~\ref{sec:summability-proof}.
The auxiliary results are proved in
Appendices~\ref{sec:variation-proof}--\ref{sec:regularity-proofs}.

\subsection{Recursive expansions of state derivatives}
We first derive recursive expansions for the mixed derivatives of the Euler--Maruyama state $ \bm{X}^{h}$
with respect to the Brownian bridge input coordinates.
Fix $\bm{z}$ and suppress the argument of $\bm{X}^h_{t_i}(\bm{z})$
and $\bm{W}_{t_i}(\bm{z})$ when convenient.
The set-partition Fa\`a di Bruno formula
\cite[Theorem~1]{levy2006why}{ (see also \cite{hardy2006combinatorics})}
gives
\begin{equation}\label{eq:chain}
 \partial_u(f\circ{\bm\Psi})=\sum_{\pi\in\PP(u)}
 {D_x^{|\pi|}}f({\bm\Psi})[\partial_v{\bm\Psi}:v\in\pi],\qquad u\ne\varnothing.
\end{equation}
For a fixed grid with $h=T/m$, set
$\bm{Y}_i(\bm{z})=\bm{X}^h_{t_i}(\bm{z})-\Sigma \bm{W}_{t_i}(\bm{z})$.
It follows from \eqref{eq:bb-em} that, after suppressing $\bm{z}$,
\begin{equation}\label{eq:Y}
 \bm{Y}_{i+1}=\bm{Y}_i+hb(t_i,\bm{Y}_i+\Sigma \bm{W}_{t_i}),\qquad \bm{Y}_0=\bm{x}_0.
\end{equation}
For every nonempty $v\subseteq[s]$, define
\[
 \bm{Y}_{i,v}=\partial_v\bm{Y}_i,\qquad
 {\bm X^h_{i,v}}=\partial_v\bm{X}^h_{t_i}.
\]
By \eqref{eq:linearbridge}, these derivatives satisfy
\[
 {\bm X^h_{i,v}}=\bm{Y}_{i,v}
 +\begin{cases}
   \bm{\xi}_j(t_i),&v=\{j\},\\
   0,&|v|\ge2.
  \end{cases}
\]
\begin{lemma}\label{lem:state-expansion}
Set
\begin{equation}\label{eq:J}
 {\Omega}_i:=\operatorname{Id}_{d_x}+hD_xb(t_i,\bm{X}^h_{t_i}),
 \qquad 0\le i<m.
\end{equation}
For $0\le {{i'}}<i$, define $P_{i,{{i'}}}={\Omega}_{i-1}\cdots {\Omega}_{{{i'}}+1}$, an empty product
being $\operatorname{Id}_{d_x}$. For $i\in\{0,\ldots,m\}$ and
$j\in[s]$, we have
\begin{equation}\label{eq:partition-first}
 \bm{Y}_{i,\{j\}}=
 \sum_{{{i'}}=0}^{i-1}P_{i,{{i'}}}hD_xb(t_{{i'}},\bm{X}^h_{t_{{i'}}})
 [\bm{\xi}_j(t_{{i'}})].
\end{equation}
For $u\subseteq[s]$ with $|u|\ge2$, we have
\begin{equation}\label{eq:correspondence-reindex}
 \bm{Y}_{i,u}={}
 \sum_{{{i'}}=0}^{i-1}\sum_{S\subseteq u}
 \sum_{\substack{\pi\in\PP(u\setminus S)\\|S|+|\pi|\ge2}}
 P_{i,{{i'}}}hD_x^{|S|+|\pi|}b(t_{{i'}},\bm{X}^h_{t_{{i'}}})
 \left[(\bm{\xi}_j(t_{{i'}}))_{j\in S},
 (\bm{Y}_{{{i'}},v})_{v\in\pi}\right].
\end{equation}
Moreover,
\begin{equation}\label{eq:propagator-bound}
 \norm{P_{i,{{i'}}}}_2\le e^{\mathsf B_1T},
\end{equation}
where $\mathsf B_1$ is defined in \eqref{eq:drift-ass}.
\end{lemma}
The proof is given in Appendix~\ref{sec:variation-proof}.

\subsection{Derivative majorant and summability estimates}
To control each $\mathcal D_{m,u}$, we define the auxiliary quantity
\begin{equation}\label{eq:eta}
 \eta_m(S)=h\sum_{i=0}^{m-1}\prod_{j\in S}|\bm{\xi}_j(t_i)|_2,
 \qquad \varnothing\ne S\subseteq[s].
\end{equation}
Repeatedly substituting the recursive expansions in
Lemma~\ref{lem:state-expansion}
into the expression for $\partial_uF_m$, and using rooted trees to keep
track of the resulting terms, yields the following estimate for
$\mathcal D_{m,u}$. The proof is given in
Appendix~\ref{sec:expansion-counting}.
\begin{proposition}\label{prop:majorant}
For $u\ne\varnothing$ let $\mathcal D_{m,u}$ be the normalized derivative envelope defined in \eqref{eq:normalized-derivative-envelope}. Under
Assumption~\ref{ass:space} there is $C_{\rm maj}\ge1$, depending only on
$T,B_0,B_1,G_0,G_1,\tau$, such that
\begin{equation}\label{eq:majorant}
 \mathcal D_{m,u}\le C_{\rm maj}^{|u|}
 \sum_{E\subseteq u\cap \mathcal E_w}(|E|!)^\tau
 \prod_{j\in E}(\sqrt T\,\sigma_j)
 \sum_{\pi\in\PP(u\setminus E)}
 C_{\rm maj}^{|\pi|}(|\pi|!)^\tau
 \prod_{S\in\pi}(|S|!)^\tau\eta_m(S).
\end{equation}
\end{proposition}

The following two technical lemmas control the summability of
$\eta_m(S)$ over $\varnothing\ne S\subseteq\mathcal J_m$ and over
all nonempty input subsets. Their proofs are given in
Appendix~\ref{sec:regularity-proofs}.

\begin{lemma}\label{lem:detail-sum}
Let $\varnothing\ne S\subseteq\mathcal J_m$, and recall the supporting
intervals $I_j$ from \eqref{eq:detail-support}. If
$\bigcap_{j\in S}\operatorname{int}(I_j)=\varnothing$, then
$\eta_m(S)=0$. Otherwise choose $j_*\in S$ such that
$|I_{j_*}|=\min_{j\in S}|I_j|$. Then
\begin{equation}\label{weighted-laminar}
 \eta_m(S)\le \frac{3}{2}\,|I_{j_*}|
       \prod_{j\in S}\left(\frac{\sigma_{r(j)}}{2}|I_j|^{1/2}\right).
\end{equation}

For every $\rho>2/3$ there exists $C_{\rho,T}>0$, depending only on
$\rho,T$, such that, for every
integer $n\ge1$,
\begin{equation}\label{eq:detail-subset-sum}
 \sup_{m\in\mathcal M}
 \sum_{\substack{S\subseteq\mathcal J_m\\|S|=n}}\eta_m(S)^\rho
       \le\frac{C_{\rho,T}^n \Theta_\rho^n}{(n-1)!}.
\end{equation}
\end{lemma}

\begin{lemma}\label{lem:all-subset-sum}
For every $2/3<\rho<1/\tau$ and $\kappa>0$,
define
\begin{equation}\label{eq:endpoint-subset-sum-definition}
 {\Upsilon}^{\mathcal E}(\kappa):=\sum_{E\subseteq\mathcal E_w}
 \kappa^{|E|}(|E|!)^{\tau\rho}\prod_{j\in E}(\sqrt T\,\sigma_j)^\rho,
\end{equation}
and
\begin{equation}\label{eq:all-subset-sum-definition}
 {\Upsilon}_m(\kappa):=\sum_{\varnothing\ne S\subseteq[s]}
     \kappa^{|S|}(|S|!)^{\tau\rho}\eta_m(S)^\rho.
\end{equation}
Then
\[
 \max\left\{{\Upsilon}^{\mathcal E}(\kappa),
 \sup_{m\in\mathcal M}{\Upsilon}_m(\kappa)\right\}
 \le C_{{\Upsilon}}(\kappa,\rho,\tau,T;\Theta_\rho)<\infty,
\]
where the finite constant $C_{{\Upsilon}}(\kappa,\rho,\tau,T;\Theta_\rho)$ depends only on
$\kappa,\rho,\tau,T,\Theta_\rho$ and can be chosen nondecreasing in
$\kappa$ and $\Theta_\rho$.
\end{lemma}

\subsection{Proof of the summability theorem}\label{sec:summability-proof}
\begin{proof}[Proof of Theorem \ref{thm:derivatives}]\label{proof:summability}
  Fix $\kappa>0$ and first take  $2/3<\rho<1/\tau$.
  Then $0<\rho<1$. Let $C_{\rm maj}$ be the constant in
Proposition~\ref{prop:majorant}. Set
  \[
 A_0=\kappa C_{\rm maj}^\rho,\qquad
 {w_{\mathcal E}}(E)=A_0^{|E|}(|E|!)^{\tau\rho}
       \prod_{j\in E}(\sqrt T\,\sigma_j)^\rho,
 \qquad
 w(S)=A_0^{|S|}(|S|!)^{\tau\rho}\eta_m(S)^\rho
\] 
  for $E\subseteq \mathcal E_w$ and nonempty $S\subseteq[s]$.
By \eqref{eq:column-monotonicity}, $\Theta_\rho\le \Theta_{2/3}^{3\rho/2}$.
Using \eqref{eq:endpoint-subset-sum-definition},
\eqref{eq:all-subset-sum-definition}, and the monotonicity of the constant
in Lemma~\ref{lem:all-subset-sum}, we have
\begin{equation}\label{eq:proof-endpoint-weight-sum}
 \sum_{E\subseteq\mathcal E_w}{w_{\mathcal E}}(E)={\Upsilon}^{\mathcal E}(A_0)
 \le C_{{\Upsilon}}(A_0,\rho,\tau,T;\Theta_\rho)
 \le C_{{\Upsilon}}(A_0,\rho,\tau,T;\Theta_{2/3}^{3\rho/2}),
\end{equation}
\begin{equation}\label{eq:proof-block-weight-sum}
 \sum_{\varnothing\ne S\subseteq[s]}w(S)={\Upsilon}_m(A_0)
 \le C_{{\Upsilon}}(A_0,\rho,\tau,T;\Theta_\rho)
 \le C_{{\Upsilon}}(A_0,\rho,\tau,T;\Theta_{2/3}^{3\rho/2}).
\end{equation}
Using $(\sum_i a_i)^\rho\le\sum_i a_i^\rho$ for $a_i\ge0$
and $0<\rho<1$ in \eqref{eq:majorant}, together with
$|u|=|E|+\sum_{S\in\pi}|S|$, gives
\[
 \kappa^{|u|}\mathcal D_{m,u}^\rho
 \le\sum_{E\subseteq u\cap \mathcal E_w}{w_{\mathcal E}}(E)
       \sum_{\pi\in\PP(u\setminus E)}
       C_{\rm maj}^{\rho|\pi|}(|\pi|!)^{\tau\rho}\prod_{S\in\pi}w(S).
\]
Summing over $u$ and setting $v=u\setminus E$ yields
\begin{align*}
 \sum_{\varnothing\ne u\subseteq[s]}\kappa^{|u|}\mathcal D_{m,u}^\rho
 &\le\sum_{\varnothing\ne u\subseteq[s]}
       \sum_{E\subseteq u\cap\mathcal E_w}{w_{\mathcal E}}(E)
       \sum_{\pi\in\PP(u\setminus E)}
       C_{\rm maj}^{\rho|\pi|}(|\pi|!)^{\tau\rho}\prod_{S\in\pi}w(S)\\
 &=\sum_{E\subseteq \mathcal E_w}{w_{\mathcal E}}(E)
       \sum_{\substack{v\subseteq[s]\setminus E\\E\cup v\ne\varnothing}}
       \sum_{\pi\in\PP(v)}
       C_{\rm maj}^{\rho|\pi|}(|\pi|!)^{\tau\rho}\prod_{S\in\pi}w(S)\\
 &\le\sum_{E\subseteq \mathcal E_w}{w_{\mathcal E}}(E)
       \sum_{v\subseteq[s]}\sum_{\pi\in\PP(v)}
       C_{\rm maj}^{\rho|\pi|}(|\pi|!)^{\tau\rho}\prod_{S\in\pi}w(S).
\end{align*}
The equality reindexes the finite sum using the bijection
$(u,E)\leftrightarrow(E,v)$, with $v=u\setminus E$ and
$u=E\mathbin{\dot\cup}v$. The last inequality drops the restriction
$v\cap E=\varnothing$ and  $E\cup v\ne\varnothing$.
For ${{n}}\ge0$, ordering the ${{n}}$ distinct blocks of a partition and using
\eqref{eq:proof-block-weight-sum} gives
\begin{align*}
 \sum_{v\subseteq[s]}\sum_{\substack{\pi\in\PP(v)\\|\pi|={{n}}}}
                      \prod_{S\in\pi}w(S)
 &=\frac{1}{{{n}}!}
   \sum_{\substack{S_1,\ldots,S_{{n}}\subseteq[s],\ S_i\ne\varnothing\\
                    S_i\cap S_j=\varnothing\ (i\ne j)}}
                    \prod_{j=1}^{{n}}w(S_j)\le\frac{1}{{{n}}!}\left(\sum_{\varnothing\ne S\subseteq[s]}w(S)\right)^{{n}}\\
 &=\frac{{\Upsilon}_m(A_0)^{{n}}}{{{n}}!}
 \le\frac{C_{{\Upsilon}}(A_0,\rho,\tau,T;\Theta_{2/3}^{3\rho/2})^{{n}}}{{{n}}!}.
\end{align*}
Using the preceding inequality and \eqref{eq:proof-endpoint-weight-sum}, we obtain
\begin{align}
 \sum_{\varnothing\ne u\subseteq[s]}\kappa^{|u|}\mathcal D_{m,u}^\rho
 &\le\left(\sum_{E\subseteq\mathcal E_w}{w_{\mathcal E}}(E)\right)
 \sum_{{{n}}=0}^{\infty}C_{\rm maj}^{\rho {{n}}}({{n}}!)^{\tau\rho}
 \sum_{v\subseteq[s]}\sum_{\substack{\pi\in\PP(v)\\|\pi|={{n}}}}
 \prod_{S\in\pi}w(S)\notag\\
 &\le C_{{\Upsilon}}(A_0,\rho,\tau,T;\Theta_{2/3}^{3\rho/2})\sum_{{{n}}=0}^{\infty}
 \bigl(C_{\rm maj}^\rho
 C_{{\Upsilon}}(A_0,\rho,\tau,T;\Theta_{2/3}^{3\rho/2})\bigr)^{{n}}
 ({{n}}!)^{\tau\rho-1}\notag\\
 &=: C_{\rm sum}(\kappa,\rho;\Theta_{2/3})<\infty.
 \label{eq:rho-derivative-sum-bound}
\end{align}
The ratio of consecutive terms in $C_{\rm sum}(\kappa,\rho;\Theta_{2/3})$ is
\[
 C_{\rm maj}^\rho C_{{\Upsilon}}(A_0,\rho,\tau,T;\Theta_{2/3}^{3\rho/2})
 ({{n}}+1)^{\tau\rho-1}
 \longrightarrow0,
\]
so $C_{\rm sum}(\kappa,\rho;\Theta_{2/3})$ is finite.
Since $A_0=\kappa C_{\rm maj}^\rho$ and $C_{\rm maj}$ depends only on
$T,B_0,B_1,G_0,G_1,\tau$, the constant
$C_{\rm sum}(\kappa,\rho;\Theta_{2/3})$ depends only on
$\kappa,\rho,\tau,T,B_0,B_1,G_0,G_1,\Theta_{2/3}$
and is nondecreasing in both $\kappa$ and $\Theta_{2/3}$.

For general $p>2/3$, take
\[
 \rho=\frac{2/3+\min\{p,1/\tau\}}{2}.
\]
Since $\tau<3/2$, this gives $2/3<\rho<\min\{p,1/\tau\}$. Put
 $x_u=\kappa^{|u|/p}\mathcal D_{m,u}$. Since $p>\rho$, we have
\[
 \sum_{u\ne\varnothing}\kappa^{|u|}\mathcal D_{m,u}^p
 =\sum_{u\ne\varnothing}x_u^p
 \le\left(\sum_{u\ne\varnothing}x_u^\rho\right)^{p/\rho}
 =\left(\sum_{u\ne\varnothing}(\kappa^{\rho/p})^{|u|}
                    \mathcal D_{m,u}^\rho\right)^{p/\rho}.
\]
Applying \eqref{eq:rho-derivative-sum-bound} with geometric factor $\kappa^{\rho/p}$
and taking the supremum over $m$ gives
\[
 \sup_{m\in\mathcal M}\sum_{u\ne\varnothing}
 \kappa^{|u|}\mathcal D_{m,u}^p
 \le \left[C_{\rm sum}\!\left(\kappa^{\rho/p},\rho;\Theta_{2/3}\right)\right]^{p/\rho}
 =:C_{\mathcal D}(\kappa,p;\Theta_{2/3}).
\]
For later use, we denote this bound by
$C_{\mathcal D}(\kappa,p;\Theta_{2/3})$, suppressing the other fixed parameters.
This proves \eqref{eq:derivative-sum}. The constant depends only on
$\kappa,p,\tau,T,B_0,B_1,G_0,G_1,\Theta_{2/3}$ and is nondecreasing in
both $\kappa$ and $\Theta_{2/3}$.
\end{proof}

\section{Lattice integration error}\label{sec:lattice}
In this section, we use Theorem~\ref{thm:derivatives} to bound the
weighted Gaussian-space norm introduced in Section~\ref{sec:spaces-rules}.
We then prove the lattice integration error bounds in
Theorem~\ref{thm:main} and Corollary~\ref{cor:dimensions}.

\subsection{Bounds on the weighted Gaussian-space norm}
Recall the space $\mathcal H_{s,\bm\gamma}$ from
Section~\ref{sec:spaces-rules}, with $\psi(t)=e^{-\alpha|t|}$.
For $u\subseteq[s]$, define
\begin{equation}\label{eq:B}
 {\mathcal S}_u(F)^2:=\int_{\R^{|u|}}
 \left|\int_{\R^{s-|u|}}\partial_uF(\bm z)
 \phi_{-u}(\bm z_{-u})\dd\bm z_{-u}\right|^2
 \prod_{j\in u}e^{-2\alpha|z_j|}\dd\bm z_u.
\end{equation}
Then
\[
 \|F\|_{\mathcal H_{s,\bm\gamma}}^2
 =\sum_{u\subseteq[s]}\frac{{\mathcal S}_u(F)^2}{\gamma_u}.
\]
In view of the exponential growth of $g$ in \eqref{eq:payoff-ass},
we set the parameter $\alpha$ in the weight function $\psi$ to be
\begin{equation}\label{eq:alpha-choice}
 \alpha=\alpha_0:=1+\beta_g e^{\mathsf B_1T}\sqrt T\,\Theta_1
 \left(2+\frac{1}{\sqrt2}\right).
\end{equation}
With this choice, we first prove the following technical lemma.

\begin{lemma}\label{lem:hybrid}
Suppose Assumption~\ref{ass:space} holds.
For each $u\subseteq[s]$, let
${\bm Z^{(u)}}=({Z_1^{(u)}},\ldots,{Z_s^{(u)}})$
have independent coordinates, standard Gaussian for
$j\notin u$ and with density $\alpha_0e^{-2\alpha_0|z|}$ for $j\in u$.
Then
\begin{align*}
 \sup_{m\in\mathcal M}\sup_{u\subseteq[s]}
 \E e^{2\beta_g|\bm X^h_{t_m}({\bm Z^{(u)}})|_2}
 \le C_{\rm aux}(|\bm x_0|_2,\Theta_1)<\infty.
\end{align*}
Here $C_{\rm aux}$ depends only on
$\beta_g,T,K_0,\mathsf B_1,|\bm x_0|_2,\Theta_1$
and is nondecreasing in $|\bm x_0|_2$ and $\Theta_1$.
\end{lemma}
The proof is given in Appendix~\ref{sec:auxiliary-moment-proof}.

Using Lemma~\ref{lem:hybrid}, we obtain the following estimate for ${\mathcal S}_u(F_m)$.

\begin{lemma}\label{lem:transfer}
Under Assumption~\ref{ass:space}, we have
\begin{equation}\label{eq:Bbound}
 {\mathcal S}_u(F_m)\le C_{\rm aux}(|\bm x_0|_2,\Theta_1)^{1/2}\mathcal D_{m,u},
 \qquad u\ne\varnothing,
\end{equation}
where $C_{\rm aux}$ is the constant in Lemma~\ref{lem:hybrid}.
\end{lemma}
The proof is given in Appendix~\ref{sec:auxiliary-moment-proof}.

\subsection{A fractional-moment bound for the lattice estimator}
The following proposition gives a fractional-moment bound for the error
of the random lattice estimator.
\begin{proposition}\label{lem:random}
Fix $1/2<\lambda<1$, $1/(2\lambda)<{\delta_*}<1$, and use
$\rho_\lambda$ from \eqref{eq:rho}.
{Let $F:\R^s\to\R$ have the weak mixed derivatives appearing in \eqref{eq:B}, with $\mathcal S_u(F)<\infty$ for every $u\subseteq[s]$,} and let $N$ be prime. Define
\begin{equation}\label{eq:wce-scale}
 \mathcal C_\lambda(F):=
 \left(\sum_{\varnothing\ne u\subseteq[s]}
       \rho_\lambda^{|u|/(1+\lambda)}
       {\mathcal S}_u(F)^{2\lambda/(1+\lambda)}\right)^{1+\lambda}.
\end{equation}
Let $\bm A$ be uniform on $\{1,\ldots,N-1\}^s$ and
$\bm\Delta$ uniform on $[0,1]^s$, independently of each other. Then
\begin{equation}\label{eq:singlemoment}
 \E_{\bm A,\bm\Delta}|Q_{\bm A,\bm\Delta}(F)-\mathcal I (F)|^{2\lambda}
 \le\frac{\mathcal C_\lambda(F)}{N-1}.
\end{equation}
\end{proposition}
\begin{proof}
Let $F^\circ:=F-\mathcal I (F)$ and choose 
\begin{equation}\label{eq:auxiliary-weights}
 \gamma_u=\left(\frac{{\mathcal S}_u(F)^2}{\rho_\lambda^{|u|}}
          \right)^{1/(1+\lambda)},\qquad \varnothing\ne u\subseteq[s].
\end{equation}
With the convention $0/0=0$ from Section~\ref{sec:spaces-rules},
\[
 \frac{{\mathcal S}_u(F)^2}{\gamma_u}
 =\gamma_u^\lambda\rho_\lambda^{|u|}
 =\rho_\lambda^{|u|/(1+\lambda)}{\mathcal S}_u(F)^{2\lambda/(1+\lambda)}.
\]
Since ${\mathcal S}_\varnothing(F^\circ)=0$ and
${\mathcal S}_u(F^\circ)={\mathcal S}_u(F)$ for $u\ne\varnothing$,
we have $\|F^\circ\|_{\mathcal H_{s,\bm\gamma}}^2
=\sum_{u\ne\varnothing}{\mathcal S}_u(F)^2/\gamma_u$.
Then \eqref{eq:shift-averaged-wce} gives
\begin{equation*}
 \E_{\bm\Delta}|Q_{\bm a,\bm\Delta}(F)-\mathcal I(F)|^2= \E_{\bm\Delta}|Q_{\bm a,\bm\Delta}(F^\circ)-\mathcal I(F^\circ)|^2
 \le\left(\sum_{u\ne\varnothing}\frac{{\mathcal S}_u(F)^2}{\gamma_u}\right)
       [e^{\rm sh}_{s,N}(\bm a;\bm\gamma)]^2.
\end{equation*}
Raising this bound to the power $\lambda$, averaging over $\bm A$,
and applying \eqref{eq:random-wce-moment} yields
\begin{align*}
 \E_{\bm A}\left(\E_{\bm\Delta}
 |Q_{\bm A,\bm\Delta}(F)-\mathcal I(F)|^2\right)^\lambda
 &\le\left(\sum_{u\ne\varnothing}\frac{{\mathcal S}_u(F)^2}{\gamma_u}\right)^\lambda
 \E_{\bm A}\left[e^{\rm sh}_{s,N}(\bm A;\bm\gamma)^{2\lambda}\right]\\
 &\le\frac{1}{N-1}
 \left(\sum_{u\ne\varnothing}\frac{{\mathcal S}_u(F)^2}{\gamma_u}\right)^\lambda
 \sum_{u\ne\varnothing}\gamma_u^\lambda\rho_\lambda^{|u|}\\
 &=\frac{1}{N-1}
 \left(\sum_{u\ne\varnothing}\rho_\lambda^{|u|/(1+\lambda)}
 {\mathcal S}_u(F)^{2\lambda/(1+\lambda)}\right)^{1+\lambda}
 =\frac{\mathcal C_\lambda(F)}{N-1}.
\end{align*}
Finally, since $0<\lambda<1$, conditional Jensen's inequality gives
\[
 \E_{\bm A,\bm\Delta}|Q_{\bm A,\bm\Delta}(F)-\mathcal I(F)|^{2\lambda}
 \le\E_{\bm A}\left(\E_{\bm\Delta}
 |Q_{\bm A,\bm\Delta}(F)-\mathcal I(F)|^2\right)^\lambda,
\]
which proves \eqref{eq:singlemoment}.
\end{proof}

\subsection{Proof of the sampling error bound}\label{sec:sampling-proof}
To prove Theorem~\ref{thm:main}, we first establish a moment bound
for the median of independent copies of an error variable.
\begin{lemma}\label{lem:median}
Let $E_1,\ldots,E_R$ be independent identically distributed real random variables,
$R=2k-1$, and suppose $\E|E_1|^{{q_0}}\le K^{{q_0}}$ for ${{q_0}}>0$ and $K\ge0$.
For every $0<q<{{q_0}} k$,
\begin{equation}\label{eq:medianmoment}
 \E|\operatorname{median}(E_1,\ldots,E_R)|^q
 \le\frac{{{q_0}} k}{{{q_0}} k-q}
       \binom Rk^{q/({{q_0}} k)}K^q.
\end{equation}
\end{lemma}
\begin{proof}
Write $M_R:=\operatorname{median}(E_1,\ldots,E_R)$. If $K=0$, then
$E_{{\nu}}=0$ almost surely for every ${\nu}$. Suppose $K>0$. If $|M_R|>t$, then
$|E_{{\nu}}|>t$ for at least $k$ indices. Therefore, by a union bound over the
$k$-element subsets ${\mathcal A}\subseteq\{1,\ldots,R\}$ and independence,
\begin{align*}
 \Prob(|M_R|>t)
 \le\sum_{|{\mathcal A}|=k}\Prob\!\left(\bigcap_{{\nu\in\mathcal A}}\{|E_{{\nu}}|>t\}\right) =\binom Rk\Prob(|E_1|>t)^k
 \le\binom Rk\left(\frac{K}{t}\right)^{{{q_0}} k},
\end{align*}
where the last inequality is Markov's inequality. Hence
\[
 \Prob(|M_R|>t)
 \le\min\{1,\tbinom Rk(K/t)^{{{q_0}} k}\}.
\]
Set $t_*=\binom Rk^{1/({{q_0}} k)}K$. Since $q<{{q_0}} k$, the tail integral
formula and the preceding bound give
\begin{align*}
 \E|M_R|^q
 &=q\int_0^\infty t^{q-1}\Prob(|M_R|>t)\dd t\le q\int_0^{t_*}t^{q-1}\dd t
   +q\binom Rk K^{{{q_0}} k}
       \int_{t_*}^\infty t^{q-1-{{q_0}} k}\dd t\\
 &=t_*^q+\frac{q}{{{q_0}} k-q}
       \binom Rk K^{{{q_0}} k}t_*^{q-{{q_0}} k}
 =\frac{{{q_0}} k}{{{q_0}} k-q}t_*^q.
\end{align*}
Inserting the value of $t_*$ proves \eqref{eq:medianmoment}.
\end{proof}

\begin{proof}[Proof of Theorem~\ref{thm:main}]\phantomsection\label{proof:sampling}
Fix $0<\varepsilon<1/2$.
Set $\alpha=\alpha_0$ as in \eqref{eq:alpha-choice}. In
Proposition~\ref{lem:random}, take $\lambda=1/[2(1-\varepsilon)]$ and
${\delta_*}=1-\varepsilon/2$; then $2{\delta_*}\lambda>1$.
  For this value of $\lambda$, we have
$2\lambda/(1+\lambda)=2/(3-2\varepsilon)>2/3$.
By \eqref{eq:Bbound} and \eqref{eq:derivative-sum},
\begin{align*}
 &\sup_{m\in\mathcal M}\sum_{u\ne\varnothing}
 \rho_\lambda^{|u|/(1+\lambda)}
 {\mathcal S}_u(F_m)^{2\lambda/(1+\lambda)}\\
 &\quad\le C_{\rm aux}(|\bm x_0|_2,\Theta_1)^{\lambda/(1+\lambda)}
 \sup_{m\in\mathcal M}\sum_{u\ne\varnothing}
 \bigl(\rho_\lambda^{1/(1+\lambda)}\bigr)^{|u|}
 \mathcal D_{m,u}^{2\lambda/(1+\lambda)}\\
 &\quad\le C_{\rm aux}(|\bm x_0|_2,\Theta_1)^{\lambda/(1+\lambda)}
 C_{\mathcal D}\!\left(
              \rho_\lambda^{1/(1+\lambda)},
              \frac{2\lambda}{1+\lambda};\Theta_{2/3}\right).
\end{align*}
Here $C_{\rm aux}$ and $C_{\mathcal D}$ are the constants in
Lemma~\ref{lem:hybrid} and Theorem~\ref{thm:derivatives}, respectively.
By the definition of $\mathcal C_\lambda$ in \eqref{eq:wce-scale}, we obtain
\begin{equation}\label{eq:uniform-wce-moment}
 \begin{aligned}
 \sup_{m\in\mathcal M}\mathcal C_\lambda(F_m)
 &\le C_{\rm aux}(|\bm x_0|_2,\Theta_1)^\lambda
 \left[C_{\mathcal D}\!\left(
              \rho_\lambda^{1/(1+\lambda)},
              \frac{2\lambda}{1+\lambda};\Theta_{2/3}\right)\right]^{1+\lambda}\\
 &=:C_{\rm lat}(|\bm x_0|_2,\Theta_1,\Theta_{2/3})<\infty.
 \end{aligned}
\end{equation}
Since $\alpha_0$ is defined by \eqref{eq:alpha-choice} and
$\rho_\lambda$ by \eqref{eq:rho}, this bound depends only on
$\varepsilon,\tau,T,K_0$, $B_0,B_1,G_0,G_1$, $\beta_g$,
$|\bm x_0|_2$, $\Theta_1$, and $\Theta_{2/3}$.
For fixed $\varepsilon$ and fixed model parameters, it is nondecreasing in
$|\bm x_0|_2$, $\Theta_1$, and $\Theta_{2/3}$. Indeed, $\alpha_0$ and
$\rho_\lambda$ are nondecreasing in $\Theta_1$, while $C_{\rm aux}$ and
$C_{\mathcal D}$ have the monotonicities stated above.

Substituting $F=F_m$ in \eqref{eq:singlemoment} and using
\eqref{eq:uniform-wce-moment} yields 
\begin{equation}\label{eq:uniform-single-copy-moment}
 \sup_{m\in\mathcal M}
 \E_{\bm{A},\bm{\Delta}}|Q_{\bm{A},\bm{\Delta}}^{N,m}-\mathcal I_m|^{2\lambda}
 \le\frac{C_{\rm lat}(|\bm x_0|_2,\Theta_1,\Theta_{2/3})}{N-1}.
\end{equation}
Since $2\lambda>1$, Jensen's inequality first gives
\[
 \|Q_{\bm{A},\bm{\Delta}}^{N,m}-\mathcal I_m\|_{L^1}
 \le\left(\E|Q_{\bm{A},\bm{\Delta}}^{N,m}-\mathcal I_m|^{2\lambda}
        \right)^{1/(2\lambda)}
 \le\left(\frac{C_{\rm lat}(|\bm x_0|_2,
 \Theta_1,\Theta_{2/3})}{N-1}\right)^{1/(2\lambda)}.
\]
Since $N-1\ge N/2$ and $1/(2\lambda)=1-\varepsilon$,
this proves \eqref{eq:single-l1-main} with
\begin{equation}\label{eq:single-lattice-error-constant}
 C_\varepsilon:=2^{1/(2\lambda)}
 C_{\rm lat}(|\bm x_0|_2,
 \Theta_1,\Theta_{2/3})^{1/(2\lambda)}.
\end{equation}

For $R=2k-1$ and $1\le q<2\lambda k = k/(1-\varepsilon)$, apply
Lemma~\ref{lem:median} to the independent errors
$Q_{\bm{A}_{{\nu}},\bm{\Delta}_{{\nu}}}^{N,m}-\mathcal I_m$, with ${q_0}=2\lambda$ and
$K=(C_{\rm lat}(|\bm x_0|_2,
\Theta_1,\Theta_{2/3})/(N-1))^{1/(2\lambda)}$. Then we have
\[
 \|M_{R,N,m}-\mathcal I_m\|_{L^q}
 \le\left(\frac{2\lambda k}{2\lambda k-q}\right)^{1/q}
       \binom Rk^{1/(2\lambda k)}
 \left(\frac{C_{\rm lat}(|\bm x_0|_2,
 \Theta_1,\Theta_{2/3})}{N-1}\right)^{1/(2\lambda)}.
\]
Using $N-1\ge N/2$ proves
\eqref{eq:median-sampling-main} with
\begin{equation}\label{eq:median-error-constant}
 C_{\varepsilon,q,R}:=
 \left(\frac{2\lambda k}{2\lambda k-q}\right)^{1/q}
 \binom Rk^{1/(2\lambda k)}C_\varepsilon.
\end{equation}
By \eqref{eq:uniform-wce-moment} 
and \eqref{eq:single-lattice-error-constant}, $C_\varepsilon$ depends only on
$\varepsilon,\tau,T,K_0$, $B_0,B_1,G_0,G_1$, $\beta_g$,
$|\bm x_0|_2$, $\Theta_1$, and $\Theta_{2/3}$.
By \eqref{eq:median-error-constant}, $C_{\varepsilon,q,R}$ additionally
depends on $q,R$.
\end{proof}

\begin{proof}[Proof of Corollary~\ref{cor:dimensions}]
Under Assumption~\ref{ass:uniform-data}, the model parameters have
common bounds across dimensions. Under Assumption~\ref{ass:columns},
\eqref{eq:column-monotonicity} gives
\[
 0\le\Theta_1\le M_{2/3}^{3/2},\qquad
 0\le\Theta_{2/3}\le M_{2/3}.
\]
By the monotonicity of $C_{\rm lat}$ established immediately after
\eqref{eq:uniform-wce-moment},
$|\bm x_0|_2\le R_0$ and \eqref{eq:column-monotonicity} give
\[
 C_{\rm lat}(|\bm x_0|_2,\Theta_1,\Theta_{2/3})
 \le C_{\rm lat}(R_0,M_{2/3}^{3/2},M_{2/3}).
\]
The right-hand side is independent of $d_x,d_w$. Therefore,
\eqref{eq:single-lattice-error-constant} and
\eqref{eq:median-error-constant} give dimension-independent bounds for
$C_\varepsilon$ and $C_{\varepsilon,q,R}$.
\end{proof}

\section{Numerical experiments}\label{sec:numerics}

In this section, we use numerical examples to examine {the sampling error of randomly shifted rank-1 lattice rules and of their median versions}
applied to the Euler--Maruyama scheme and to test the conclusions of
Theorem~\ref{thm:main} and Corollary~\ref{cor:dimensions}.
We focus on {their} advantage over standard Monte Carlo (MC), the stability
of the QMC error as the number of time steps $m$ increases, and, under
the stated summability condition on the columns of $\Sigma$, its stability
as the {Brownian path dimension} $d_w$ increases.
All models use $T=1$, and the target is the Euler--Maruyama expectation $\mathcal I_m$. The code and data used for these numerical experiments
are available at \url{https://github.com/ZiyangYe2000/lattice4SDE}.

We consider four estimators: MC, single QMC, median-of-three QMC, and
median-of-thirteen QMC. In all experiments, $d_x=d_w=d$.
For a fixed pair $(d,m)$, let $s=dm$ and write
$F_{d,m}(\bm z)=g(\bm X^h_{t_m}(\bm z))$ for the integrand
$F_m$ in \eqref{eq:integrand}, displaying its dependence on $d$.
For a prescribed budget of ${n}$ simulated Euler--Maruyama paths,
the MC estimator is
\[
 \widehat{{\mathcal I}}^{\rm MC}_{m,{n}}=\frac1{{n}}\sum_{\ell=1}^{{n}}F_{d,m}({\bm Z_\ell}).
\]
Here ${\bm Z_1},\ldots,{\bm Z_n}\overset{\mathrm{i.i.d.}}{\sim}N(0,{\operatorname{Id}_s})$.
For the lattice-rule estimators, each constituent rule
independently randomizes both the generating vector and the shift: $\bm A_{{\nu}}$ is sampled uniformly from
$\{1,\ldots,N-1\}^s$, while
$\bm\Delta_{{\nu}}$ {is uniform on $[0,1]^s$}.  The ${\nu}$-th randomized {rank-1} lattice
estimate, $Q_{m,N}^{({\nu})}=Q_{\bm A_{{\nu}},\bm\Delta_{{\nu}}}^{N,m}$ in \eqref{eq:Q}, is
\[
 Q_{m,N}^{({\nu})}=\frac1N\sum_{\ell=0}^{N-1}
 F_{d,m}\!\left(\Phi^{-1}\!\left(
 \left\{\frac{\ell \bm A_{{\nu}}}{N}+\bm\Delta_{{\nu}}\right\}\right)\right).
\]
Single QMC uses one such rule,
\[
 \widehat{{\mathcal I}}^{\rm SQMC}_{m,{n}}=Q_{m,N_1}^{(1)},
 \qquad N_1=\text{the nearest prime to }{n}.
\]
The other two estimators are
\[
 \widehat{{\mathcal I}}^{\operatorname{med}({R})}_{m,{n}}
 =\operatorname{median}\{Q_{m,N_{{R}}}^{(1)},\ldots,Q_{m,N_{{R}}}^{({R})}\},
 \qquad {R}\in\{3,13\},
\]
where $N_{{R}}$ is the nearest prime to ${n}/{R}$. Hence, every estimator uses
approximately ${n}$ Euler--Maruyama paths in total.
The figures use the actual number of path evaluations on the horizontal
axis: ${n}$ for MC, $N_1$ for single QMC, and ${R}N_{{R}}$ for the median rules.

For each $(d,m)$, ${\mathcal I}_m^{\rm ref}$ is a high-accuracy estimate of the
Euler--Maruyama expectation $\mathcal I_m$. Using $n_{\rm rep}=32$ independent
replications, we compute
\[
 \widehat{\operatorname{RMSE}}({n})
 =\left\{\frac1{n_{\rm rep}}\sum_{\ell=1}^{n_{\rm rep}}
 \left(\widehat{{\mathcal I}}_{m,{n}}^{(\ell)}-{\mathcal I}_m^{\rm ref}\right)^2\right\}^{1/2}.
\]

\subsection{Time-modulated asymmetric Kuramoto system}\label{sec:numerics-p2}

The first model is a noisy population of coupled phase oscillators,
\begin{equation}\label{eq:kuramoto equation}
 \begin{aligned}
 \dd X_j(t)&=\left[\omega_j+
 \left(1+\tfrac12\sin(2\pi t)\right)
 \sum_{\substack{{j'}=1\\{j'}\ne j}}^d
 K_{j{j'}}\sin\bigl(X_{{j'}}(t)-X_j(t)\bigr)\right]\dd t\\
 &\quad+0.4\,\dd W_j(t),\qquad j=1,\ldots,d,
 \end{aligned}
\end{equation}
where
\[
 \begin{gathered}
 \omega_j=-0.5+\frac{j-1}{d-1},\qquad
 K_{j{j'}}=\frac{1+0.15\sin(j({j'}+1))}{d}\quad(j\ne {j'}),
 \\
 K_{jj}=0,
 \qquad X_j(0)=\frac{2\pi(j-1)}d.
 \end{gathered}
\]
This is a stochastic Kuramoto-type system with heterogeneous frequencies and
asymmetric coupling \cite{acebron2005kuramoto}.  The factor
$1+\tfrac12\sin(2\pi t)$ periodically modulates the common interaction
strength and tests the time-dependent setting of Theorem~\ref{thm:main}.
We take the payoff function
\[
 g(\bm x)=\left|\frac1d\sum_{j=1}^d e^{\mathrm i x_j}\right|^2,
 \qquad \mathrm i^2=-1,
\]
which is the squared Kuramoto order parameter.  It is close to one for coherent
phases and close to zero for a dispersed population.
It is straightforward to verify Assumption~\ref{ass:space} with $\tau=1$
for each fixed $d$. Here $\Sigma=0.4{\operatorname{Id}_d}$, so
$\Theta_{2/3}=0.4^{2/3}d$. {Hence $\Theta_{2/3}$ grows linearly in $d$, so Assumption~\ref{ass:columns} does not hold uniformly in $d$, and any stability in $d$ observed for this model is not covered by Corollary~\ref{cor:dimensions}. Since the factor $1+\tfrac12\sin(2\pi t)$ is Lipschitz continuous in $t$, Theorem~\ref{thm:euler-bias} gives a time discretization error of order $m^{-1}$.}

\begin{figure}[H]
 \centering
 \includegraphics[width=0.7\linewidth]{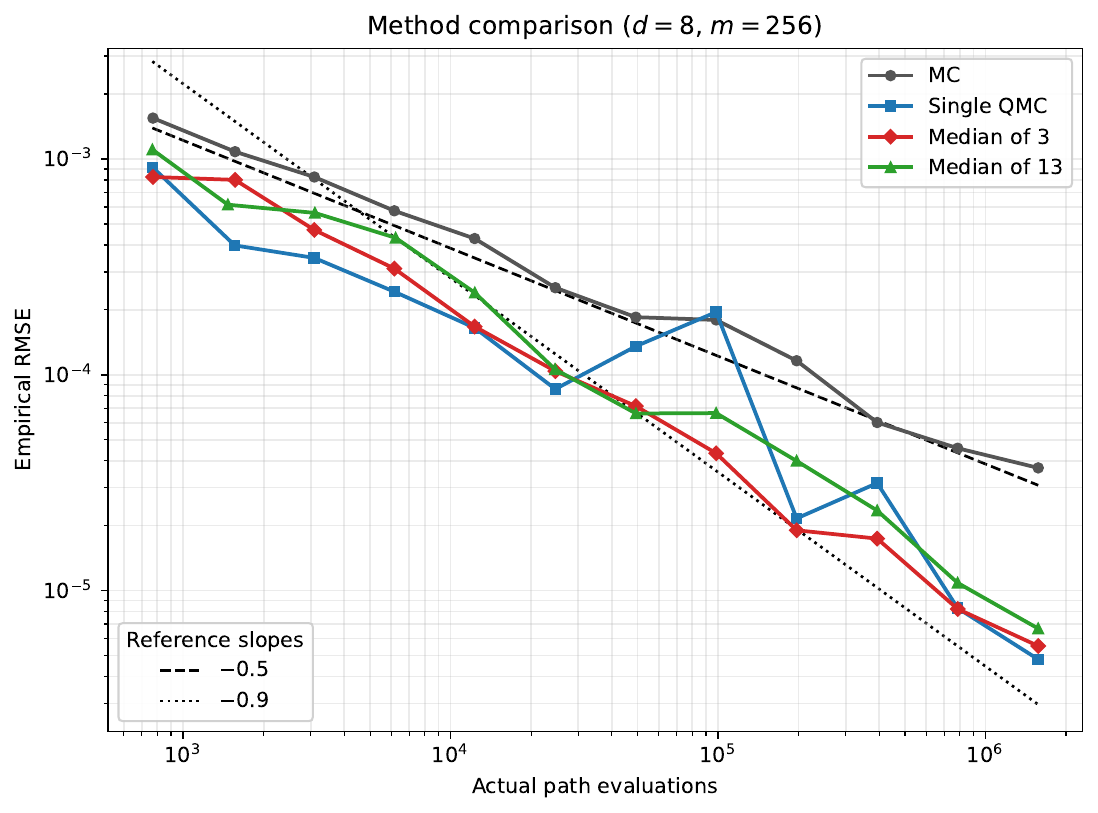}
 \caption{Method comparison for the time-modulated asymmetric Kuramoto
 system \eqref{eq:kuramoto equation} with $d=8$ and $m=256$.}
 \label{fig:p2-methods}
\end{figure}

Figure~\ref{fig:p2-methods} compares all four estimators at
$(d,m)=(8,256)$.  At comparable values of ${n}$, QMC generally improves upon
MC.  Although single QMC sometimes attains the smallest empirical RMSE, its
curve contains pronounced upward fluctuations, reflecting occasional poor
estimates from an individual randomized rule.  Median-of-three and
median-of-thirteen suppress this instability and have a smoother error
decay.  At the largest budgets, the median-QMC errors exhibit an empirical
rate close to ${n}^{-0.9}$, substantially steeper than the MC reference rate
${n}^{-0.5}$.

\begin{figure}[H]
 \centering
 \begin{minipage}[t]{0.49\linewidth}
  \centering
  \includegraphics[width=\linewidth]{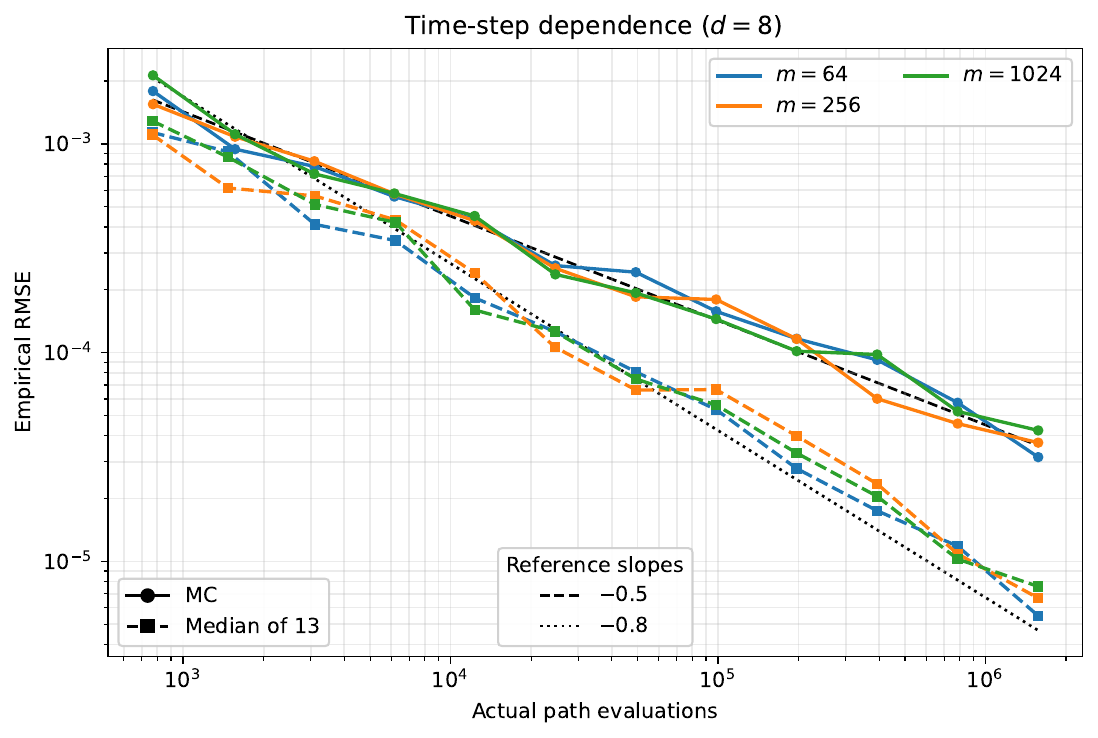}
  \par\smallskip\textbf{(a)} Fixed $d=8$; varying $m$.
 \end{minipage}\hfill
 \begin{minipage}[t]{0.49\linewidth}
  \centering
  \includegraphics[width=\linewidth]{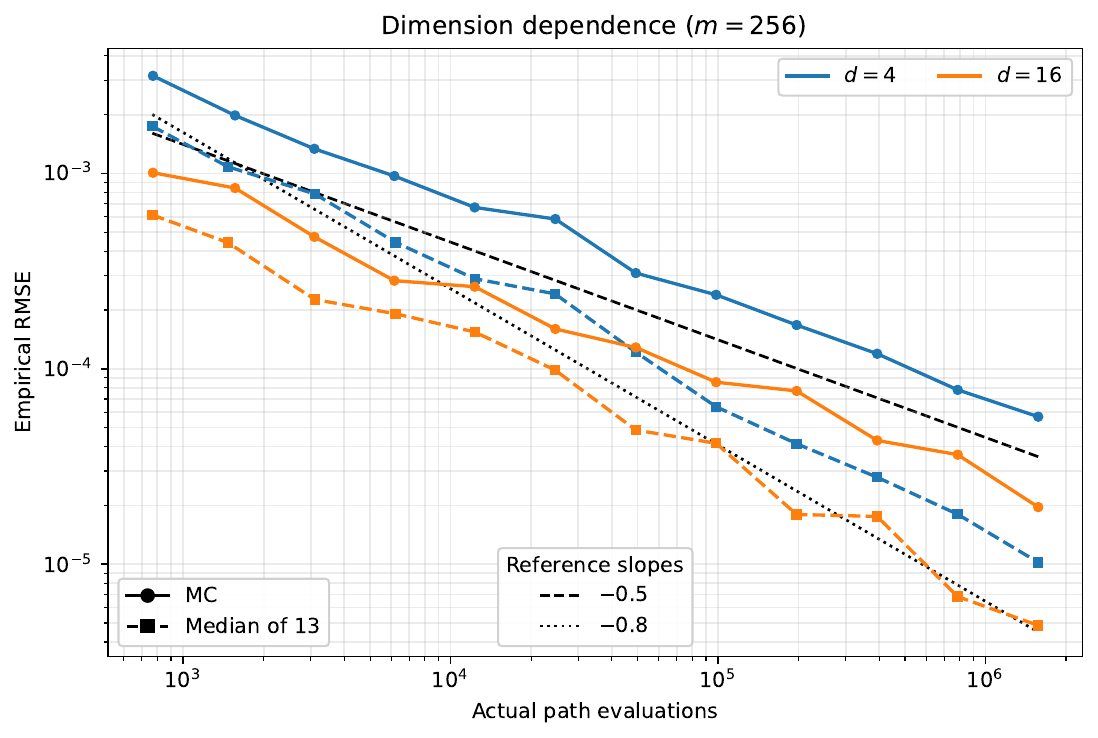}
  \par\smallskip\textbf{(b)} Fixed $m=256$; varying $d$.
 \end{minipage}
 \caption{Dependence on the number of time steps and the state dimension for
 the time-modulated asymmetric Kuramoto system \eqref{eq:kuramoto equation}.}
 \label{fig:p2-parameters}
\end{figure}

Figure~\ref{fig:p2-parameters} examines the effects of time refinement and
state dimension.  In both panels, median-of-thirteen has smaller RMSE than MC
at comparable values of ${n}$, and its large-budget decay is close to the
${n}^{-0.8}$ reference rate.  In panel~(a), we fix $d=8$ and vary
$m\in\{64,256,1024\}$.  For both MC and median-of-thirteen, the three RMSE
curves remain close as $m$ increases, with no visible loss of accuracy under
time refinement.  In panel~(b), we fix $m=256$
and compare $d=4$ with $d=16$.  The RMSE is smaller at $d=16$ for both
methods.  The decrease may be related to the $1/d$ normalization in the payoff.

\subsection{Time-dependent overdamped Langevin system}\label{sec:numerics-p5}

The second model is a system of interacting particles governed
by time-dependent overdamped Langevin dynamics \cite{pavliotis2014stochastic},
\begin{equation}\label{eq:langevin}
 \begin{aligned}
 \dd X_j(t)&=\left[-X_j(t)-{a}(t)
 \left\{\sin X_j(t)+\frac{{0.5}}{d}
 \sum_{{j'}\ne j}\sin\bigl(X_j(t)-X_{{j'}}(t)\bigr)\right\}\right]\dd t\\
 &\quad+{0.4}\,\dd W_j(t),\qquad j=1,\ldots,d,
 \end{aligned}
\end{equation}
with
\[
 {a}(t)=1+\frac12|2t-1|^{1/2},\qquad
 X_j(0)=0.2\cos\left(\frac{2\pi(j-1)}d\right).
\]
Its drift is $-\nabla_xV(t,\bm x)$ for the time-dependent potential
\[
 V(t,\bm x)={\frac12}|\bm x|_2^2+{a}(t)
 \left\{\sum_{j=1}^d(1-\cos x_j)
 +\frac{{0.5}}{d}\sum_{j<{j'}}(1-\cos(x_j-x_{{j'}}))\right\}.
\]
We take the payoff function
\[
 g(\bm x)=\exp\left(\frac{0.2}{\sqrt d}\sum_{j=1}^d x_j\right).
\]
This is the exponential of a normalized collective coordinate; the factor
$d^{-1/2}$ keeps its scale comparable as $d$ varies.
Assumption~\ref{ass:space} is easily verified with $\tau=1$ for each fixed
$d$. As in the first model, $\Theta_{2/3}=0.4^{2/3}d$. {Hence Assumption~\ref{ass:columns} again does not hold uniformly in $d$. The modulation $a$ is $1/2$-H\"older continuous but not Lipschitz continuous at $t=1/2$, so Theorem~\ref{thm:euler-strong-half} applies and gives a time discretization error of order $m^{-1/2}$, whereas Theorem~\ref{thm:euler-bias} does not apply.}

\begin{figure}[H]
 \centering
 \includegraphics[width=0.63\linewidth]{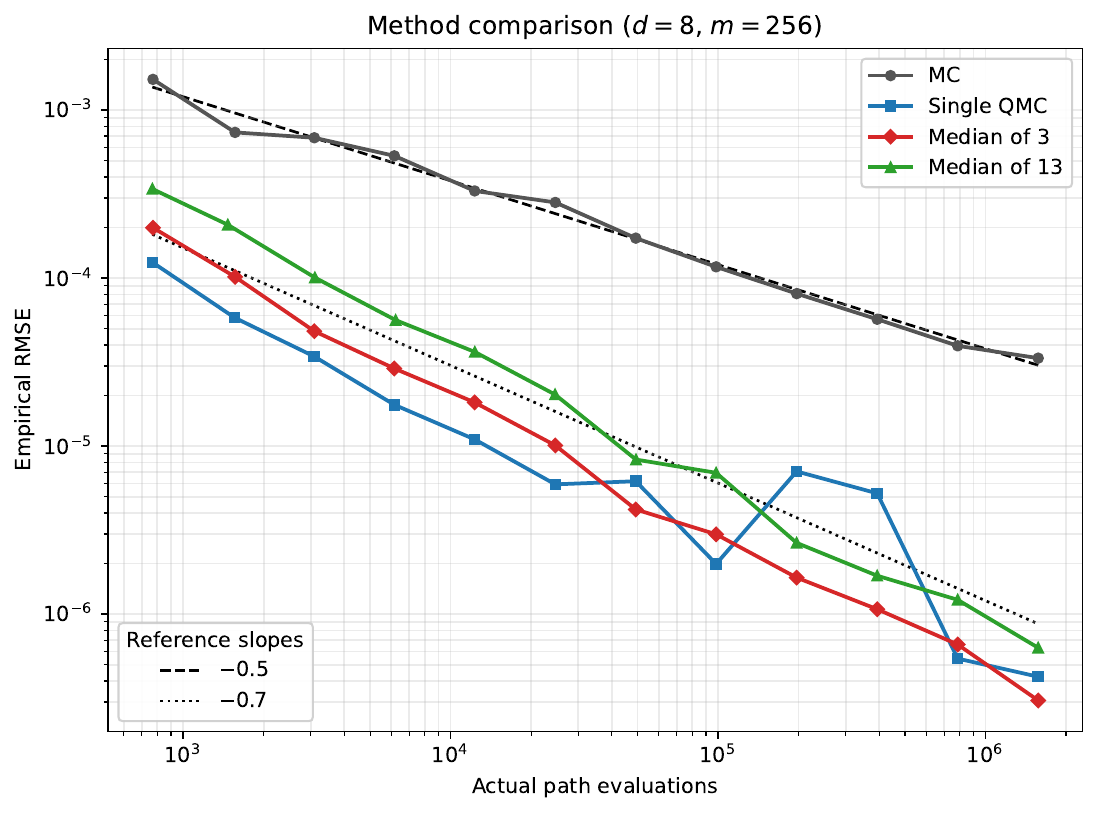}
 \caption{Method comparison for the time-dependent overdamped Langevin
 system \eqref{eq:langevin} with $d=8$ and $m=256$.}
 \label{fig:p5-methods}
\end{figure}

Figure~\ref{fig:p5-methods} compares the four estimators at
$(d,m)=(8,256)$.  As in the previous example, all three QMC estimators have
substantially smaller RMSE than MC.  Single QMC has occasional upward jumps,
whereas the two median estimators decay more steadily.  At large budgets,
median-of-thirteen has an empirical rate close to ${n}^{-0.7}$.

\begin{figure}[H]
 \centering
 \begin{minipage}[t]{0.49\linewidth}
  \centering
  \includegraphics[width=\linewidth]{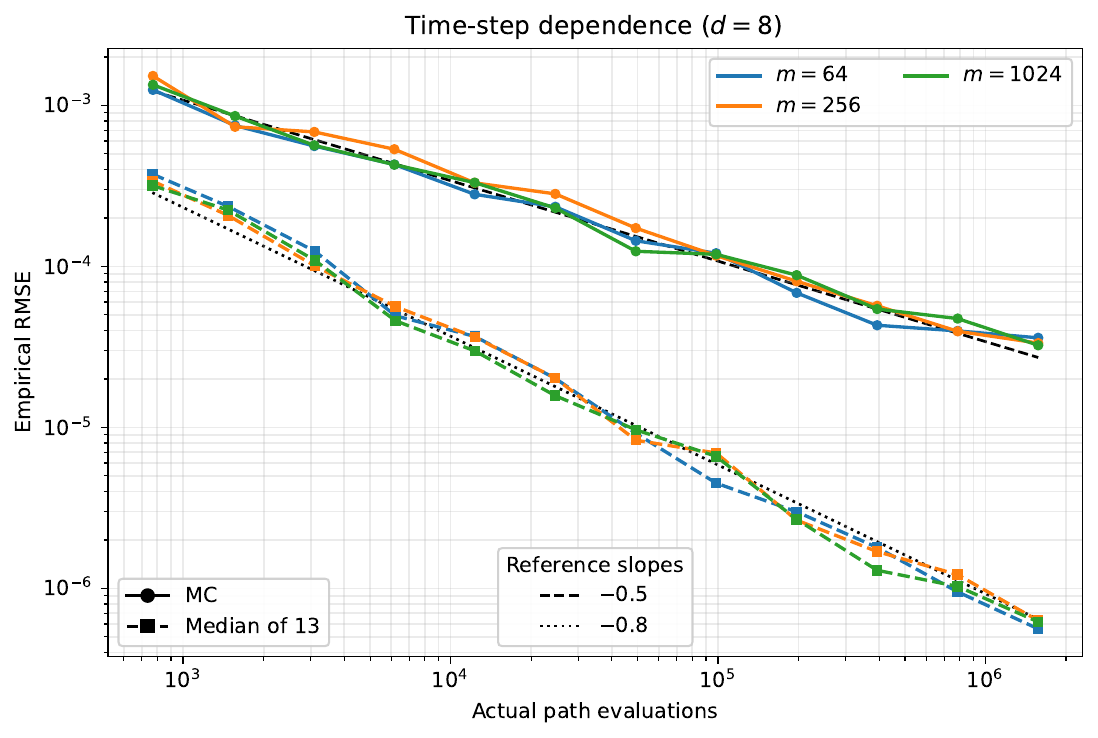}
  \par\smallskip\textbf{(a)} Fixed $d=8$; varying $m$.
 \end{minipage}\hfill
 \begin{minipage}[t]{0.49\linewidth}
  \centering
  \includegraphics[width=\linewidth]{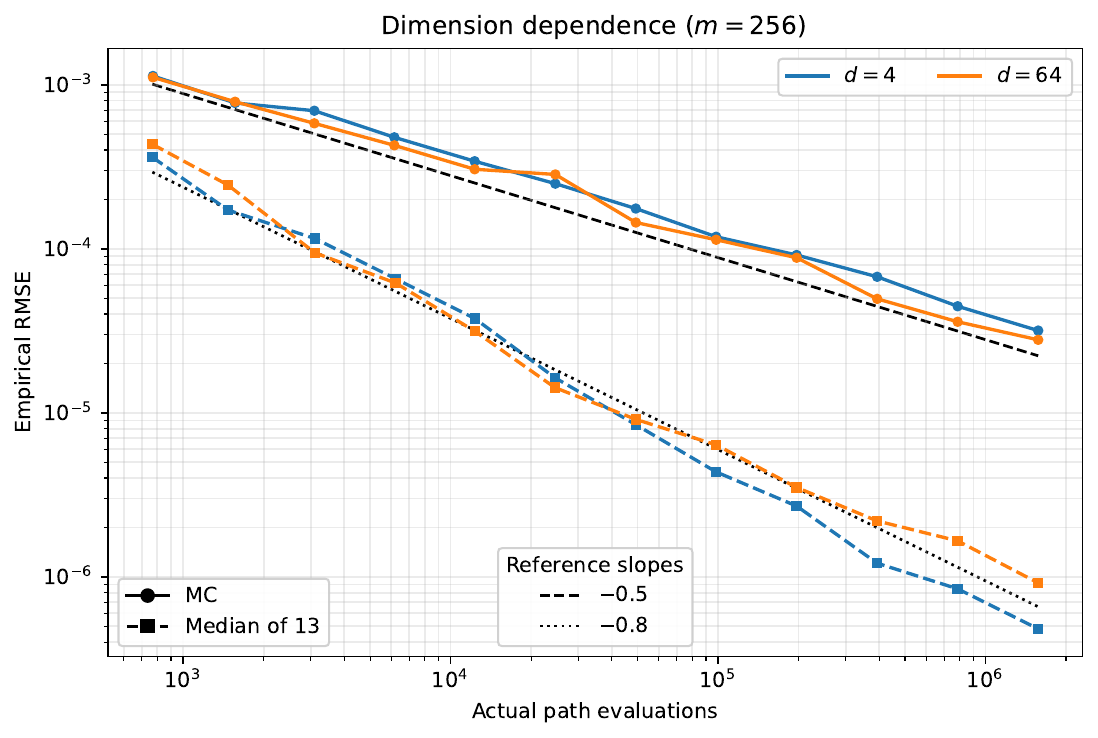}
  \par\smallskip\textbf{(b)} Fixed $m=256$; varying $d$.
 \end{minipage}
 \caption{Dependence on the number of time steps and the state dimension for
 the time-dependent overdamped Langevin system \eqref{eq:langevin}.}
 \label{fig:p5-parameters}
\end{figure}

Figure~\ref{fig:p5-parameters} shows the dependence on $m$ and $d$.  In both
panels, median-of-thirteen remains below MC and its decay is close to the
${n}^{-0.8}$ reference rate.  For fixed $d=8$, the RMSE curves change little as
$m$ increases from $64$ to $1024$.  For fixed $m=256$, the curves for $d=4$
and $d=64$ are also close, with no visible deterioration at the larger
dimension. {This robustness in $d$ is not covered by Corollary~\ref{cor:dimensions}, since $\Theta_{2/3}$ grows linearly in $d$ for this model; it suggests that Assumption~\ref{ass:columns} is sufficient but not necessary for dimension-stable behavior.}

\subsection{Dimension-uniform design and coordinate rotation}
\label{sec:numerics-p3-p4}

The final experiment is designed around the dimension-independent conclusion
of Corollary~\ref{cor:dimensions}.  The aligned model is
\[
 \dd \bm X_t=\left[-({\operatorname{Id}_d}+0.25L_d)\bm X_t+0.5\sin(\bm X_t)\right]\dd t
 +\Sigma_d^{\rm al}\dd \bm W_t,
 \qquad g(\bm x)=|\bm x|_2^2,
\]
where $\bm X_0=0$ and
\[
 L_d=\begin{pmatrix}
 2 & -1 & & & \\
 -1 & 2 & -1 & & \\
 & \ddots & \ddots & \ddots & \\
 & & -1 & 2 & -1\\
 & & & -1 & 2
 \end{pmatrix}\in\R^{d\times d}.
\]
The diffusion matrix is
\[
 \Sigma_d^{\rm al}=0.4\operatorname{diag}(1,2^{-2},\ldots,d^{-2}),
 \qquad
 \sum_{r=1}^d|\Sigma_d^{\rm al}\bm e_r|_2^{2/3}
 \le 0.4^{2/3}\zeta(4/3).
\]
The model satisfies Assumptions~\ref{ass:space} and
\ref{ass:uniform-data} with $\tau=1$: $\|L_d\|_2\le4$, the higher
derivative operator norms of the componentwise sine are at most one,
and the quadratic payoff has common derivative bounds.
The displayed column bound verifies Assumption~\ref{ass:columns}.
For the rotated model we set
\[
 \Sigma_d^{\rm rot}=\Sigma_d^{\rm al}H_d,
\]
where $d$ is a power of two and $H_d$ is the normalized Sylvester
Hadamard matrix defined recursively by
\[
 H_1=(1),\qquad
 H_{2d}=\frac{1}{\sqrt2}
 \begin{pmatrix}
 H_d&H_d\\
 H_d&-H_d
 \end{pmatrix}.
\]
Then $H_dH_d^\top={\operatorname{Id}_d}$, and hence
\[
 \Sigma_d^{\rm rot}(\Sigma_d^{\rm rot})^\top
 =\Sigma_d^{\rm al}(\Sigma_d^{\rm al})^\top.
\]
The two Euler--Maruyama schemes therefore have the same law.  On the other hand,
\[
 \sum_{r=1}^d|\Sigma_d^{\rm rot}\bm e_r|_2
 =0.4d^{1/2}
 \left(\sum_{j=1}^d j^{-4}\right)^{1/2}\longrightarrow\infty.
\]
Thus the column norms are not even uniformly $\ell^1$-summable in $d$, and
the summability condition is no longer satisfied as $d\to\infty$.
The QMC rule depends on the Gaussian coordinate representation, and the
rotation spreads the dominant noise directions over all coordinates.

\begin{figure}[H]
 \centering
 \begin{minipage}[t]{0.49\linewidth}
  \centering
  \includegraphics[width=\linewidth]{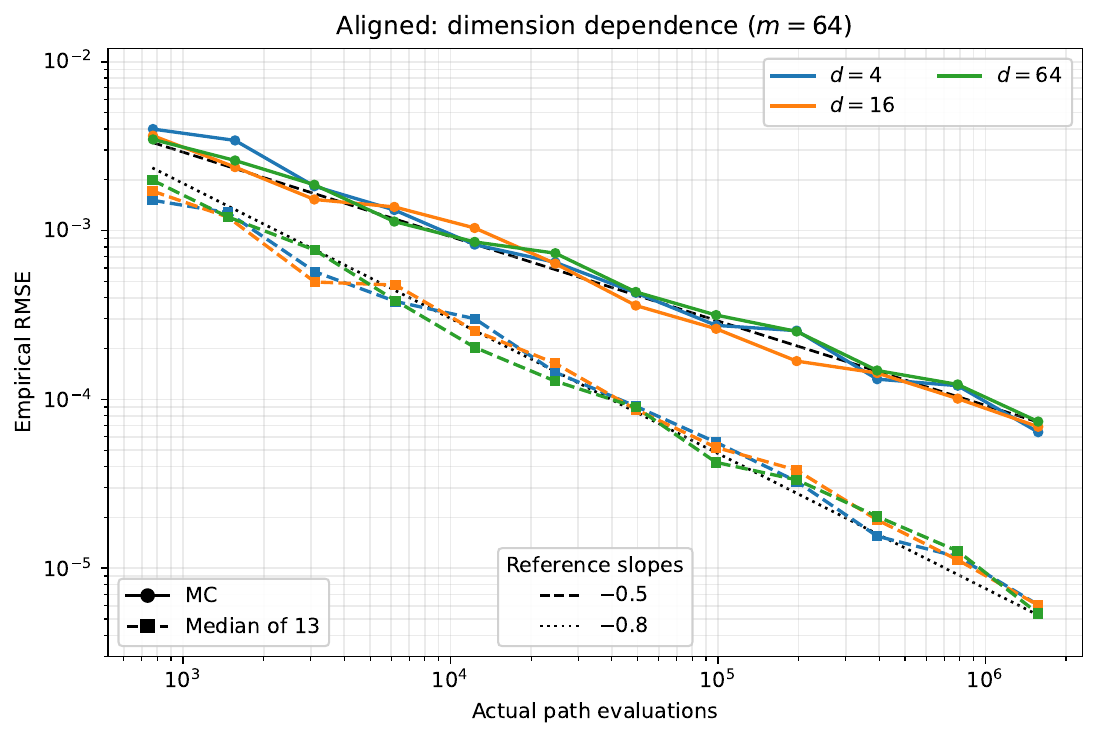}
  \par\smallskip\textbf{(a)} Aligned, fixed $m=64$.
 \end{minipage}\hfill
 \begin{minipage}[t]{0.49\linewidth}
  \centering
  \includegraphics[width=\linewidth]{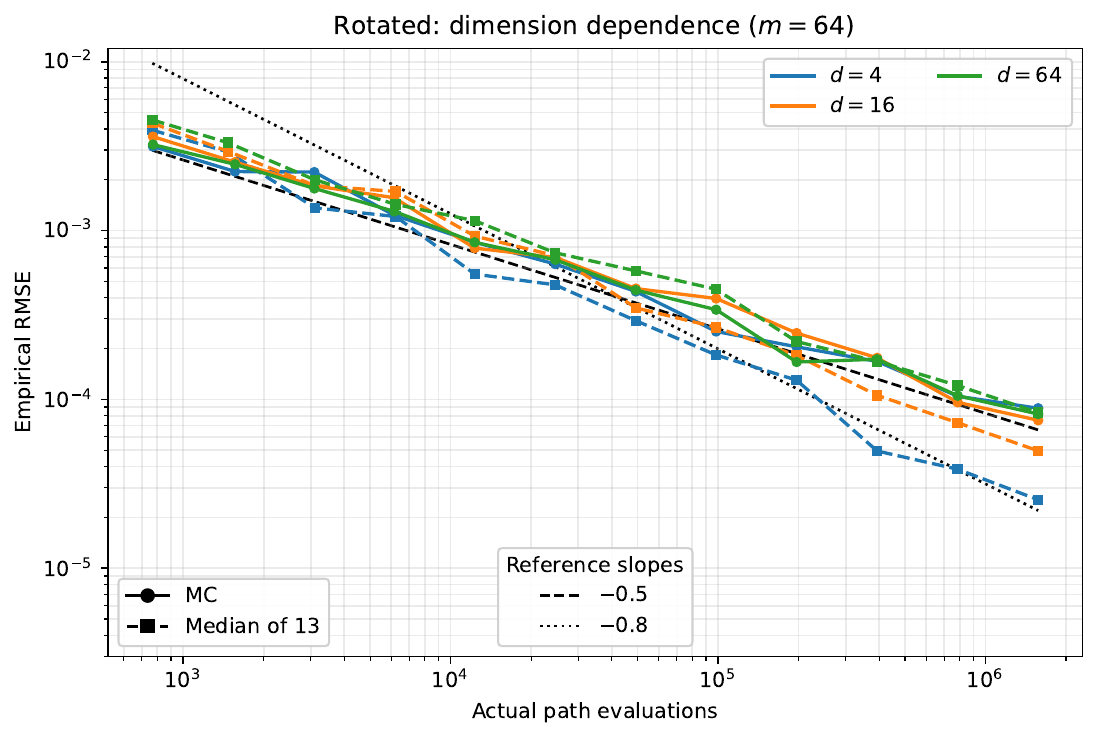}
  \par\smallskip\textbf{(b)} Rotated, fixed $m=64$.
 \end{minipage}
 \par\medskip
 \begin{minipage}[t]{0.49\linewidth}
  \centering
  \includegraphics[width=\linewidth]{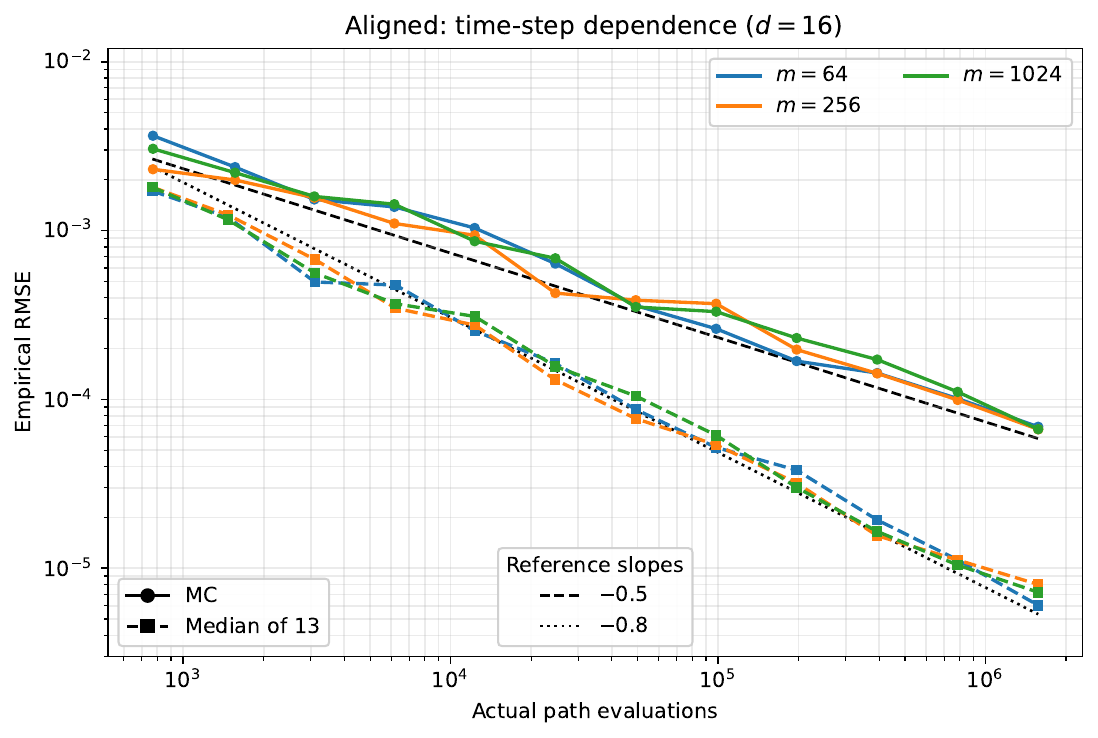}
  \par\smallskip\textbf{(c)} Aligned, fixed $d=16$.
 \end{minipage}\hfill
 \begin{minipage}[t]{0.49\linewidth}
  \centering
  \includegraphics[width=\linewidth]{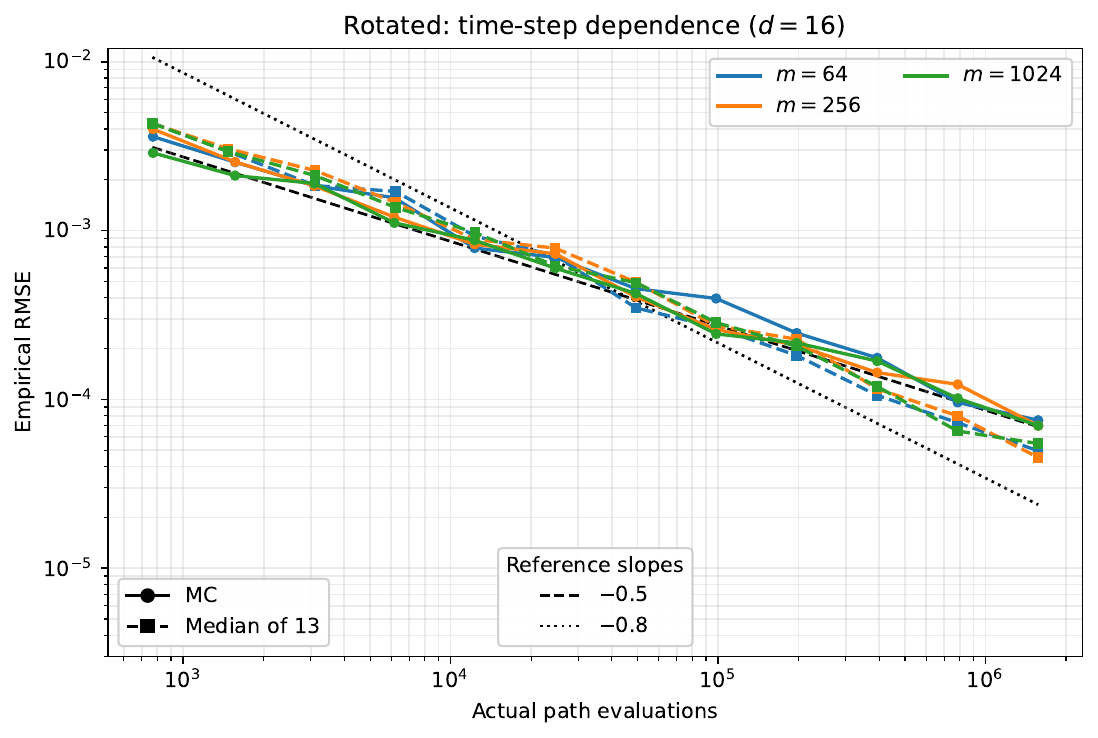}
  \par\smallskip\textbf{(d)} Rotated, fixed $d=16$.
 \end{minipage}
 \caption{Dimension and time-step comparisons for the aligned and rotated
 representations.  The underlying Euler laws are identical at common
 $(d,m)$, whereas the QMC coordinate representations differ.}
 \label{fig:p3-p4-results}
\end{figure}

Figure~\ref{fig:p3-p4-results} compares the aligned and rotated
representations using the same slices.  Panels (a)--(b) fix $m=64$ and use
$d=4,16,64$; panels (c)--(d) fix $d=16$ and use
$m=64,256,1024$.  The MC error levels are comparable between the aligned and
rotated representations.  For
$\Sigma_d^{\rm al}$, median-of-thirteen has smaller RMSE than MC and
approaches the ${n}^{-0.8}$ reference rate at large budgets.  Its error does
not visibly deteriorate as either $m$ or $d$ increases, consistent with Corollary~\ref{cor:dimensions}.  For
$\Sigma_d^{\rm rot}$, the advantage of QMC over MC is pronounced at $d=4$
but largely disappears at $d=16$ and $d=64$, where the QMC and MC errors are
comparable.
At fixed $d=16$, increasing $m$ has little effect on the QMC error. {These results illustrate that rotating $\Sigma$ as in Remark~\ref{rem:rotation} so that Assumption~\ref{ass:columns} holds preserves the target expectation and can substantially improve the efficiency of the lattice estimators.}

\section{Conclusion}\label{sec:discussion}
We studied randomly shifted rank-1 lattice rules for computing expectations of solutions
to SDEs with a constant diffusion matrix. Using the Euler--Maruyama
scheme and Brownian bridge construction, we established almost
first-order sampling error bounds. A rooted-tree expansion of the
mixed derivatives yields summability estimates uniform in the number
of time steps. These estimates give an $O(N^{-1+\varepsilon})$ mean absolute error
for a single random lattice rule and an $L^q$ error of the same order
for the median estimator, with admissible $q>1$.
We also identified conditions for error bounds uniform in the state
and {Brownian path dimensions}. The numerical experiments demonstrate the advantage of median lattice
estimators over Monte Carlo and their stability under time grid refinement.
Under the diffusion-column summability condition, the observed convergence
also remains stable as the state and {Brownian path dimensions} increase.

The current bounds retain a dependence on the time horizon $T$,
which limits their use in assessing long-time computations.
Addressing this dependence by identifying conditions for error bounds
uniform in $T$ is one direction for future work. Another is to develop
methods {whose sampling errors decay faster than $O(N^{-1})$}.
Together, these directions would extend QMC approximation of SDE
expectations to longer time horizons and higher convergence orders.

 \appendix
\section{Proof of the recursive state-derivative expansions}\label{sec:variation-proof}
\begin{proof}[Proof of Lemma~\ref{lem:state-expansion}]
Fix the input $\bm{z}$. We first consider $u=\{j\}$. Applying
\eqref{eq:chain} to \eqref{eq:Y} and then using
${\bm X^h_{i,\{j\}}}=\bm{Y}_{i,\{j\}}+\bm{\xi}_j(t_i)$ gives
\begin{align}
 \bm{Y}_{i+1,\{j\}}
 &=\bm{Y}_{i,\{j\}}
   +hD_xb(t_i,\bm{X}^h_{t_i})[{\bm X^h_{i,\{j\}}}]\notag\\
 &=\bm{Y}_{i,\{j\}}
   +hD_xb(t_i,\bm{X}^h_{t_i})
      [\bm{Y}_{i,\{j\}}+\bm{\xi}_j(t_i)]\notag\\
 &={\Omega}_i\bm{Y}_{i,\{j\}}
   +hD_xb(t_i,\bm{X}^h_{t_i})[\bm{\xi}_j(t_i)],
 \qquad \bm{Y}_{0,\{j\}}=0.
 \label{eq:singleton-recursion}
\end{align}
Since $P_{i+1,{{i'}}}={\Omega}_iP_{i,{{i'}}}$ for ${{i'}}<i$ and
$P_{i+1,i}=\operatorname{Id}_{d_x}$, induction on $i$ gives
\begin{equation}\label{eq:state-duhamel}
 \bm{Y}_{i,\{j\}}
 =\sum_{{{i'}}=0}^{i-1}P_{i,{{i'}}}hD_xb(t_{{i'}},\bm{X}^h_{t_{{i'}}})
 [\bm{\xi}_j(t_{{i'}})].
\end{equation}
This proves \eqref{eq:partition-first}.

Next let $|u|\ge2$. Applying \eqref{eq:chain} to \eqref{eq:Y} gives
\begin{align}
 \bm{Y}_{i+1,u}
 &=\bm{Y}_{i,u}+h\sum_{\rho\in\PP(u)}
 D_x^{|\rho|}b(t_i,\bm{X}^h_{t_i})
 [{\bm X^h_{i,v}}:v\in\rho]\notag\\
 &=\bm{Y}_{i,u}+hD_xb(t_i,\bm{X}^h_{t_i})[{\bm X^h_{i,u}}]
  +h\sum_{\substack{\rho\in\PP(u)\\|\rho|\ge2}}
 D_x^{|\rho|}b(t_i,\bm{X}^h_{t_i})
 [{\bm X^h_{i,v}}:v\in\rho]\notag\\
 &={\Omega}_i\bm{Y}_{i,u}
  +h\sum_{\substack{\rho\in\PP(u)\\|\rho|\ge2}}
 D_x^{|\rho|}b(t_i,\bm{X}^h_{t_i})
 [{\bm X^h_{i,v}}:v\in\rho],
 \qquad \bm{Y}_{0,u}=0,
 \label{eq:higher-recursion}
\end{align}
where the last equality uses ${\bm X^h_{i,u}}=\bm{Y}_{i,u}$. The same induction, again
using $P_{i+1,{{i'}}}={\Omega}_iP_{i,{{i'}}}$ and
$P_{i+1,i}=\operatorname{Id}_{d_x}$, yields
\begin{equation}\label{eq:state-duhamel-high}
 \bm{Y}_{i,u}
 =\sum_{{{i'}}=0}^{i-1}P_{i,{{i'}}}h
 \sum_{\substack{\rho\in\PP(u)\\|\rho|\ge2}}
 D_x^{|\rho|}b(t_{{i'}},\bm{X}^h_{t_{{i'}}})
 [{\bm X^h_{{{i'}},v}}:v\in\rho].
\end{equation}

To rewrite \eqref{eq:state-duhamel-high}, fix $\rho\in\PP(u)$ with
$|\rho|\ge2$. For $S\subseteq\operatorname{sing}(\rho)$, with
$\operatorname{sing}(\rho)$ defined in \eqref{eq:singleton-blocks}, set
\[
 \pi_S:=\rho\setminus\{\{j\}:j\in S\}.
\]
By multilinearity,
\begin{align*}
 D_x^{|\rho|}b(t_{{i'}},\bm{X}^h_{t_{{i'}}})
 [{\bm X^h_{{{i'}},v}}:v\in\rho]
 &=D_x^{|\rho|}b(t_{{i'}},\bm{X}^h_{t_{{i'}}})
 \left[(\bm{Y}_{{{i'}},\{j\}}+\bm{\xi}_j(t_{{i'}}))_{j\in\operatorname{sing}(\rho)},
       (\bm{Y}_{{{i'}},v})_{v\in\rho,\ |v|\ge2}\right]\\
 &=
 \sum_{S\subseteq\operatorname{sing}(\rho)}
 D_x^{|S|+|\pi_S|}b(t_{{i'}},\bm{X}^h_{t_{{i'}}})
 \left[(\bm{\xi}_j(t_{{i'}}))_{j\in S},
       (\bm{Y}_{{{i'}},v})_{v\in\pi_S}\right].
\end{align*}
We note that
$\pi_S\in\PP(u\setminus S)$ and
$|S|+|\pi_S|=|\rho|\ge2$. Conversely, every
$S\subseteq u$ and $\pi\in\PP(u\setminus S)$ satisfying
$|S|+|\pi|\ge2$ uniquely determine
$\rho=\{\{j\}:j\in S\}\cup\pi$ and the corresponding choices in the
multilinear expansion.
Reindexing the finite sum in \eqref{eq:state-duhamel-high} therefore
gives \eqref{eq:correspondence-reindex}.

Finally, \eqref{eq:drift-ass} and submultiplicativity give
\[
 \|{\Omega_i}\|_2\le1+h\mathsf B_1,
 \qquad
 \|P_{i,{{i'}}}\|_2
 \le(1+h\mathsf B_1)^{i-{{i'}}-1}
 \le e^{\mathsf B_1T},
\]
including the empty product. This proves \eqref{eq:propagator-bound}.
\end{proof}

\section{A tree expansion for the derivative majorant}\label{sec:expansion-counting}
This appendix establishes the estimate for $\mathcal D_{m,u}$ stated in
Proposition~\ref{prop:majorant}. We first derive an expansion of
$\partial_uF_m$ that separates the direct Brownian bridge inputs from
the state derivatives. We then use labelled rooted trees to track the
terms generated by repeatedly substituting the recursive expansions in
Lemma~\ref{lem:state-expansion}. After estimating the value of each tree
and counting the compatible tree structures, we sum over the admissible
time assignments and complete the proof of
Proposition~\ref{prop:majorant}.

\subsection{Payoff-derivative expansion}\label{sec:partition-expansion}
The following lemma separates the endpoint inputs that enter the payoff
derivative directly from the state derivatives that will be expanded
recursively.

\begin{lemma}\label{lem:payoff-expansion}
For nonempty $u\subseteq[s]$, the payoff derivative is
\begin{equation}\label{eq:partition-payoff}
 \partial_uF_m=\sum_{E\subseteq u\cap \mathcal E_w}
 \sum_{\pi\in\PP(u\setminus E)}
 {D_x^{|E|+|\pi|}}g(\bm{X}^h_{t_m})
\left[(\bm{\xi}_j(T))_{j\in E},
 (\bm{Y}_{m,v})_{v\in\pi}\right]{,}
\end{equation}
where $\mathcal E_w$ is defined in \eqref{eq:bb-index-sets}.
\end{lemma}
\begin{proof}
The partition chain rule \cite[Theorem 1]{levy2006why} gives
\[
 \partial_uF_m=\sum_{\rho\in\PP(u)}
 {D_x^{|\rho|}}g(\bm{X}^h_{t_m})[{\bm X^h_{m,v}}:v\in\rho].
\]
As in the proof of Lemma~\ref{lem:state-expansion}, expand the singleton
factors by multilinearity. For $\rho\in\PP(u)$ and
$E\subseteq\operatorname{sing}(\rho)$, set
\[
 \pi_E:=\rho\setminus\{\{j\}:j\in E\}.
\]
Then
\[
 \partial_uF_m=
 \sum_{\rho\in\PP(u)}
 \sum_{E\subseteq\operatorname{sing}(\rho)}
 {D_x^{|E|+|\pi_E|}}g(\bm{X}^h_{t_m})
 \left[(\bm{\xi}_j(T))_{j\in E},
       (\bm{Y}_{m,v})_{v\in\pi_E}\right].
\]
Any term with $E\not\subseteq u\cap\mathcal E_w$ vanishes, since it
contains some $j\in E\cap\mathcal J_m$ and $\bm{\xi}_j(T)=0$.
Reindexing the remaining
terms by
\[
 \rho=\{\{j\}:j\in E\}\cup\pi
\]
gives \eqref{eq:partition-payoff}.
\end{proof}

\subsection{Rooted-tree representation}\label{sec:tree-representation}
We encode the repeated substitutions of Lemma~\ref{lem:state-expansion} into \eqref{eq:partition-payoff} by finite rooted trees, using the standard terminology of Levin~\cite[Section~4.2]{levin2021discrete}. The root represents one summand in \eqref{eq:partition-payoff}, while each nonroot vertex $y$ is associated with an occurrence of a state derivative $\bm Y_{i,v}$. The edges record how the state derivatives in \eqref{eq:partition-payoff} are recursively expanded using Lemma~\ref{lem:state-expansion}. We now make this representation precise.

\begin{definition}\label{def:tree}
For $\varnothing\ne u\subseteq[s]$, a \textbf{labelled payoff tree}
$\mathfrak T$ has a payoff root $y_*$ and drift vertices $y\ne y_*$.  The
root carries a direct label set $S(y_*)=E\subseteq u\cap\mathcal E_w$,
while each drift vertex carries a direct label set $S(y)$, possibly empty.
These sets are pairwise disjoint and satisfy
\[
 u= \mathop{\dot\bigcup}_{y\in\mathfrak T }S(y).
\]
The root-only tree is included; in this case $E=u$.  For a nonroot vertex
$y$, let $\operatorname{pa}(y)$ denote its parent.  For every vertex $y$,
let $\operatorname{ch}(y)$ denote its set of children, and define
\begin{equation}\label{eq:child-count}
 c_y:=|\operatorname{ch}(y)|.
\end{equation}
At every drift vertex require
\begin{equation}\label{eq:drift-vertex-condition}
 (|S(y)|,c_y)=(1,0)\quad\text{or}\quad |S(y)|+c_y\ge2.
\end{equation}
Thus $S(y)\ne\varnothing$ for every drift leaf $y$.  The children of a vertex
are unordered, so permuting them does not produce a new labelled tree.  Let
$\mathcal T(u)$ denote the set of these trees, and let
$\mathcal V(\mathfrak T)$ denote the vertex set of $\mathfrak T$.
\end{definition}

\begin{definition}\label{def:time-assignment}
Fix $\mathfrak T\in\mathcal T(u)$.
A \textbf{time assignment} is a family
${\bm i}=({i}_y)_{y\in\mathcal V(\mathfrak T)}$ in
\[
 \mathcal R_m(\mathfrak T)=
 \left\{({i}_y)\in\{0,\ldots,m\}^{|\mathcal V(\mathfrak T)|}:
 {i}_{y_*}=m,\quad {i}_y<{i}_{\operatorname{pa}(y)}
 \text{ for every } y_*\ne y\in\mathcal V(\mathfrak T) \right\}.
\]
\end{definition}

\begin{definition}\label{def:tree-value}
Fix $\mathfrak T\in\mathcal T(u)$ and
${\bm i}\in\mathcal R_m(\mathfrak T)$.  Evaluate its drift vertices from the
leaves toward the root by setting
\begin{equation}\label{eq:tree-value}
 \bm V_y=P_{{i}_{\operatorname{pa}(y)},{i}_y}h
 D_x^{|S(y)|+c_y}b(t_{{i}_y},\bm X^h_{t_{{i}_y}})
 [ (\bm\xi_j(t_{{i}_y}))_{j\in S(y)},
   (\bm V_{y'})_{y'\in\operatorname{ch}(y)}].
\end{equation}
For a leaf $y$, $\operatorname{ch}(y)=\varnothing$, so $c_y=0$ and
the family $(\bm V_{y'})_{y'\in\operatorname{ch}(y)}$ is empty.  The
\textbf{value of the payoff tree} is
\begin{equation}\label{eq:root-value}
 \mathcal G(\mathfrak T,{\bm i})
 ={D_x^{|E|+c_{y_*}}}g(\bm X^h_{t_m})
 [ (\bm\xi_j(T))_{j\in E},
   (\bm V_{y'})_{y'\in\operatorname{ch}(y_*)}].
\end{equation}
Since the derivatives of $b$ and $g$ are symmetric,
$\mathcal G(\mathfrak T,{\bm i})$ is independent of the ordering of the
children.  The following lemma identifies $\partial_uF_m$ with the sum of
these tree values.
\end{definition}

\begin{lemma}\label{lem:correspondence}
For $\varnothing\ne u\subseteq[s]$, we have
\begin{equation}\label{eq:tree-payoff-sum}
 \partial_uF_m
 =\sum_{\mathfrak T\in\mathcal T(u)}
   \sum_{{\bm i}\in\mathcal R_m(\mathfrak T)}
       \mathcal G(\mathfrak T,{\bm i}).
\end{equation}
\end{lemma}
\begin{proof}
\textbf{Construction of the labelled tree and its time assignment.}
Repeatedly substitute Lemma~\ref{lem:state-expansion} whenever a state
derivative $\bm Y_{i,v}$ with $i\geq 1$ appears in
\eqref{eq:partition-payoff}, and
distribute the resulting finite sums.  Each newly generated state
derivative has a strictly smaller time index.  Hence, after at most $m$
generations, every branch  
contains no state derivative $\bm Y_{i,v}$ with $i\geq 1$.  We remove
every term containing a factor $\bm Y_{0,v}$, since $\bm Y_{0,v}=0$.
The remaining choices are complete in the sense
that no state derivative remains, so all terminal vector arguments are
direct inputs $\bm\xi_j(t)$.

Let $\mathfrak C_{m,u}$ denote the finite set of these complete expansion
choices, and let $\operatorname{val}(\chi)$ be the fully expanded scalar
term corresponding to $\chi\in\mathfrak C_{m,u}$.  Then
\begin{equation}\label{eq:formal-expansion-sum}
 \partial_uF_m
 =\sum_{\chi\in\mathfrak C_{m,u}}\operatorname{val}(\chi).
\end{equation}
Fix $\chi\in\mathfrak C_{m,u}$.  We construct from $\chi$ a labelled
payoff tree $\mathfrak T$ and a time assignment ${\bm i}$.

First consider the summand of \eqref{eq:partition-payoff} selected by
$\chi$.  It is indexed by
\[
 E\subseteq u\cap\mathcal E_w,
 \qquad \pi\in\PP(u\setminus E).
\]
Create the payoff root $y_*$, set
$S(y_*)=E$ and ${i}_{y_*}=m$.  For each $v\in\pi$, create a distinct child
$y$ of $y_*$ and set $\mathsf Y(y):=\bm Y_{m,v}$. So $|\operatorname{ch}(y_*)|=|\pi|$.  We extend the notation $\mathsf Y(y)$
recursively below.  The summand selected by $\chi$ at the payoff root can
then be written as
\begin{equation}\label{eq:selected-root-summand}
 {D_x^{|S(y_*)|+c_{y_*}}}g(\bm X^h_{t_m})
 \left[(\bm\xi_j(T))_{j\in S(y_*)},
       (\mathsf Y(y))_{y\in\operatorname{ch}(y_*)}\right].
\end{equation}
If $\pi=\varnothing$, no drift vertex is introduced and the tree consists
only of the payoff root.

Suppose that $\mathsf Y(y)=\bm Y_{i,v}$ with $i\ge 1$ for a drift vertex $y$. Then ${i}_{\operatorname{pa}(y)}=i$.
We distinguish two cases.  If $v=\{j\}$, the expansion choice $\chi$ selects a summand of
\eqref{eq:partition-first}, indexed by $0\le {{i'}}<i$, namely
\[
 P_{i,{{i'}}}hD_xb(t_{{i'}},\bm X^h_{t_{{i'}}})[\bm\xi_j(t_{{i'}})].
\]
We set $ {i}_y={{i'}},
S(y)=\{j\}.$
This summand requires no further expansion, so no child is attached to
$y$.  Hence $c_y=0$, $(|S(y)|,c_y)=(1,0)$, and $y$ is a leaf.
If $|v|\ge2$, the expansion choice $\chi$ selects a summand of
\eqref{eq:correspondence-reindex}, indexed by
\[
 0\le {{i'}}<i,
 \qquad
 S\subseteq v,
 \qquad
 \pi\in\PP(v\setminus S),
 \qquad
 |S|+|\pi|\ge2.
\]
This summand is
\begin{equation}\label{eq:tree-proof-local-term}
 P_{i,{{i'}}}hD_x^{|S|+|\pi|}b(t_{{i'}},\bm X^h_{t_{{i'}}})
 \left[(\bm\xi_j(t_{{i'}}))_{j\in S},
       (\bm Y_{{{i'}},v'})_{v'\in\pi}\right].
\end{equation}
Set ${i}_y={{i'}}$ and $S(y)=S$.  For each $v'\in\pi$, create a distinct child
$w$ of $y$ and set $\mathsf Y(w):=\bm Y_{{{i'}},v'}$.  Thus $c_y=|\pi|$ and
$|S(y)|+c_y\ge2$, while $y$ is a leaf
exactly when $\pi=\varnothing$.  Moreover, for a retained term,
${i}_y=0$ forces $\pi=\varnothing$, since every term containing a factor
$\bm Y_{0,v'}$ has been removed.  In both cases,
${i}_y<i={i}_{\operatorname{pa}(y)}$, and the local summand selected by
$\chi$ at $y$ can be written uniformly as
\begin{equation}\label{eq:selected-drift-summand}
 P_{{i}_{\operatorname{pa}(y)},{i}_y}h
 D_x^{|S(y)|+c_y}b(t_{{i}_y},\bm X^h_{t_{{i}_y}})
 \left[(\bm\xi_j(t_{{i}_y}))_{j\in S(y)},
       (\mathsf Y(w))_{w\in\operatorname{ch}(y)}\right].
\end{equation}
Applying the same rule recursively to every new child completes the
construction of $\mathfrak T$ and ${\bm i}$ from $\chi$.

\textbf{Admissibility and preservation of values.}
We first verify that $ \mathfrak T\in\mathcal T(u),
{\bm i}\in\mathcal R_m(\mathfrak T).$
The two alternatives in the construction above give \eqref{eq:drift-vertex-condition}.  If
$\mathsf Y(y)=\bm Y_{i,v}$, let $\mathfrak T_y$ denote the
subtree of $\mathfrak T$ rooted at $y$.  Induction over $\mathfrak T_y$
gives
\begin{equation}\label{eq:subtree-label-identity}
 v=\mathop{\dot\bigcup}_{w\in\mathcal V(\mathfrak T_y)}S(w).
\end{equation}
Applying this identity to the children of the payoff root and recalling
their association with the state-derivative factors in the selected
summand of \eqref{eq:partition-payoff}, we obtain
\[
 u=\mathop{\dot\bigcup}_{w\in\mathcal V(\mathfrak T)}S(w).
\]
Together with
${i}_{y_*}=m$ and the strict time inequalities, this proves that the
construction yields
$\mathfrak T\in\mathcal T(u)$ and
${\bm i}\in\mathcal R_m(\mathfrak T)$.

It remains to prove
$\operatorname{val}(\chi)=\mathcal G(\mathfrak T,{\bm i})$.
For every drift vertex $y$, we claim that recursively applying
Lemma~\ref{lem:state-expansion} to the state-derivative factors in
\eqref{eq:selected-drift-summand}, with summands selected according
to $\chi$, yields $\bm V_y$.
We prove this by induction from the leaves toward the root.

For a leaf $y$, $\operatorname{ch}(y)=\varnothing$, so the claim
follows directly by comparing \eqref{eq:selected-drift-summand}
with \eqref{eq:tree-value}.
For a nonleaf vertex $y$, each state-derivative factor
$\mathsf Y(w)$ in \eqref{eq:selected-drift-summand} is associated
with a child $w$. By the construction and the induction hypothesis,
the term selected by $\chi$ in the expansion of $\mathsf Y(w)$
equals $\bm V_w$.
By multilinearity, the selected term at $y$ is therefore obtained
by replacing each $\mathsf Y(w)$ in
\eqref{eq:selected-drift-summand} with $\bm V_w$, which gives
\eqref{eq:tree-value}.

Applying the same argument to \eqref{eq:selected-root-summand}
and comparing with \eqref{eq:root-value} gives
\begin{equation}\label{eq:expansion-tree-value}
 \operatorname{val}(\chi)=\mathcal G(\mathfrak T,{\bm i}).
\end{equation}
This identity holds directly in the root-only case.

\textbf{Inverse construction.}
Conversely, fix $\mathfrak T\in\mathcal T(u)$ and
${\bm i}\in\mathcal R_m(\mathfrak T)$. For every drift vertex $y$, recall that
$\mathfrak T_y$ {denotes} the subtree of $\mathfrak T$ rooted at $y$,
and define its coordinate set by
\begin{equation}\label{eq:subtree-coordinate-set}
 v(y):=\mathop{\dot\bigcup}_{w\in\mathcal V(\mathfrak T_y)}S(w).
\end{equation}
The child coordinate sets are nonempty and pairwise disjoint by
\eqref{eq:drift-vertex-condition} and the disjointness of the direct
label sets. For every vertex $y$, set
\[
 \pi_y:=\{v(w):w\in\operatorname{ch}(y)\}.
\]
At the payoff root, $E=S(y_*)$ and $\pi=\pi_{y_*}$ determine a unique
summand of \eqref{eq:partition-payoff}. At a drift vertex $y$, the
associated derivative is $\bm Y_{{i}_{\operatorname{pa}(y)},v(y)}$ and
$\pi_y\in\PP(v(y)\setminus S(y))$. If $(|S(y)|,c_y)=(1,0)$,
then $v(y)=S(y)$ is a singleton, which together with ${i}_y$
determines a unique summand of \eqref{eq:partition-first}. Otherwise
$|v(y)|\ge |S(y)|+c_y\ge2$, so $S(y)$, $\pi_y$, and
${i}_y<{i}_{\operatorname{pa}(y)}$ select a unique summand of
\eqref{eq:correspondence-reindex}.

The strict time inequalities force every vertex at time zero to be a
leaf, so no recovered summand contains a factor $\bm Y_{0,v}$.
These local choices therefore form a complete choice
$\chi\in\mathfrak C_{m,u}$.

\textbf{Verification of the bijection.}
First, let $(\mathfrak T,{\bm i})$ be constructed from $\chi$.
Comparing \eqref{eq:subtree-label-identity} with
\eqref{eq:subtree-coordinate-set} gives $v(y)=v$ whenever
$\mathsf Y(y)=\bm Y_{i,v}$. Since $i={i}_{\operatorname{pa}(y)}$, we have $\mathsf Y(y)=\bm Y_{{i}_{\operatorname{pa}(y)},v(y)}.$
Thus the inverse construction recovers the same associated
derivatives. Applying the root and drift-vertex construction
rules recursively recovers the same selected summands, hence
the original $\chi$.

Conversely, fix $(\mathfrak T,{\bm i})$ and let $\chi$ be obtained
by the inverse construction. The same vertex-by-vertex comparison
shows that the forward construction applied to this $\chi$
returns the original pair $(\mathfrak T,{\bm i})$.
Thus the two constructions are mutually inverse and give a bijection.

Combining this bijection with \eqref{eq:expansion-tree-value} and
\eqref{eq:formal-expansion-sum} proves \eqref{eq:tree-payoff-sum}.
\end{proof}

\subsection{Tree counts and coefficient bounds}\label{sec:expression-counting}
We establish three technical lemmas controlling the size,
weighted count, and coefficients of trees with a fixed direct partition.

For a labelled payoff tree $\mathfrak T\in\mathcal T(u)$, define
its global direct partition by
\[
 \Pi(\mathfrak T)
 :=\bigl\{S(y):
 y\in\mathcal V(\mathfrak T)\setminus\{y_*\},\ S(y)\ne\varnothing\bigr\}.
\]
By Definition~\ref{def:tree},
$\Pi(\mathfrak T)\in\PP(u\setminus S(y_*))$.
Throughout this subsection, fix $E\subseteq u\cap\mathcal E_w$
such that ${{\bar u}}:=u\setminus E\ne\varnothing$, and fix $\pi\in\PP({{\bar u}})$.
Define
\[
 \Gamma(E,\pi):=\bigl\{\mathfrak T\in\mathcal T(u):
 S(y_*)=E,\ \Pi(\mathfrak T)=\pi\bigr\}.
\]
Write
\[
 {{n_\pi}}:=|\pi|=|\Pi(\mathfrak T)|
 \quad\text{for }\mathfrak T\in\Gamma(E,\pi).
\]

\begin{lemma}\label{lem:expression-structure}
Let $\mathfrak T\in\Gamma(E,\pi)$ and define
\begin{equation}\label{eq:empty-drift-vertex-count}
 {{n_\varnothing}}:=\bigl|\{y\in\mathcal V(\mathfrak T)\setminus\{y_*\}:
 S(y)=\varnothing\}\bigr|.
\end{equation}
Then
\[
 {{n_\varnothing}}\le {{n_\pi}}-1,
 \qquad
 |\mathcal V(\mathfrak T)|={{n_\pi}}+{{n_\varnothing}}+1\le2{{n_\pi}}.
\]
Moreover,
\begin{equation}\label{eq:global-input-factorial}
 {\mathfrak c}(\mathfrak T):=\prod_{y\in\mathcal V(\mathfrak T)}c_y!\le {{n_\pi}}!.
\end{equation}
\end{lemma}
\begin{proof}
There are ${{n_\pi}}+{{n_\varnothing}}+1$ vertices and ${{n_\pi}}+{{n_\varnothing}}$ edges. Each edge joins a child
to its parent, so
$\sum_{y\in\mathcal V(\mathfrak T)}c_y={{n_\pi}}+{{n_\varnothing}}$. By
\eqref{eq:drift-vertex-condition}, if $y\ne y_*$ and
$S(y)=\varnothing$, then $c_y\ge2$. Since ${{\bar u}}\ne\varnothing$,
$c_{y_*}\ge1$. Consequently
\[
 {{n_\pi}}+{{n_\varnothing}}=\sum_{y\in\mathcal V(\mathfrak T)}c_y
 \ge2{{n_\varnothing}}+1,\qquad {{n_\varnothing}}\le {{n_\pi}}-1.
\]
We then prove an auxiliary bound valid for every finite rooted tree
$\mathfrak S$, without imposing the conditions defining $\mathfrak T$.
Let $L(\mathfrak S)$ denote its number of leaves. We claim that
$\prod_{y\in\mathcal V(\mathfrak S)}c_y!\le L(\mathfrak S)!$.
We argue by induction on the number of vertices. For the one-vertex
tree, the root has no children and is itself a leaf, so both sides
equal $1$. For the induction step, let the root have $c\ge1$ child
subtrees $\mathfrak S_1,\ldots,\mathfrak S_c$,
with respective leaf counts $L_1,\ldots,L_c\ge1$. Applying the
induction hypothesis to these subtrees gives
\[
 \prod_{y\in\mathcal V(\mathfrak S)}c_y!
 =c!\prod_{a=1}^c\prod_{y'\in\mathcal V(\mathfrak S_a)}c_{y'}!
 \le c!\prod_{a=1}^cL_a!
 \le\left(\sum_{a=1}^cL_a\right)!=L(\mathfrak S)!.
\]
The second inequality holds because, for $L_a\ge1$, the multinomial
coefficient $(\sum_aL_a)!/\prod_aL_a!$ is at least its value $c!$
at $L_1=\cdots=L_c=1$: increasing $L_j$ multiplies it by
$(\sum_aL_a+1)/(L_j+1)\ge1$.
Finally, apply this auxiliary bound to $\mathfrak T$. Since
${{\bar u}}\ne\varnothing$, the {payoff} root is not a leaf. By
\eqref{eq:drift-vertex-condition}, every drift leaf has a nonempty
label set. Thus all leaves are drift vertices with nonempty label
sets, and their number $L_*$ is at most ${{n_\pi}}$. Hence
\[
 \prod_{y\in\mathcal V(\mathfrak T)}c_y!\le L_*!\le {{n_\pi}}!,
\]
which proves \eqref{eq:global-input-factorial}.
\end{proof}

\begin{lemma}\label{lem:weighted-pattern-count}
With ${\mathfrak c}(\mathfrak T)$ defined in \eqref{eq:global-input-factorial}, we have
\begin{equation}\label{eq:weighted-pattern-count}
 \sum_{\mathfrak T\in\Gamma(E,\pi)}{\mathfrak c}(\mathfrak T)\le64^{{n_\pi}}{{n_\pi}}!.
\end{equation}
\end{lemma}
\begin{proof}
Temporarily order the children at each vertex. Since sibling subtrees carry pairwise disjoint nonempty coordinate sets, each \(\mathfrak T\) admits exactly \({\mathfrak c}(\mathfrak T)\) such orderings.

By Lemma~\ref{lem:expression-structure}, there are at most $2{{n_\pi}}$
vertices. By the classical Catalan enumeration of plane rooted trees
\cite[Section~I.5]{flajolet2009analytic}, the number of ordered rooted
shapes with $n$ vertices is
\[
 \mathrm{Cat}_{n-1}
 =\frac{1}{n}\binom{2n-2}{n-1}\le4^{n-1}.
\]
For each shape, the root carries $E$; choose the ${{n_\pi}}$ other vertices
carrying nonempty direct sets in $\binom{n-1}{{{n_\pi}}}$ ways and assign the
distinct blocks of $\pi$ in ${{n_\pi}}!$ ways. The subtree coordinate sets
determine all local partitions, so no further choices are needed.
Including inadmissible choices only enlarges the bound. Thus
\begin{equation*}
 \sum_{\mathfrak T\in\Gamma(E,\pi)}{\mathfrak c}(\mathfrak T)
\le {{n_\pi}}!\sum_{n={{n_\pi}}+1}^{2{{n_\pi}}}4^{n-1}\binom{n-1}{{{n_\pi}}}
\le {{n_\pi}}!\sum_{n=1}^{2{{n_\pi}}}8^{n-1}
={{n_\pi}}!\frac{8^{2{{n_\pi}}}-1}{7}
\le64^{{n_\pi}}{{n_\pi}}!. 
\end{equation*}
\end{proof}

\begin{lemma}\label{lem:count}
For $\mathfrak T\in\Gamma(E,\pi)$ and
${\bm i}\in\mathcal R_m(\mathfrak T)$, define
\[
 {i}_{\mathfrak T}(j):={i}_y
 \quad\text{for the unique $y\in\mathcal V(\mathfrak T)$ with $j\in S(y)$},
 \qquad j\in {{\bar u}}.
\]
Recall from Assumption~\ref{ass:space} that
\[
 \mathsf B_n=B_0B_1^n(n!)^\tau,\qquad
 \mathsf G_n=G_0G_1^n(n!)^\tau,\qquad n\ge1.
\]
Define
\[
 {\mathcal K}_{\mathfrak T}:=\mathsf G_{|E|+c_{y_*}}
       \prod_{j\in E}|\bm{\xi}_j(T)|_2
       \prod_{y\in\mathcal V(\mathfrak T)\setminus\{y_*\}}
       e^{\mathsf B_1T}\mathsf B_{|S(y)|+c_y}.
\]
Then we have,
\begin{equation}\label{eq:single-coefficient-bound}
 e^{-\beta_g|\bm{X}^h_{t_m}|_2}|\mathcal G(\mathfrak T,{\bm i})|
 \le {\mathcal K}_{\mathfrak T} h^{|\mathcal V(\mathfrak T)|-1}
       \prod_{j\in {{\bar u}}}|\bm{\xi}_j(t_{{i}_{\mathfrak T}(j)})|_2
\end{equation}
 Moreover, there exists $C_{\rm cnt}\ge1$, depending only on
$T,B_0,B_1,G_0,G_1,\tau$, such that
\begin{equation}\label{eq:count}
 \sum_{\mathfrak T\in\Gamma(E,\pi)}{\mathcal K}_{\mathfrak T}
 \le C_{\rm cnt}^{|u|+{{n_\pi}}}(|E|!)^\tau
       \prod_{j\in E}(\sqrt T\,\sigma_j)\,({{n_\pi}}!)^\tau
       \prod_{S\in\pi}(|S|!)^\tau.
\end{equation}
\end{lemma}
\begin{proof}
By \eqref{eq:tree-value}, \eqref{eq:drift-ass}, and
\eqref{eq:propagator-bound}, every drift vertex satisfies
\[
 |\bm V_y|_2
 \le he^{\mathsf B_1T}\mathsf B_{|S(y)|+c_y}
      \prod_{j\in S(y)}|\bm\xi_j(t_{{i}_y})|_2
      \prod_{w\in\operatorname{ch}(y)}|\bm V_w|_2.
\]
At the payoff root, \eqref{eq:root-value} and \eqref{eq:payoff-ass} give
\[
 e^{-\beta_g|\bm X^h_{t_m}|_2}|\mathcal G(\mathfrak T,{\bm i})|
 \le \mathsf G_{|E|+c_{y_*}}
      \prod_{j\in E}|\bm\xi_j(T)|_2
      \prod_{w\in\operatorname{ch}(y_*)}|\bm V_w|_2.
\]
Recursively substituting the drift bounds into the root bound gives
\begin{align*}
 e^{-\beta_g|\bm X^h_{t_m}|_2}|\mathcal G(\mathfrak T,{\bm i})|
 &\le \mathsf G_{|E|+c_{y_*}}
       \prod_{j\in E}|\bm\xi_j(T)|_2
       \prod_{y\in\mathcal V(\mathfrak T)\setminus\{y_*\}}\left(
         he^{\mathsf B_1T}\mathsf B_{|S(y)|+c_y}
         \prod_{j\in S(y)}|\bm\xi_j(t_{{i}_y})|_2
       \right)\\
 &={\mathcal K}_{\mathfrak T}h^{|\mathcal V(\mathfrak T)|-1}
       \prod_{y\in\mathcal V(\mathfrak T)\setminus\{y_*\}}
       \prod_{j\in S(y)}|\bm\xi_j(t_{{i}_y})|_2\\
 &={\mathcal K}_{\mathfrak T}h^{|\mathcal V(\mathfrak T)|-1}
       \prod_{j\in {{\bar u}}}|\bm\xi_j(t_{{i}_{\mathfrak T}(j)})|_2.
\end{align*}
This proves \eqref{eq:single-coefficient-bound}.

For the coefficient sum, put $n_y=|S(y)|$ and let ${{n_\varnothing}}$ count drift
vertices with empty direct sets. Then
\[
 n_{y_*}=|E|,\qquad
 \sum_{y\in\mathcal V(\mathfrak T)} n_y=|u|,\qquad
 \sum_{y\in\mathcal V(\mathfrak T)} c_y
 =|\mathcal V(\mathfrak T)|-1={{n_\pi}}+{{n_\varnothing}}.
\]
Define
\[
 R(E):=(|E|!)^\tau\prod_{j\in E}|\bm\xi_j(T)|_2
      =(|E|!)^\tau\prod_{j\in E}(\sqrt T\,\sigma_j).
\]
At every vertex $y$,
\[
 ((n_y+c_y)!)^\tau
 \le 2^{\tau(n_y+c_y)}(n_y!)^\tau(c_y!)^\tau.
\]
Expanding the coefficient bounds and applying this inequality gives
\begin{align*}
 {\mathcal K}_{\mathfrak T}
 &=G_0G_1^{n_{y_*}+c_{y_*}}
       ((n_{y_*}+c_{y_*})!)^\tau
       \prod_{j\in E}|\bm\xi_j(T)|_2
       \prod_{y\in\mathcal V(\mathfrak T)\setminus\{y_*\}}
       e^{\mathsf B_1T}B_0B_1^{n_y+c_y}((n_y+c_y)!)^\tau\\
 &\le G_0G_1^{n_{y_*}+c_{y_*}}
       2^{\tau(n_{y_*}+c_{y_*})}
       (n_{y_*}!)^\tau(c_{y_*}!)^\tau
       \prod_{j\in E}|\bm\xi_j(T)|_2\\
 &\qquad{}\times\prod_{y\in\mathcal V(\mathfrak T)\setminus\{y_*\}}
       e^{\mathsf B_1T}B_0B_1^{n_y+c_y}
       2^{\tau(n_y+c_y)}(n_y!)^\tau(c_y!)^\tau\\
 &=R(E)G_0(e^{B_0B_1T}B_0)^{{{n_\pi}}+{{n_\varnothing}}}
       G_1^{n_{y_*}+c_{y_*}}B_1^{\sum_{y\in\mathcal V(\mathfrak T)\setminus\{y_*\}}(n_y+c_y)}\\
 &\qquad{}\times(2^\tau)^{\sum_{y\in\mathcal V(\mathfrak T)}(n_y+c_y)}
       \prod_{S\in\pi}(|S|!)^\tau {\mathfrak c}(\mathfrak T)^\tau.
\end{align*}
Here the nonempty drift direct sets are exactly the blocks of $\pi$;
empty sets contribute $0!=1$, and the child factorials give
$\prod_{y\in\mathcal V(\mathfrak T)}(c_y!)^\tau={\mathfrak c}(\mathfrak T)^\tau$
by \eqref{eq:global-input-factorial}.
Set
\[
 K=\max\{1,G_0,G_1,e^{B_0B_1T}B_0,B_1,2^\tau\}.
\]
Since $\sum_{y\in\mathcal V(\mathfrak T)}(n_y+c_y)=|u|+{{n_\pi}}+{{n_\varnothing}}$ and ${{n_\varnothing}}\le {{n_\pi}}-1$,
\begin{align*}
 {\mathcal K}_{\mathfrak T}
\le R(E)K^{1+2|u|+3({{n_\pi}}+{{n_\varnothing}})}
       \prod_{S\in\pi}(|S|!)^\tau {\mathfrak c}(\mathfrak T)^\tau\le R(E)K^{2|u|+6{{n_\pi}}}
       \prod_{S\in\pi}(|S|!)^\tau {\mathfrak c}(\mathfrak T)^\tau.
\end{align*}
Taking $ C_0=K^6$,
we have $K^{2|u|+6{{n_\pi}}}\le C_0^{|u|+{{n_\pi}}}$.
Moreover, Lemma~\ref{lem:expression-structure} and $\tau\ge1$ give
\[
 {\mathfrak c}(\mathfrak T)^\tau
 ={\mathfrak c}(\mathfrak T)^{\tau-1}{\mathfrak c}(\mathfrak T)
 \le({{n_\pi}}!)^{\tau-1}{\mathfrak c}(\mathfrak T).
\]
Consequently,
\begin{equation}\label{eq:splitfactorial}
 {\mathcal K}_{\mathfrak T}\le C_0^{|u|+{{n_\pi}}}R(E)
       \prod_{S\in\pi}(|S|!)^\tau({{n_\pi}}!)^{\tau-1}{\mathfrak c}(\mathfrak T).
\end{equation}
Summing and applying \eqref{eq:weighted-pattern-count}, we obtain
\begin{align*}
 \sum_{\mathfrak T\in\Gamma(E,\pi)}{\mathcal K}_{\mathfrak T}
 &\le C_0^{|u|+{{n_\pi}}}R(E)
       \prod_{S\in\pi}(|S|!)^\tau({{n_\pi}}!)^{\tau-1}
       \sum_{\mathfrak T\in\Gamma(E,\pi)}{\mathfrak c}(\mathfrak T)\\
 &\le C_0^{|u|+{{n_\pi}}}64^{{n_\pi}} R(E)({{n_\pi}}!)^\tau
       \prod_{S\in\pi}(|S|!)^\tau.
\end{align*}
Hence, taking
\[
 C_{\rm cnt}=64C_0=64K^6,
\]
proves \eqref{eq:count}. This constant depends only on
$T,B_0,B_1,G_0,G_1,\tau$.
\end{proof}

\subsection{Proof of the derivative majorant}\label{sec:majorant-proof}
\begin{proof}[Proof of Proposition~\ref{prop:majorant}]
Fix $E\subseteq u\cap\mathcal E_w$, put ${{\bar u}}=u\setminus E$, and
fix $\pi\in\PP({{\bar u}})$. If ${{\bar u}}\ne\varnothing$, put ${{n_\pi}}=|\pi|$.
For $\mathfrak T\in\Gamma(E,\pi)$, let ${{n_\varnothing}}$ be the number of drift
vertices with empty direct label sets, as defined in
\eqref{eq:empty-drift-vertex-count}. Then $\mathfrak T$ has ${{n_\pi}}+{{n_\varnothing}}$ drift vertices.
At each drift vertex $y$ put
\[
 f_y(i)=\prod_{j\in S(y)}|\bm\xi_j(t_i)|_2,
\]
with the empty product equal to one.
Enlarging the admissible times to independent choices gives
\begin{align}
 h^{{{n_\pi}}+{{n_\varnothing}}}\sum_{{\bm i}\in\mathcal R_m(\mathfrak T)}
       \prod_{y\in\mathcal V(\mathfrak T)\setminus\{y_*\}}f_y({i}_y)
 &\le\prod_{y\in\mathcal V(\mathfrak T)\setminus\{y_*\}}
       \left(h\sum_{i=0}^{m-1}f_y(i)\right)
 =T^{{{n_\varnothing}}}\prod_{S\in\pi}\eta_m(S).
 \label{eq:time-factor}
\end{align}
The equality follows because each nonempty direct set $S\in\pi$
contributes $\eta_m(S)$ by \eqref{eq:eta}, while each empty direct set
contributes $mh=T$.

Let $C_{\rm cnt}\ge1$ be the constant in \eqref{eq:count} and choose
\[
 C_{\rm maj}=\max\{C_{\rm cnt}\max(1,T),\,G_1\max(1,G_0)\}.
\]
This constant depends only on $T,B_0,B_1,G_0,G_1,\tau$.
Combining
\eqref{eq:single-coefficient-bound}, \eqref{eq:time-factor}, and
\eqref{eq:count} gives
\begin{align}
 &e^{-\beta_g|\bm X^h_{t_m}|_2}
 \sum_{\mathfrak T\in\Gamma(E,\pi)}
 \sum_{{\bm i}\in\mathcal R_m(\mathfrak T)}|\mathcal G(\mathfrak T,{\bm i})|\notag\\
 &\le\sum_{\mathfrak T\in\Gamma(E,\pi)}
       {\mathcal K}_{\mathfrak T}h^{|\mathcal V(\mathfrak T)|-1}
       \sum_{{\bm i}\in\mathcal R_m(\mathfrak T)}
       \prod_{y\in\mathcal V(\mathfrak T)\setminus\{y_*\}}f_y({i}_y)\notag\\
 &\le\prod_{S\in\pi}\eta_m(S)
       \sum_{\mathfrak T\in\Gamma(E,\pi)}{\mathcal K}_{\mathfrak T} T^{{{n_\varnothing}}}\notag\\
 &\le\max(1,T)^{{n_\pi}}\prod_{S\in\pi}\eta_m(S)
       \sum_{\mathfrak T\in\Gamma(E,\pi)}{\mathcal K}_{\mathfrak T}\notag\\
 &\le C_{\rm maj}^{|u|+{{n_\pi}}}(|E|!)^\tau
       \prod_{j\in E}(\sqrt T\,\sigma_j)\,({{n_\pi}}!)^\tau
       \prod_{S\in\pi}(|S|!)^\tau\eta_m(S).
 \label{eq:fixed-partition-contribution}
\end{align}
Here the third inequality uses ${{n_\varnothing}}\le {{n_\pi}}-1$, which gives
$T^{{{n_\varnothing}}}\le\max(1,T)^{{n_\pi}}$. The last inequality absorbs
$\max(1,T)^{{n_\pi}}$ into $C_{\rm maj}^{|u|+{{n_\pi}}}$.

If ${{\bar u}}=\varnothing$, then $E=u\subseteq\mathcal E_w$ and the tree
consists only of its payoff root. By \eqref{eq:payoff-ass},
\begin{align*}
 e^{-\beta_g|\bm X^h_{t_m}|_2}
 \left|{D_x^{|u|}}g(\bm X^h_{t_m})[\bm\xi_j(T):j\in u]\right|\le G_0G_1^{|u|}(|u|!)^\tau
       \prod_{j\in u}(\sqrt T\,\sigma_j)
 \le C_{\rm maj}^{|u|}(|u|!)^\tau\prod_{j\in u}(\sqrt T\,\sigma_j).
\end{align*}
Thus the final bound in \eqref{eq:fixed-partition-contribution} also
holds for ${{\bar u}}=\varnothing,\pi = \varnothing$.

Finally, the triangle inequality in \eqref{eq:tree-payoff-sum} and
the final bound in \eqref{eq:fixed-partition-contribution} give
\begin{align*}
 e^{-\beta_g|\bm X^h_{t_m}|_2}|\partial_uF_m(\bm z)| &\le e^{-\beta_g|\bm X^h_{t_m}|_2}
       \sum_{E\subseteq u\cap\mathcal E_w}
       \sum_{\pi\in\PP(u\setminus E)}
       \sum_{\mathfrak T\in\Gamma(E,\pi)}
       \sum_{{\bm i}\in\mathcal R_m(\mathfrak T)}
       |\mathcal G(\mathfrak T,{\bm i})|\\
 & \le\sum_{E\subseteq u\cap\mathcal E_w}
       \sum_{\pi\in\PP(u\setminus E)}
       C_{\rm maj}^{|u|+|\pi|}(|E|!)^\tau
       \prod_{j\in E}(\sqrt T\,\sigma_j)
       (|\pi|!)^\tau\prod_{S\in\pi}(|S|!)^\tau\eta_m(S)\\
 & =C_{\rm maj}^{|u|}\sum_{E\subseteq u\cap\mathcal E_w}
       (|E|!)^\tau\prod_{j\in E}(\sqrt T\,\sigma_j)
       \sum_{\pi\in\PP(u\setminus E)}
       C_{\rm maj}^{|\pi|}(|\pi|!)^\tau
       \prod_{S\in\pi}(|S|!)^\tau\eta_m(S).
\end{align*}
The right-hand side is independent of $\bm z$. Taking the supremum
over $\bm z\in\R^s$ proves \eqref{eq:majorant}.
\end{proof}

\section{Proofs of the auxiliary summability estimates}
\label{sec:regularity-proofs}
This appendix proves Lemmas~\ref{lem:detail-sum} and
\ref{lem:all-subset-sum}.

\begin{proof}[Proof of Lemma~\ref{lem:detail-sum}]
Fix $m\in\mathcal M$ and $\varnothing\ne S\subseteq\mathcal J_m$,
and put $n=|S|$. If $\bigcap_{j\in S}\operatorname{int}(I_j)=\varnothing$, then
at every grid point $t_i$ at least one $\bm\xi_j(t_i)$ vanishes,
since $\bm\xi_j$ is zero outside $\operatorname{int}(I_j)$
for every $j\in\mathcal J_m$. Consequently,
\[
 \eta_m(S)=h\sum_{i=0}^{m-1}\prod_{j\in S}|\bm{\xi}_j(t_i)|_2=0.
\]
Otherwise choose $j_*\in S$ such that
$|I_{j_*}|=\min_{j\in S}|I_j|$.
By \eqref{eq:detail-support}, for $j,j'\in\mathcal J_m$,
either $\operatorname{int}(I_j)\cap\operatorname{int}(I_{j'})
=\varnothing${, or $I_j\subseteq I_{j'}$, or $I_{j'}\subseteq I_j$}. Since the intervals indexed by $S$ have a nonempty
interior intersection, every $I_j$, $j\in S$, contains $I_{j_*}$.
Thus a summand in $\eta_m(S)$ can be nonzero only when
$t_i\in I_{j_*}$.

The interval $I_{j_*}$ contains at most $1+|I_{j_*}|/h$ grid points.
Write $|I_{j_*}|=T2^{-\ell_*}$, where
$\ell_*\in\{0,\ldots,J-1\}$. Since $h=T2^{-J}$, we have
\[
 |I_{j_*}|\ge T2^{-(J-1)}=2h,
 \qquad
 h\left(1+\frac{|I_{j_*}|}{h}\right)
 =h+|I_{j_*}|\le\frac{3}{2}|I_{j_*}|.
\]
For each $j\in S$, the height bound in \eqref{eq:bridgebounds} gives $ |\bm{\xi}_j(t_i)|_2\le\frac{\sigma_{r(j)}}{2}|I_j|^{1/2}. $
Applying these bounds to the defining sum, we obtain
\begin{align*}
 \eta_m(S)
 &=h\sum_{i=0}^{m-1}\prod_{j\in S}|\bm{\xi}_j(t_i)|_2 =h\sum_{\substack{0\le i\le m-1\\t_i\in I_{j_*}}}\prod_{j\in S}|\bm{\xi}_j(t_i)|_2 \le h\left(1+\frac{|I_{j_*}|}{h}\right)
 \prod_{j\in S}\left(\frac{\sigma_{r(j)}}{2}|I_j|^{1/2}\right)\\
 &\le\frac{3}{2}|I_{j_*}|
       \prod_{j\in S}\left(\frac{\sigma_{r(j)}}{2}|I_j|^{1/2}\right).
\end{align*}
This proves \eqref{weighted-laminar}.

We next prove \eqref{eq:detail-subset-sum}.
Fix an integer $n\ge1$. Recall that each $j\in\mathcal J_m$
has indices $(\ell(j),{\mu}(j),r(j))$ determined by
\eqref{eq:detail-index}, with $|I_j|=T2^{-\ell(j)}$.
For $0\le\ell_*\le J-1$ and $j_*\in\mathcal J_m$ with
$\ell(j_*)=\ell_*$, define
\[
 \mathfrak S_{\ell_*,j_*}
 :=\left\{S\subseteq\mathcal J_m:
 \begin{array}{c}
 |S|=n,\quad \displaystyle\bigcap_{j\in S}\operatorname{int}(I_j)
                    \ne\varnothing,\\
 \displaystyle\max_{j\in S}\ell(j)=\ell_*,\quad
 j_*=\min\{j\in S:\ell(j)=\ell_*\}
 \end{array}\right\}.
\]
Here, choosing the smallest index $j_*$ among those at level
$\ell_*$ makes the classes pairwise disjoint. Every $n$-element subset $S\subseteq\mathcal J_m$ with
$\bigcap_{j\in S}\operatorname{int}(I_j)\ne\varnothing$
belongs to exactly one such class; all other $n$-element
subsets satisfy $\eta_m(S)=0$.
If $J=0$, then $\mathcal J_m=\varnothing$ and the claimed sum is zero,
so below we may assume $J\ge1$.

For $S\in\mathfrak S_{\ell_*,j_*}$, the interval $I_{j_*}$ is shortest.
Raising \eqref{weighted-laminar} to the power $\rho$ and substituting
$|I_j|=T2^{-\ell(j)}$ gives
\begin{align}\label{eq:detail-set-level-bound}
 \eta_m(S)^\rho
 &\le (3/2)^\rho2^{-n\rho}
       \sigma_{r(j_*)}^\rho |I_{j_*}|^{3\rho/2}
       \prod_{j\in S\setminus\{j_*\}}
                \sigma_{r(j)}^\rho |I_j|^{\rho/2}\notag\\
 &= (3/2)^\rho2^{-n\rho}T^{\rho(n+2)/2}
       \sigma_{r(j_*)}^\rho 2^{-3\rho\ell_*/2}
       \prod_{j\in S\setminus\{j_*\}}
                \sigma_{r(j)}^\rho2^{-\rho\ell(j)/2}.
\end{align}
To bound the sum over $S\in\mathfrak S_{\ell_*,j_*}$ of the products
in \eqref{eq:detail-set-level-bound}, define the list
\[
 \mathcal A_{\ell_*}:=
 \bigl(\sigma_1^\rho,\ldots,\sigma_{d_w}^\rho,
  2^{-\rho/2}\sigma_1^\rho,\ldots,2^{-\rho/2}\sigma_{d_w}^\rho,
  \ldots,
  2^{-\rho\ell_*/2}\sigma_1^\rho,\ldots,
  2^{-\rho\ell_*/2}\sigma_{d_w}^\rho\bigr).
\]
Write $N_{\ell_*}=d_w(\ell_*+1)$ and denote the entries of
$\mathcal A_{\ell_*}$ by
$\alpha_q$, $1\le q\le N_{\ell_*}$. Define
\[
 e_{n-1}(\mathcal A_{\ell_*})
 :=\sum_{1\le q_1<\cdots<q_{n-1}\le N_{\ell_*}}
       \prod_{{i}=1}^{n-1}\alpha_{q_{{i}}},
 \qquad e_0(\mathcal A_{\ell_*})=1.
\]
The sum is zero when $n-1>N_{\ell_*}$.

For each $0\le\ell\le\ell_*$, let ${I^*_\ell}$ be the unique
dyadic interval of length $T2^{-\ell}$ containing $I_{j_*}$.
For $S\in\mathfrak S_{\ell_*,j_*}$, every $j\in S$ with
$\ell(j)=\ell$ satisfies $I_j={I^*_\ell}$.
By \eqref{eq:detail-index} and \eqref{eq:detail-support},
exactly $d_w$ indices $j\in\mathcal J_m$ have support ${I^*_\ell}$,
with pairwise distinct values $r(j)$ ranging over $\{1,\ldots,d_w\}$.
For each $j\in S\setminus\{j_*\}$, define
$q(j):=d_w\ell(j)+r(j)$. By the definition of
$\mathcal A_{\ell_*}$,
\[
 \alpha_{q(j)}=\sigma_{r(j)}^\rho2^{-\rho\ell(j)/2}.
\]
Each position $q$ determines $(\ell,r)$, which together with
${I^*_\ell}$ uniquely determines $j$ by \eqref{eq:detail-index} and
\eqref{eq:detail-support}. Thus the position set
$\{q(j):j\in S\setminus\{j_*\}\}$ uniquely determines
$S\in\mathfrak S_{\ell_*,j_*}$.
Summing over all $(n-1)$-element subsets of positions therefore
only enlarges the sum. Hence
\begin{equation}\label{eq:detail-class-weight-sum}
 \sum_{S\in\mathfrak S_{\ell_*,j_*}}
 \prod_{j\in S\setminus\{j_*\}}
       \sigma_{r(j)}^\rho2^{-\rho\ell(j)/2}
 \le e_{n-1}(\mathcal A_{\ell_*}).
\end{equation}
Note that
\begin{equation}\label{eq:deepest-index-weight-sum}
 \sum_{\substack{j_*\in\mathcal J_m\\\ell(j_*)=\ell_*}}
       \sigma_{r(j_*)}^\rho
 =\sum_{{\mu}=0}^{2^{\ell_*}-1}\sum_{r=1}^{d_w}\sigma_r^\rho
 =2^{\ell_*}\Theta_\rho.
\end{equation}
Using \eqref{eq:detail-set-level-bound},
\eqref{eq:detail-class-weight-sum}, and
\eqref{eq:deepest-index-weight-sum}, we obtain
\begin{align*}
 \sum_{\substack{S\subseteq\mathcal J_m\\|S|=n}}\eta_m(S)^\rho
 &=\sum_{\ell_*=0}^{J-1}
   \sum_{\substack{j_*\in\mathcal J_m\\\ell(j_*)=\ell_*}}
   \sum_{S\in\mathfrak S_{\ell_*,j_*}}\eta_m(S)^\rho\\
 &\le (3/2)^\rho2^{-n\rho}T^{\rho(n+2)/2}\Theta_\rho
       \sum_{\ell_*=0}^{J-1}2^{(1-3\rho/2)\ell_*}
       e_{n-1}(\mathcal A_{\ell_*}).
\end{align*}
Moreover,
\begin{equation}\label{eq:elementary-symmetric-bound}
 e_{n-1}(\mathcal A_{\ell_*})
 \le\frac{1}{(n-1)!}
       \left(\sum_{q=1}^{N_{\ell_*}}\alpha_q\right)^{n-1} =\frac{1}{(n-1)!}
       \left(\Theta_\rho\sum_{\ell=0}^{\ell_*}2^{-\rho\ell/2}\right)^{n-1}.
\end{equation}
Consequently,
\begin{align*}
 \sum_{\substack{S\subseteq\mathcal J_m\\|S|=n}}\eta_m(S)^\rho
 &\le\frac{(3/2)^\rho2^{-n\rho}T^{\rho(n+2)/2}\Theta_\rho}{(n-1)!}
       \sum_{\ell_*=0}^{J-1}2^{(1-3\rho/2)\ell_*}
       \left(\Theta_\rho\sum_{\ell=0}^{\ell_*}2^{-\rho\ell/2}\right)^{n-1}\\
 &\le\frac{(3/2)^\rho2^{-n\rho}T^{\rho(n+2)/2}\Theta_\rho}{(n-1)!}
       \left(\sum_{\ell_*=0}^{\infty}2^{(1-3\rho/2)\ell_*}\right)
       \left(\Theta_\rho\sum_{\ell=0}^{\infty}2^{-\rho\ell/2}\right)^{n-1}\\
 &=\frac{(3/2)^\rho2^{-n\rho}T^{\rho(n+2)/2}\Theta_\rho^n}{(1-2^{1-3\rho/2})(1-2^{-\rho/2})^{n-1}(n-1)!}.
\end{align*}
Both geometric series converge since $\rho>2/3$. To make the constant
explicit, set
\[
 {C'_{\rho,T}}=\frac{2^{-\rho}T^{\rho/2}}{1-2^{-\rho/2}},
 \qquad
 {C''_{\rho,T}}=\frac{(3/2)^\rho T^\rho(1-2^{-\rho/2})}{1-2^{1-3\rho/2}},
 \qquad
 C_{\rho,T}={C'_{\rho,T}}\max\{1,{C''_{\rho,T}}\}.
\]
The last bound is ${C''_{\rho,T}}{C'_{\rho,T}}^n \Theta_\rho^n/(n-1)!$. Since $n\ge1$, we have $ {C''_{\rho,T}}{C'_{\rho,T}}^n \Theta_\rho^n
\le C_{\rho,T}^n \Theta_\rho^n.$
The constant $C_{\rho,T}$ depends only on $\rho,T$. Taking the supremum over $m\in\mathcal M$
proves \eqref{eq:detail-subset-sum}.
\end{proof}

\begin{proof}[Proof of Lemma~\ref{lem:all-subset-sum}]
Fix $2/3<\rho<1/\tau$ and $\kappa>0$.
Similarly to \eqref{eq:elementary-symmetric-bound}, 
\begin{equation}\label{endpoint-series}
 \begin{aligned}
 {\Upsilon}^{\mathcal E}(\kappa)
=\sum_{\ell=0}^{d_w}(\kappa T^{\rho/2})^\ell(\ell!)^{\tau\rho}
   \sum_{\substack{E\subseteq\mathcal E_w\\|E|=\ell}}
   \prod_{j\in E}\sigma_j^\rho\le \sum_{\ell=0}^{\infty}
 (\kappa T^{\rho/2}\Theta_\rho)^\ell(\ell!)^{\tau\rho-1}=:{\Upsilon}_*^{\mathcal E}(\kappa).
 \end{aligned}
\end{equation}
For the detail inputs, Lemma~\ref{lem:detail-sum}  yields
\begin{align*}
 {\Upsilon}_m^{\mathcal J}(\kappa)
 &:=\sum_{\varnothing\ne D\subseteq\mathcal J_m}
 \kappa^{|D|}(|D|!)^{\tau\rho}\eta_m(D)^\rho=\sum_{n\ge1}\kappa^n(n!)^{\tau\rho}
   \sum_{\substack{D\subseteq\mathcal J_m\\|D|=n}}\eta_m(D)^\rho\\
 &\le 
 \sum_{n\ge1}(\kappa C_{\rho,T}\Theta_\rho)^n n(n!)^{\tau\rho-1}=:{\Upsilon}_*^{\mathcal J}(\kappa).
\end{align*}
If $\Theta_\rho>0$, the ratios of consecutive terms in ${\Upsilon}_*^{\mathcal E}(\kappa)$ and
${\Upsilon}_*^{\mathcal J}(\kappa)$ are
$\kappa T^{\rho/2}\Theta_\rho(\ell+1)^{\tau\rho-1}$ and
$\kappa C_{\rho,T}\Theta_\rho(1+1/n)(n+1)^{\tau\rho-1}$, respectively.
Both tend to zero since $\tau\rho<1$, so both series converge for every
$\kappa>0$. If $\Theta_\rho=0$, then ${\Upsilon}_*^{\mathcal E}(\kappa)=1$ and ${\Upsilon}_*^{\mathcal J}(\kappa)=0$.

Every nonempty $S\subseteq[s]$ has a unique decomposition
$S=D\mathbin{\dot\cup}E$, with $D\subseteq\mathcal J_m$ and
$E\subseteq\mathcal E_w$. The endpoint bound $|\bm\xi_j(t_i)|_2\le\sqrt T\,\sigma_j$
for $j\in E$, together with $mh=T$, gives
\begin{align*}
 \eta_m(D\cup E)
 =h\sum_{i=0}^{m-1}
   \left(\prod_{j\in D}|\bm\xi_j(t_i)|_2\right)
   \left(\prod_{j\in E}|\bm\xi_j(t_i)|_2\right)\le
 \begin{cases}
 \displaystyle\eta_m(D)\prod_{j\in E}(\sqrt T\,\sigma_j),
   &D\ne\varnothing,\\[2pt]
 \displaystyle T\prod_{j\in E}(\sqrt T\,\sigma_j),
   &D=\varnothing.
 \end{cases}
\end{align*}
The inequality $\binom{|D|+|E|}{|D|}\le2^{|D|+|E|}$ gives,
with $\kappa'=2^{\tau\rho}\kappa$,
\[
 \kappa^{|D|+|E|}((|D|+|E|)!)^{\tau\rho}
 \le (\kappa')^{|D|}(|D|!)^{\tau\rho}
     (\kappa')^{|E|}(|E|!)^{\tau\rho}.
\]
Summing over $D\ne\varnothing$ and all $E$ gives
\begin{align*}
 &\sum_{\varnothing\ne D\subseteq\mathcal J_m}
   \sum_{E\subseteq\mathcal E_w}
   \kappa^{|D|+|E|}((|D|+|E|)!)^{\tau\rho}\eta_m(D\cup E)^\rho\\
 &\quad\le
   \left[\sum_{\varnothing\ne D\subseteq\mathcal J_m}
   (\kappa')^{|D|}(|D|!)^{\tau\rho}\eta_m(D)^\rho\right]
   \left[\sum_{E\subseteq\mathcal E_w}
   (\kappa')^{|E|}(|E|!)^{\tau\rho}
   \prod_{j\in E}(\sqrt T\,\sigma_j)^\rho\right]\\
 &\quad={\Upsilon}_m^{\mathcal J}(\kappa'){\Upsilon}^{\mathcal E}(\kappa').
\end{align*}
For $D=\varnothing$, the condition $S\ne\varnothing$ requires
$E\ne\varnothing$, and
\begin{align*}
 \sum_{\varnothing\ne E\subseteq\mathcal E_w}
 \kappa^{|E|}(|E|!)^{\tau\rho}\eta_m(E)^\rho
 &\le T^\rho\sum_{\varnothing\ne E\subseteq\mathcal E_w}
 \kappa^{|E|}(|E|!)^{\tau\rho}
 \prod_{j\in E}(\sqrt T\,\sigma_j)^\rho\\
 &=T^\rho({\Upsilon}^{\mathcal E}(\kappa)-1).
\end{align*}
Combining these two contributions yields
\begin{align*}
 {\Upsilon}_m(\kappa)&\le {\Upsilon}_m^{\mathcal J}(\kappa'){\Upsilon}^{\mathcal E}(\kappa')
                  +T^\rho({\Upsilon}^{\mathcal E}(\kappa)-1)\\
 &\le {\Upsilon}_*^{\mathcal J}(\kappa'){\Upsilon}_*^{\mathcal E}(\kappa')
                  +T^\rho({\Upsilon}_*^{\mathcal E}(\kappa)-1)<\infty.
\end{align*}
Together with ${\Upsilon}^{\mathcal E}(\kappa)\le {\Upsilon}_*^{\mathcal E}(\kappa)$,
this proves the lemma with
\[
 C_{{\Upsilon}}(\kappa,\rho,\tau,T;\Theta_\rho):={\Upsilon}_*^{\mathcal E}(\kappa)
 +{\Upsilon}_*^{\mathcal J}(\kappa'){\Upsilon}_*^{\mathcal E}(\kappa')
 +T^\rho({\Upsilon}_*^{\mathcal E}(\kappa)-1).
\]
Every term in the defining series is nonnegative and nondecreasing in
$\kappa$ and $\Theta_\rho$. Hence this constant depends only on
$\kappa,\rho,\tau,T,\Theta_\rho$ and is nondecreasing in both
$\kappa$ and $\Theta_\rho$.
\end{proof}

\section{Proofs of the weighted Gaussian norm estimates}\label{sec:auxiliary-moment-proof}
\begin{proof}[Proof of Lemma~\ref{lem:hybrid}]
For each $1\le r\le d_w$ and $0\le\ell<J$, define
\[
 M_{\ell,r}:=
 \max_{\substack{j\in\mathcal J_m\\\ell(j)=\ell,\ r(j)=r}}|z_j|.
\]
For fixed $0\le\ell<J$, $1\le r\le d_w$, and $0\le i\le m$,
\eqref{eq:detail-support} and
\eqref{eq:diffusion-coefficients} give
\[
 \left|\{j\in\mathcal J_m:\ell(j)=\ell,\ r(j)=r,
                  \bm\xi_j(t_i)\ne\bm0\}\right|\le1.
\]
Hence, by \eqref{eq:linearbridge} and \eqref{eq:bridgebounds},
\begin{align}
 |\Sigma\bm W_{t_i}(\bm z)|_2
 &\le\sum_{r=1}^{d_w}|\bm\xi_r(t_i)|_2|z_r|
 +\sum_{\ell=0}^{J-1}\sum_{r=1}^{d_w}
 \sum_{\substack{j\in\mathcal J_m\\\ell(j)=\ell,\ r(j)=r}}
 |\bm\xi_j(t_i)|_2|z_j|\notag\\
 &\le\sqrt T\sum_{r=1}^{d_w}\sigma_r|z_r|
 +\frac{\sqrt T}{2}\sum_{\ell=0}^{J-1}\sum_{r=1}^{d_w}
 \sigma_r2^{-\ell/2}M_{\ell,r}\notag\\
 &=\sqrt T\sum_{r=1}^{d_w} \sigma_r
       \left(|z_r|+\frac{1}{2}\sum_{\ell=0}^{J-1}
                              2^{-\ell/2}M_{\ell,r}\right).
 \label{eq:auxiliary-noise-path-bound}
\end{align}
The same bound holds for $\max_{0\le i\le m}|\Sigma\bm W_{t_i}(\bm z)|_2$.
The telescoped Euler recursion and the drift bound in \eqref{eq:growth} give
\[
 |\bm{X}^h_{t_i}(\bm{z})|_2\le |\bm{x}_0|_2+K_0T+
 \max_{0\le {i'}\le m}|\Sigma \bm{W}_{t_{{i'}}}(\bm{z})|_2
 +\mathsf B_1h\sum_{{i'}=0}^{i-1}|\bm{X}^h_{t_{{i'}}}(\bm{z})|_2.
\]
Discrete Gronwall yields
\begin{equation}\label{eq:euler-state-gronwall}
 \max_i|\bm{X}^h_{t_i}(\bm{z})|_2\le e^{\mathsf B_1T}
 \left(|\bm{x}_0|_2+K_0T+\max_i|\Sigma \bm{W}_{t_i}(\bm{z})|_2\right).
\end{equation}

In the following moment estimates, evaluate these quantities at
$\bm z={\bm Z^{(u)}}$. Set
\[
 \lambda_0:=2\beta_g e^{\mathsf B_1T}\sqrt T\,\Theta_1
 \left(2+\frac{1}{\sqrt2}\right).
\]
Then \eqref{eq:alpha-choice} gives $2\alpha_0-\lambda_0=2$, and each
coordinate satisfies
\[
 \E e^{\lambda_0|{Z_j^{(u)}}|}\le K_{\alpha_0}(\lambda_0)
       :=\max\left(2e^{\lambda_0^2/2},
                    \frac{2\alpha_0}{2\alpha_0-\lambda_0}\right)
       =\max\left(2e^{\lambda_0^2/2},\alpha_0\right).
\]
Since \eqref{eq:detail-index} gives $2^\ell$ indices with
$\ell(j)=\ell$ and $r(j)=r$, we have
\[
 \E e^{\lambda_0 M_{\ell,r}}
 \le\sum_{\substack{j\in\mathcal J_m\\\ell(j)=\ell,\ r(j)=r}}
       \E e^{\lambda_0|{Z_j^{(u)}}|}
 \le2^\ell K_{\alpha_0}(\lambda_0).
\]
Assume $\Theta_1>0$, with $\Theta_1$ defined in \eqref{eq:column-sums}.
Define the weights
\[
 w_r=\frac{\sigma_r}{\Theta_1\left(2+\frac{1}{\sqrt2}\right)},\qquad
 w_{\ell,r}=\frac{\sigma_r2^{-\ell/2}}{2\Theta_1\left(2+\frac{1}{\sqrt2}\right)}.
\]
Since $1+\frac12\sum_{\ell\ge0}2^{-\ell/2}=2+1/\sqrt2$, we have
\[
 \omega:=\sum_{r=1}^{d_w}w_r+
       \sum_{\ell=0}^{J-1}\sum_{r=1}^{d_w}w_{\ell,r}
 =\frac{1+\frac12\sum_{\ell=0}^{J-1}2^{-\ell/2}}
        {2+1/\sqrt2}\le1.
\]
Put $w_*=1-\omega$. Then
$\lambda_0 w_r=2\beta_g e^{\mathsf B_1T}\sqrt T\,\sigma_r$ and
$\lambda_0 w_{\ell,r}=2\beta_g e^{\mathsf B_1T}\sqrt T\,\sigma_r2^{-\ell/2}/2$.
Omitting zero weights, generalized H\"older gives
\begin{align*}
 &\E\exp\left\{2\beta_g e^{\mathsf B_1T}\sqrt T\sum_{r=1}^{d_w}\sigma_r
       \left(|{Z_r^{(u)}}|+\frac12\sum_{\ell=0}^{J-1}
                   2^{-\ell/2}M_{\ell,r}\right)\right\}\\
 &=\E\left[
    \prod_{r=1}^{d_w}\left(e^{\lambda_0|{Z_r^{(u)}}|}\right)^{w_r}
    \prod_{\ell=0}^{J-1}\prod_{r=1}^{d_w}
         \left(e^{\lambda_0 M_{\ell,r}}\right)^{w_{\ell,r}}
    1^{w_*}\right]\\
 &\le\prod_{r=1}^{d_w}
       \left(\E e^{\lambda_0|{Z_r^{(u)}}|}\right)^{w_r}
       \prod_{\ell=0}^{J-1}\prod_{r=1}^{d_w}
       \left(\E e^{\lambda_0 M_{\ell,r}}\right)^{w_{\ell,r}}\\
 &\le K_{\alpha_0}(\lambda_0)^\omega
       2^{\sum_{\ell=0}^{J-1}\sum_{r=1}^{d_w}\ell w_{\ell,r}}\\
 &\le K_{\alpha_0}(\lambda_0)
       \exp\left\{\frac{\log2}{4+\sqrt2}
                    \sum_{\ell=0}^{J-1}\ell2^{-\ell/2}\right\}\\
 &\le K_{\alpha_0}(\lambda_0)
       \exp\left\{\frac{\log2}{4+\sqrt2}
                    \sum_{\ell=0}^{\infty}\ell2^{-\ell/2}\right\}
 =K_{\alpha_0}(\lambda_0)2^{\frac{4+3\sqrt2}{4+\sqrt2}}.
\end{align*}
Here we used $\omega\le1$, $K_{\alpha_0}(\lambda_0)\ge1$,
$\sum_{r=1}^{d_w}\sigma_r=\Theta_1$, and
$\sum_{\ell\ge0}\ell2^{-\ell/2}=4+3\sqrt2$. Define
\begin{equation}\label{eq:auxiliary-moment-constant}
 C_{\rm aux}(|\bm x_0|_2,\Theta_1)
 :=\exp\left\{2\beta_g e^{\mathsf B_1T}
       (|\bm x_0|_2+K_0T)\right\}K_{\alpha_0}(\lambda_0)
              2^{\frac{4+3\sqrt2}{4+\sqrt2}}<\infty.
\end{equation}
Combining \eqref{eq:auxiliary-noise-path-bound},
\eqref{eq:euler-state-gronwall}, and the preceding estimate gives
\[
 \sup_{m\in\mathcal M}\sup_{u\subseteq[s]}
 \E e^{2\beta_g|\bm X^h_{t_m}({\bm Z^{(u)}})|_2}
 \le C_{\rm aux}(|\bm x_0|_2,\Theta_1).
\]
If $\Theta_1=0$, then $\Sigma=0$, and
\eqref{eq:euler-state-gronwall} directly gives
$\E e^{2\beta_g|\bm X^h_{t_m}({\bm Z^{(u)}})|_2}
\le\exp\left\{2\beta_g e^{\mathsf B_1T}
       (|\bm x_0|_2+K_0T)\right\}$, so the same bound follows from
\eqref{eq:auxiliary-moment-constant} with $\lambda_0=0$.
Both entries of the maximum defining $K_{\alpha_0}(\lambda_0)$ are
nondecreasing in $\Theta_1$. Hence \eqref{eq:auxiliary-moment-constant}
is nondecreasing in $|\bm x_0|_2$ and $\Theta_1$.
\end{proof}

\begin{proof}[Proof of Lemma~\ref{lem:transfer}]
Let ${\bm Z^{(u)}}$ be as defined in Lemma~\ref{lem:hybrid}.
By \eqref{eq:normalized-derivative-envelope},
$|\partial_uF_m(\bm{z})|\le\mathcal D_{m,u}e^{\beta_g|\bm{X}^h_{t_m}(\bm{z})|_2}$.
Jensen's inequality and $\alpha_0\ge1$ give
\begin{align*}
 {\mathcal S}_u(F_m)^2
 &\le\int_{\R^s}|\partial_uF_m(\bm{z})|^2\phi_{-u}(\bm{z}_{-u})
                         \prod_{j\in u}e^{-2\alpha_0|z_j|}\dd \bm{z}\\
 &\le\alpha_0^{-|u|}\mathcal D_{m,u}^2
       \E e^{2\beta_g|\bm{X}^h_{t_m}({\bm Z^{(u)}})|_2}\le\mathcal D_{m,u}^2
       \E e^{2\beta_g|\bm{X}^h_{t_m}({\bm Z^{(u)}})|_2}\\
&\le C_{\rm aux}(|\bm x_0|_2,\Theta_1)\mathcal D_{m,u}^2.
\end{align*}
The last inequality uses Lemma~\ref{lem:hybrid}.
Taking square roots proves \eqref{eq:Bbound}.
\end{proof}

\section{Moment bounds and the Euler bias}\label{sec:moments}
Here $\bm{X}_t^{\bm{x}}$ and $\widetilde{\bm{X}}_{t_i}^{\bm{x}}$ denote,
respectively, the solution of the SDE~\eqref{eq:sde} and its
Euler--Maruyama approximation \eqref{eq:em}, both with initial
state $\bm{x}$. Throughout, $h=T/m$.
\begin{lemma}\label{lem:noise-moments}
For every $v>0$,
\begin{equation}\label{eq:noise-explicit-moment}
 \E\exp\left(v\sup_{t\le T}|\Sigma \bm{W}_t|_2\right)
 \le 5\exp\left(v^2T(1+\|\Sigma\|_{\rm HS}^2)\right).
\end{equation}
In particular, this bound is {uniform in $d_x,d_w$ under Assumption~\ref{ass:columns}, since then $\|\Sigma\|_{\rm HS}^2\le M_{2/3}^3$ by \eqref{eq:column-monotonicity}}.
\end{lemma}
\begin{proof}
The case $\Sigma=0$ is immediate. By orthogonal diagonalization of
$\Sigma\Sigma^\top$, the process $\Sigma \bm{W}$ has the same law as
an orthogonal image of $({\varsigma_r\widetilde W_r})_{r=1}^{{\operatorname{rank}(\Sigma)}}$, where
{$\varsigma_r>0$ are the nonzero singular values of $\Sigma$}, $\sum_{r=1}^{{\operatorname{rank}(\Sigma)}}{\varsigma_r^2}=\|\Sigma\|_{\rm HS}^2$, and the ${\widetilde W_r}$ are independent
one-dimensional Brownian motions. This identity holds for the
whole process, since the centered Gaussian covariances agree.
Thus $\sup_{t\le T}|\Sigma \bm{W}_t|_2^2$ has the same distribution as
$\sup_{t\le T}\sum_{r=1}^{{\operatorname{rank}(\Sigma)}}{\varsigma_r^2} {\widetilde W_r}(t)^2$.
Writing ${\overline W_r}=\sup_{t\le T}|{\widetilde W_r}(t)|$, we have
\[
 \sup_{t\le T}\sum_{r=1}^{{\operatorname{rank}(\Sigma)}}{\varsigma_r^2} {\widetilde W_r}(t)^2
 \le \sum_{r=1}^{{\operatorname{rank}(\Sigma)}}{\varsigma_r^2} {\overline W_r}^2.
\]
Convexity of the exponential, with weights
${\varsigma_r^2}/\|\Sigma\|_{\rm HS}^2$, gives
\[
 \E\exp\left(c\sum_{r=1}^{{\operatorname{rank}(\Sigma)}}{\varsigma_r^2\overline W_r^2}\right)
 \le\sum_{r=1}^{{\operatorname{rank}(\Sigma)}}\frac{{\varsigma_r^2}}{\|\Sigma\|_{\rm HS}^2}
       \E e^{c\|\Sigma\|_{\rm HS}^2{\overline W_r}^2}
 =\E e^{c\|\Sigma\|_{\rm HS}^2{\overline W_1}^2}.
\]
The reflection principle \cite[Theorem~2.21]{morters2010brownian} gives
$\Prob({\overline W_1}>u)\le4e^{-u^2/(2T)}$. Hence
\[
 \E e^{a{\overline W_1}^2}
 \le1+8a\int_0^\infty u e^{-(1/(2T)-a)u^2}\,\mathrm du
 <\infty,\qquad 0<a<1/(2T).
\]
For $0<a\le1/(4T)$, the right-hand side is at most
$1+4a/(1/(2T)-a)\le5$.
Take $c=[4T(1+\|\Sigma\|_{\rm HS}^2)]^{-1}$.
Then $c\|\Sigma\|_{\rm HS}^2\le1/(4T)$, so the preceding display gives
\[
 \E\exp\left(c\sup_{t\le T}|\Sigma\bm W_t|_2^2\right)\le5.
\]
Using $vy\le cy^2+v^2/(4c)$ proves
\eqref{eq:noise-explicit-moment}.
All these inequalities use scalar constants independent of the rank.
\end{proof}

\begin{lemma}\label{lem:moments}
Under Assumption~\ref{ass:space}, \eqref{eq:sde} has a unique
global strong solution. Define
\begin{equation}\label{eq:state-moment-constant}
 C_X\bigl(|\bm{x}|_2,\|\Sigma\|_{\rm HS}\bigr)
 :=e^{\mathsf B_1T}
       \left(|\bm{x}|_2+K_0T+2\sqrt T\,\|\Sigma\|_{\rm HS}\right).
\end{equation}
Then
\begin{equation}\label{eq:moments}
 \left(\E\sup_{t\le T}|\bm{X}_t^{\bm{x}}|_2^2\right)^{1/2}
 \le C_X\bigl(|\bm{x}|_2,\|\Sigma\|_{\rm HS}\bigr).
\end{equation}
For $v>0$, define
\begin{equation}\label{eq:state-exp-moment-constant}
 C_{\rm exp}\bigl(v,\|\Sigma\|_{\rm HS}\bigr):=10\exp\left(
 v e^{\mathsf B_1T}K_0T
 +v^2e^{2\mathsf B_1T}T(1+\|\Sigma\|_{\rm HS}^2)\right).
\end{equation}
For every $m\in\mathcal M$,
\begin{equation}\label{eq:state-exp-moments}
 \E e^{v\sup_{t\le T}|\bm{X}_t^{\bm{x}}|_2}
 +\E e^{v\max_{i\le m}|\widetilde{\bm{X}}_{t_i}^{\bm{x}}|_2}
 \le C_{\rm exp}\bigl(v,\|\Sigma\|_{\rm HS}\bigr)
       e^{v e^{\mathsf B_1T}|\bm{x}|_2}.
\end{equation}
\end{lemma}
\begin{proof}
Existence and uniqueness follow from global spatial Lipschitz
continuity, linear growth, and Borel time dependence; see
\cite[Theorem~7.1]{pages2018numerical}.
Put
$A=|\bm{x}|_2+K_0T+\sup_{t\le T}|\Sigma \bm{W}_t|_2$.
The integral equation, the telescoped Euler recursion, and the
continuous and discrete Gronwall inequalities give
\[
 \max\left\{\sup_{t\le T}|\bm{X}_t^{\bm{x}}|_2,
             \max_{i\le m}|\widetilde{\bm{X}}_{t_i}^{\bm{x}}|_2\right\}
 \le e^{\mathsf B_1T}A
\]
on the same Brownian path. Doob's $L^2$ inequality gives
$\|\sup_{t\le T}|\Sigma\bm W_t|_2\|_{L^2}
\le2\sqrt T\|\Sigma\|_{\rm HS}$, which proves \eqref{eq:moments}.
The left side of \eqref{eq:state-exp-moments} is bounded by
\[
 2e^{v e^{\mathsf B_1T}(|\bm{x}|_2+K_0T)}
 \E e^{v e^{\mathsf B_1T}\sup_{t\le T}|\Sigma \bm{W}_t|_2}.
\]
Applying \eqref{eq:noise-explicit-moment} with
$v e^{\mathsf B_1T}$ proves \eqref{eq:state-exp-moments} with the
choice \eqref{eq:state-exp-moment-constant}.
\end{proof}

For the following two proofs, couple the solution $\bm X$ of
\eqref{eq:sde} and its Euler--Maruyama approximation $\bm X^h$ in
\eqref{eq:bb-em} using the same Brownian motion $\bm W$.
Define
\begin{equation}\label{eq:payoff-transfer-constant}
 C_{\rm pay}\bigl(|\bm x_0|_2,\|\Sigma\|_{\rm HS}\bigr)
 :=G_0G_1\sqrt{2C_{\rm exp}(2\beta_g,\|\Sigma\|_{\rm HS})}
       \exp\left(\beta_g e^{\mathsf B_1T}|\bm x_0|_2\right).
\end{equation}
The first derivative bound in \eqref{eq:payoff-ass} gives
\[
 |g(\bm y)-g(\bm x)|
 \le G_0G_1|\bm y-\bm x|_2
       \bigl(e^{\beta_g|\bm x|_2}+e^{\beta_g|\bm y|_2}\bigr).
\]
Cauchy--Schwarz and \eqref{eq:state-exp-moments} imply
\begin{equation}\label{eq:bias-from-strong}
 \begin{aligned}
 |\mathcal I_m-\mathcal I|
 &\le G_0G_1\|\bm X_T-\bm X^h_{t_m}\|_{L^2}
 \left[
  \bigl(\E e^{2\beta_g|\bm X_T|_2}\bigr)^{1/2}
  +\bigl(\E e^{2\beta_g|\bm X^h_{t_m}|_2}\bigr)^{1/2}
 \right]\\
 &\le C_{\rm pay}\bigl(|\bm x_0|_2,\|\Sigma\|_{\rm HS}\bigr)
       \|\bm X_T-\bm X^h_{t_m}\|_{L^2}.
 \end{aligned}
\end{equation}

\begin{proof}[Proof of Theorem~\ref{thm:euler-strong-half}]
{Define}
\begin{equation}\label{eq:euler-increment-constant}
 C_{\rm inc}\bigl(|\bm x_0|_2,\|\Sigma\|_{\rm HS}\bigr)
 :=\sqrt T\left[K_0+\mathsf B_1
 C_X\bigl(|\bm x_0|_2,\|\Sigma\|_{\rm HS}\bigr)\right]
 +\|\Sigma\|_{\rm HS}{.}
\end{equation}
{For $s\in[t_i,t_{i+1}]$, the integral equation, It\^o's isometry, \eqref{eq:moments}, and $s-t_i\le T$ give}
\begin{align*}
 \|\bm{X}_s-\bm{X}_{t_i}\|_{L^2}
 &\le\int_{t_i}^s
       \bigl(K_0+\mathsf B_1\|\bm X_{{t'}}\|_{L^2}\bigr)\dd {{t'}}
       +\|\Sigma\|_{\rm HS}(s-t_i)^{1/2}\\
 &\le C_{\rm inc}\bigl(|\bm x_0|_2,\|\Sigma\|_{\rm HS}\bigr)
       (s-t_i)^{1/2}.
\end{align*}
Split $b(s,\bm{X}_s)-b(t_i,\bm{X}_{t_i})$ into
$b(s,\bm{X}_s)-b(s,\bm{X}_{t_i})$ and
$b(s,\bm{X}_{t_i})-b(t_i,\bm{X}_{t_i})$. The spatial Lipschitz
bound and \eqref{eq:strong-time-regularity} then yield
\begin{align*}
 \|b(s,\bm{X}_s)-b(t_i,\bm{X}_{t_i})\|_{L^2}
 &\le\mathsf B_1\|\bm X_s-\bm X_{t_i}\|_{L^2}
       +{\widetilde L_t}\bigl(1+\|\bm X_{t_i}\|_{L^2}\bigr)(s-t_i)^{1/2}\\
 &\le\left[\mathsf B_1
       C_{\rm inc}\bigl(|\bm x_0|_2,\|\Sigma\|_{\rm HS}\bigr)
       +{\widetilde L_t}\left(1+C_X\bigl(|\bm x_0|_2,\|\Sigma\|_{\rm HS}\bigr)\right)\right]
       (s-t_i)^{1/2}.
\end{align*}
Let $\bm{E}_i=\bm{X}_{t_i}-\bm{X}^h_{t_i}$. Subtracting the
Euler recursion from the integral equation gives
\[
 \bm{E}_j
 =h\sum_{i<j}\bigl[b(t_i,\bm{X}_{t_i})
                  -b(t_i,\bm{X}^h_{t_i})\bigr]
 +\sum_{i<j}\int_{t_i}^{t_{i+1}}
        \bigl[b(s,\bm{X}_s)-b(t_i,\bm{X}_{t_i})\bigr]\dd s.
\]
Discrete Gronwall applied pathwise implies
\[
 \max_{j\le m}|\bm{E}_j|_2
 \le e^{\mathsf B_1T}\sum_{i=0}^{m-1}\int_{t_i}^{t_{i+1}}
        |b(s,\bm{X}_s)-b(t_i,\bm{X}_{t_i})|_2\dd s.
\]
Define
\begin{align}\label{eq:euler-half-constants}
 C_{{\rm str},1/2}\bigl(|\bm x_0|_2,\|\Sigma\|_{\rm HS}\bigr)
 &:=\frac23 e^{\mathsf B_1T}T\sqrt T
 \left[\mathsf B_1C_{\rm inc}\bigl(|\bm x_0|_2,\|\Sigma\|_{\rm HS}\bigr)
       +{\widetilde L_t}\left(1+C_X\bigl(|\bm x_0|_2,\|\Sigma\|_{\rm HS}\bigr)\right)\right],\notag\\
 C_{{\rm bias},1/2}\bigl(|\bm x_0|_2,\|\Sigma\|_{\rm HS}\bigr)
 &:=C_{\rm pay}\bigl(|\bm x_0|_2,\|\Sigma\|_{\rm HS}\bigr)
 C_{{\rm str},1/2}\bigl(|\bm x_0|_2,\|\Sigma\|_{\rm HS}\bigr).
\end{align}
Taking $L^2$ norms and using Minkowski's inequality gives
\begin{align*}
 \left\|\max_{j\le m}|\bm{E}_j|_2\right\|_{L^2}
 &\le e^{\mathsf B_1T}
 \left[\mathsf B_1C_{\rm inc}\bigl(|\bm x_0|_2,\|\Sigma\|_{\rm HS}\bigr)
       +{\widetilde L_t}\left(1+C_X\bigl(|\bm x_0|_2,\|\Sigma\|_{\rm HS}\bigr)\right)\right]\\
 &\quad\times
 \sum_{i=0}^{m-1}\int_{t_i}^{t_{i+1}}(s-t_i)^{1/2}\dd s\\
 &=C_{{\rm str},1/2}\bigl(|\bm x_0|_2,\|\Sigma\|_{\rm HS}\bigr)m^{-1/2}.
\end{align*}
This proves \eqref{eq:strong-half-main}. Applying
\eqref{eq:bias-from-strong} gives
$|\mathcal I_m-\mathcal I|
 \le C_{{\rm bias},1/2}\bigl(|\bm x_0|_2,\|\Sigma\|_{\rm HS}\bigr)m^{-1/2}$ and proves
\eqref{eq:bias-half-main}.
Under the additional hypotheses in (ii), the model constants,
$|\bm x_0|_2$ and ${\widetilde L_t}$ have common bounds, and
\eqref{eq:column-monotonicity} gives
$\|\Sigma\|_{\rm HS}^2\le M_{2/3}^3$, hence
$\|\Sigma\|_{\rm HS}\le M_{2/3}^{3/2}$.
The constants in \eqref{eq:state-moment-constant},
\eqref{eq:state-exp-moment-constant},
\eqref{eq:payoff-transfer-constant},
\eqref{eq:euler-increment-constant}, and
\eqref{eq:euler-half-constants} are nondecreasing in
$|\bm x_0|_2$ and $\|\Sigma\|_{\rm HS}$.
Thus the dimension-uniform constants are obtained by evaluating
\eqref{eq:euler-half-constants} at
$(|\bm x_0|_2,\|\Sigma\|_{\rm HS})=(R_0,M_{2/3}^{3/2})$.
\end{proof}

\begin{proof}[Proof of Theorem~\ref{thm:euler-bias}]
Let $\bm{E}_i=\bm{X}_{t_i}-\bm{X}^h_{t_i}$.
The additive noise cancels
in the error recursion, giving
\begin{align}
 \bm{E}_{i+1}&=\bm{E}_i
 +h\bigl[b(t_i,\bm{X}_{t_i})-b(t_i,\bm{X}^h_{t_i})\bigr]
 +\bm{R}_i,\label{eq:bias-error-recursion}\\
 \bm{R}_i&=\int_{t_i}^{t_{i+1}}
   \bigl[b(s,\bm{X}_s)-b(t_i,\bm{X}_{t_i})\bigr]\dd s.\notag
\end{align}
Split the integrand in $\bm{R}_i$ into
$b(s,\bm{X}_s)-b(t_i,\bm{X}_s)$ and
$b(t_i,\bm{X}_s)-b(t_i,\bm{X}_{t_i})$.
Applying It\^o's formula to the fixed function $b(t_i,\cdot)$ and
then stochastic Fubini gives $\bm{R}_i=\bm{A}_i+\bm{M}_i$, where
\begin{align*}
 \bm{A}_i
 &=\int_{t_i}^{t_{i+1}}
    \bigl[b(s,\bm{X}_s)-b(t_i,\bm{X}_s)\bigr]\dd s\\
 &\quad+\int_{t_i}^{t_{i+1}}(t_{i+1}-{{t'}})
 \Bigl(D_xb(t_i,\bm{X}_{{t'}})b({{t'}},\bm{X}_{{t'}})
 +\frac12\sum_{{{r}}=1}^{d_w}
 D_x^2b(t_i,\bm{X}_{{t'}})[\Sigma\bm{e}_{{r}},\Sigma\bm{e}_{{r}}]\Bigr)\dd {{t'}},\\
 \bm{M}_i
 &=\int_{t_i}^{t_{i+1}}(t_{i+1}-{{t'}})
       D_xb(t_i,\bm{X}_{{t'}})\Sigma\dd\bm{W}_{{t'}}.
\end{align*}
The bounded first and second spatial derivatives and
Lemma~\ref{lem:moments} justify these formulas by localization;
the stochastic integrand is square integrable. The operator-norm
bounds give
\begin{align}
 \left|\sum_{{{r}}=1}^{d_w}
 D_x^2b(t_i,\bm{X}_{{t'}})[\Sigma\bm{e}_{{r}},\Sigma\bm{e}_{{r}}]\right|_2
 &\le\mathsf B_2\|\Sigma\|_{\rm HS}^2,\label{eq:bias-trace-bound}\\
 \|D_xb(t_i,\bm{X}_{{t'}})\Sigma\|_{\mathrm{HS}}^2
 &\le\mathsf B_1^2\|\Sigma\|_{\rm HS}^2.\label{eq:bias-integrand-bound}
\end{align}
Thus \eqref{eq:bias-time-regularity}, linear growth of $b$,
Minkowski's inequality, and \eqref{eq:moments} imply the bound below.
Define
\begin{equation}\label{eq:euler-drift-remainder-constant}
 C_A\bigl(|\bm x_0|_2,\|\Sigma\|_{\rm HS}\bigr)
 :=\frac12\left(L_t+\mathsf B_1K_0
                 +\frac12\mathsf B_2\|\Sigma\|_{\rm HS}^2\right)
   +\frac12(L_t+\mathsf B_1^2)
    C_X\bigl(|\bm x_0|_2,\|\Sigma\|_{\rm HS}\bigr).
\end{equation}
Then
\begin{equation}\label{eq:bias-drift-remainder}
 \begin{aligned}
 \|\bm{A}_i\|_{L^2}
 &\le\frac{h^2}{2}
 \left(L_t+\mathsf B_1K_0
       +\frac12\mathsf B_2\|\Sigma\|_{\rm HS}^2\right)+\frac{h^2}{2}(L_t+\mathsf B_1^2)
       \left(\E\sup_{t\le T}|\bm X_t|_2^2\right)^{1/2}\\
 &\le C_A\bigl(|\bm x_0|_2,\|\Sigma\|_{\rm HS}\bigr)h^2.
 \end{aligned}
\end{equation}

The partial sums of $\bm{M}_i$ form a square-integrable martingale
at the grid points. It\^o's isometry, orthogonality of its increments,
and \eqref{eq:bias-integrand-bound} give
\[
 \E\left|\sum_{i=0}^{m-1}\bm{M}_i\right|_2^2
 \le\mathsf B_1^2\|\Sigma\|_{\rm HS}^2\sum_{i=0}^{m-1}
       \int_{t_i}^{t_{i+1}}(t_{i+1}-{{t'}})^2\dd {{t'}}
 =\frac{\mathsf B_1^2\|\Sigma\|_{\rm HS}^2T}{3}h^2.
\]
Doob's inequality therefore yields
\begin{equation}\label{eq:bias-martingale-remainder}
 \left\|\max_{j\le m}
       \left|\sum_{i<j}\bm{M}_i\right|_2\right\|_{L^2}
 \le 2\mathsf B_1\|\Sigma\|_{\rm HS}\sqrt{\frac{T}{3}}\,h.
\end{equation}
Summing \eqref{eq:bias-error-recursion} and using the spatial
Lipschitz bound, discrete Gronwall gives
\[
 \max_{j\le m}|\bm{E}_j|_2
 \le e^{\mathsf B_1T}\max_{j\le m}
       \left|\sum_{i<j}(\bm{A}_i+\bm{M}_i)\right|_2.
\]
Since the maximum of the partial sums of $\bm{A}_i$ is bounded
by $\sum_i|\bm{A}_i|_2$, estimates
\eqref{eq:bias-drift-remainder} and
\eqref{eq:bias-martingale-remainder} give the following choices:
\begin{align}\label{eq:euler-order-one-constants}
 C_{\rm str}\bigl(|\bm x_0|_2,\|\Sigma\|_{\rm HS}\bigr)
 &:=e^{\mathsf B_1T}
   \left[TC_A\bigl(|\bm x_0|_2,\|\Sigma\|_{\rm HS}\bigr)
   +2\mathsf B_1\|\Sigma\|_{\rm HS}\sqrt{\frac T3}\right],\notag\\
 C_{\rm bias}\bigl(|\bm x_0|_2,\|\Sigma\|_{\rm HS}\bigr)
 &:=T C_{\rm pay}\bigl(|\bm x_0|_2,\|\Sigma\|_{\rm HS}\bigr)
 C_{\rm str}\bigl(|\bm x_0|_2,\|\Sigma\|_{\rm HS}\bigr).
\end{align}
Indeed,
\begin{equation}\label{eq:bias-strong-order-one}
 \left\|\max_{j\le m}|\bm{E}_j|_2\right\|_{L^2}
 \le C_{\rm str}\bigl(|\bm x_0|_2,\|\Sigma\|_{\rm HS}\bigr)h.
\end{equation}
Combining \eqref{eq:bias-strong-order-one} with
\eqref{eq:bias-from-strong} and $h=T/m$ gives
\[
 |\mathcal I_m-\mathcal I|
 \le C_{\rm bias}\bigl(|\bm x_0|_2,\|\Sigma\|_{\rm HS}\bigr)m^{-1},
\]
which proves (i).

For (ii), Assumption~\ref{ass:uniform-data} bounds the model
constants and $|\bm x_0|_2$ uniformly in $d_x,d_w$, while
Assumption~\ref{ass:columns} and \eqref{eq:column-monotonicity} give
$\|\Sigma\|_{\rm HS}\le M_{2/3}^{3/2}$. For fixed common model bounds and $L_t$,
the constants in \eqref{eq:euler-drift-remainder-constant} and
\eqref{eq:euler-order-one-constants} are nondecreasing in
$|\bm x_0|_2$ and $\|\Sigma\|_{\rm HS}$.
Thus the dimension-uniform constants are obtained by evaluating
\eqref{eq:euler-order-one-constants} at
$(|\bm x_0|_2,\|\Sigma\|_{\rm HS})=(R_0,M_{2/3}^{3/2})$.
\end{proof}

\paragraph*{Acknowledgement} 

The work of Josef Dick was supported in part by a grant of access to OpenAI models through the ChatGPT for Academic Researchers program. The work of Ziyang Ye was supported by the China Scholarship Council
program (Project ID: 202606210192).

\paragraph*{Declaration of generative AI use}

The authors used ChatGPT 6 Astra for literature searches, exploratory development, and the preparation of portions of the exposition and LaTeX source. All mathematical arguments, calculations, references, and conclusions were independently checked and verified by the authors, who take full responsibility for the contents of the paper.

\paragraph*{Lean 4 formalization.}
A Lean 4 formalization of Theorem~\ref{thm:main},
Corollary~\ref{cor:dimensions}, 
Theorem~\ref{thm:euler-strong-half}, and Theorem~\ref{thm:euler-bias}
is available at
\url{https://github.com/ZiyangYe2000/Lattice4SDE-Lean}.
The formalization is carried out under
Assumptions~\ref{ass:space}, \ref{ass:uniform-data}, and
\ref{ass:columns}. It uses Lemmas~\ref{lem:random-wce} and
\ref{lem:moments}, as well as standard stochastic and distributional
identities, as external inputs.
A complete account of these inputs is provided in the repository
blueprint. The project can be checked by running \texttt{lake build}.

\bibliographystyle{siam}
\bibliography{ref}

\end{document}